\input amstex
\documentstyle{amsppt}
\pageheight{194mm}
\pagewidth{133mm}
\magnification\magstep1

\document

%%%%% \amsusui.tex %%%%%

\def\nologo{\let\logo@\empty}

\def\Ad{\operatorname{Ad}}

\def\Aut{\operatorname{Aut}}
\def\Aut{\operatorname{Aut}}

\def\BS{\operatorname{BS}}

\def\End{\operatorname{End}}

\def\face{\operatorname{face}}

\def\fil{\operatorname{fil}}

\def\gr{\operatorname{gr}}
\def\Hom{\operatorname{Hom}}

\def\Im{\operatorname{Im}}
\def\Image{\operatorname{Image}}
\def\Int{\operatorname{Int}}
\def\Isom{\operatorname{Isom}}
\def\Ker{\operatorname{Ker}}

\def\Map{\;\operatorname{Map}}

\def\Mor{\operatorname{Mor}}

\def\rank{\operatorname{rank}}
\def\Re{\operatorname{Re}}

\def\SL{\operatorname{SL}}

\def\Spec{\operatorname{Spec}}

\def\spl{\operatorname{spl}}

\def\toric{\operatorname{toric}}
\def\abtoric{\operatorname{|toric|}}
\def\torus{\operatorname{torus}}
\def\abtorus{\operatorname{|torus|}}

\def\tCu{\tsize\bigcup}

\def\tOp{\tsize\bigoplus}
\def\tp{\tsize\prod}
\def\ts{\tsize\sum}
\def\tsqrt{\tsize\sqrt}

\def\bC{\bold C}

\def\be{\bold e}
\def\bF{\bold F}
\def\bG{\bold G}

\def\bi{\bold i}

\def\bN{\bold N}
\def\bP{\bold P}
\def\bQ{\bold Q}
\def\br{\bold r}
\def\bR{\bold R}
\def\bS{\bold S}

\def\bZ{\bold Z}

\def\cA{{\Cal A}}
\def\cB{{\Cal B}}
\def\cC{{\Cal C}}

\def\cE{{\Cal E}}

\def\cI{{\Cal I}}

\def\cL{{\Cal L}}

\def\cO{{\Cal O}}

\def\fg{{\frak g}}
\def\fh{{\frak h}}

\def\fsl{{\frak s\frak l}}

\def\a{\alpha}
\def\b{\beta}

\def\g{\gamma}
\def\G{\Gamma}
\def\lam{\lambda}
\def\Lam{\Lambda}

\def\sig{\sigma}
\def\Sig{\Sigma}

\def\vf{\varphi}

\def\.{$.\;$}
\def\an{{\text{\rm an}}}
\def\add{{\text{\rm add}}}

\def\gp{{\text{\rm gp}}}
\def\loga{{\text{\rm log}}}
\def\mult{{\text{\rm mult}}}

\def\resp.{\text{resp}.\;}

\def\tra{\overset{\sim}\to{\to}}

\def\triv{{\text{\rm triv}}}

\def\val{{\text{\rm val}}}

\def\O^logten{\cO\^log\otimes}

\let\bs=\backslash

\let\da=\downarrow

\let\la=\leftarrow

\let\lan=\langle
\let\lan=\langle

\let\hra=\hookrightarrow

\let\ox=\otimes
\let\op=\oplus

\let\p=\prime

\let\ran=\rangle

\let\sub=\subset

\let\t=\tilde

\let\x=\times
\let\v=\vee

\def\Fb{\overline{F}}

\def\Dc{\check{D}}
\def\Ec{\check{E}}

\define\LMH{\text{LMH}}

%%%%%

\topmatter

\title
Classifying spaces of degenerating mixed Hodge structures, III: Spaces of nilpotent orbits
\endtitle

\author
Kazuya Kato\footnote{\text{Partially supported by NFS grant DMS 1001729.}},
Chikara Nakayama\footnote{\text{Partially supported by JSPS. KAKENHI (C) No. 18540017,
(C) No. 22540011.}},
Sampei Usui\footnote{\text{Partially supported by JSPS. KAKENHI  (B) No. 19340008.}}
\endauthor

%\date
%(June 12, 2010)
%\enddate

\address
\newline
{\rm Kazuya KATO}
\newline
Department of mathematics
\newline
University of Chicago
\newline
Chicago, Illinois, 60637, USA
\newline
{\tt kkato\@math.uchicago.edu}
\endaddress

\address
\newline
{\rm Chikara NAKAYAMA}
\newline
Graduate School of Science
\newline
Tokyo Institute of Technology 
\newline
Meguro-ku, Tokyo, 152-8551, Japan
\newline
{\tt cnakayam\@math.titech.ac.jp}
\endaddress

\address
\newline
{\rm Sampei USUI}
\newline
Graduate School of Science
\newline
Osaka University
\newline
Toyonaka, Osaka, 560-0043, Japan
\newline
{\tt usui\@math.sci.osaka-u.ac.jp}
\endaddress

\abstract
We construct toroidal partial compactifications 
of the moduli spaces of mixed Hodge structures with polarized graded quotients.
They are moduli spaces of log mixed Hodge structures with polarized graded 
quotients. We construct them as the spaces of nilpotent orbits. 
\endabstract
\footnote"{}"{2000 {\it Mathematics Subject Classification}.
Primary 14C30; Secondary 14D07, 32G20.}

\magnification\magstep1

\endtopmatter
\NoRunningHeads

{\bf Contents}

\medskip

\S1. Log mixed Hodge strucutres

\S2. Moduli spaces of log mixed Hodge structures with polarized graded quotients

\S3. Moduli spaces of log mixed Hodge structures with given graded quotients

\S4. Nilpotent orbits and SL(2)-orbits

\S5. Proofs of the main results

\S6. N\'eron models

\S7. Examples and complements

\vskip20pt

{\bf Introduction}

\medskip

  This is Part III of our series of papers to study degeneration of mixed Hodge structures, and an expanded, full detailed version  
of our previous three announcements \cite{KNU10a}, \cite{KNU10b}, and 
\cite{KNU10c}. 
  In this part, we construct toroidal partial compactifications 
of the moduli spaces of mixed Hodge structures with polarized graded quotients, 
study their properties, and prove that they are moduli spaces of log mixed Hodge structures with polarized graded 
quotients.
  We also apply them to investigate 
 N\'eron models including degenerations of intermediate Jacobians.

\medskip

  We explain the above more precisely. 
  For a complex analytic space $S$ with an fs log structure, a log mixed Hodge structure (LMH, for short) on $S$ is an analytic family parametrized by $S$ of \lq\lq mixed Hodge structures which may have logarithmic degeneration''. 

  The main subject of this paper is to construct a moduli space of LMH. 
A morphism from $S$ to this moduli space is in one-to-one correspondence with an isomorphism class of LMH on $S$.
  This space gives a toroidal partial compactification of the moduli of 
mixed Hodge structure.  
  This moduli space is not necessarily a complex analytic space, but it 
belongs to a category $\cB(\log)$ which contains the category of fs log analytic spaces as a full subcategory. 
  In the category $\cB(\log)$, the space is a fine moduli. 
  See Theorems 2.6.6 and 3.3.3.

  This paper is the mixed Hodge theoretic version of the main part of 
\cite{KU09} (a summary of which is \cite{KU99}), where the degeneration of Hodge structures was studied. 
  We generalize results in \cite{KU09} to the mixed Hodge theoretic 
versions.  
  In particular, we prove that the above moduli space has nice properties, 
especially, that it is a log manifold. 
  See Theorems 2.5.1--2.5.6 and 3.2.8.

  Every proof goes in a parallel way to the pure case \cite{KU09}, and 
we further enhance the story by replacing \lq\lq fan'' in \cite{KU09} 
with \lq\lq weak fan'', which is introduced in \cite{KNU10c}, 
for the reason of applications 
including the study 
of intermediate Jacobians. 
  In fact, with this relaxed concept of fan, the theory can be applied 
more vastly.  
  In this paper, we illustrate a few of such applications, but 
this important subject of the degeneration of intermediate Jacobians 
should be investigated more in a forthcoming paper. 
  Thus, this paper is more general than \cite{KU09} even in the pure case, 
though the proofs are still parallel. 

\medskip

  We organize this paper as follows.  
  First, we present the main results in \S1--\S3. 
  Next, we prove them in \S4 and \S5. 
  Lastly, we gather applications, examples, and complements in \S6--\S7. 

More closely, each section is as follows.

In \S1, we review the notion of log mixed Hodge structures, explaining 
the concept admissibility of the local monodromy, which appears in the mixed 
Hodge case. 

In \S2, we consider the moduli space of log mixed Hodge structures with polarized graded quotients for the weight filtration. 
Let $D$ be a classifying space of mixed Hodge structures with polarized graded quotients (\cite{U84}). 
It is a complex analytic manifold, and  is a mixed Hodge theoretic generalization  of a classifying space of polarized Hodge structures defined by Griffiths (\cite{G68a}). 
In \S2, we consider toroidal partial compactifications $\G\bs D_\Sig$ of $\G\bs D$ for discrete subgroups $\G$ of $\Aut(D)$, and state a theorem 
(2.6.6) that it is a moduli space of LMH with polarized graded quotients. 
Here $D_\Sig$ is a space of nilpotent orbits in the directions in $\Sig$. 
This space sits in the following fundamental diagram. 

$$
\matrix
&&&&&D_{\SL(2),\val}&\hra
&D_{\BS,\val}&\\
&&&&&&&&\\
&&&&&\da&&\da &\\
&&&&&&&&\\
&D_{\Sig,\val} &\la &D_{\Sig,\val}^\sharp
&\to&D_{\SL(2)}
&&D_{\BS}&\\
&&&&&&&&\\
&\da &&\da&&&&&\\
&&&&&&&&\\
&D_{\Sig}&\la
&D_{\Sig}^\sharp&&&&&
\endmatrix
$$
We constructed $D_{\BS}$ in Part I of this series of papers 
(\cite{KNU09}), and $D_{\SL(2)}$ in Part II (\cite{KNU.p}). 

For the degeneration of polarized (pure) Hodge structures, this diagram was constructed in \cite{KU09}.

  The results of this section was announced in \cite{KNU10b}, but 
only in the case of fan, not weak fan.  
  Hence, this section is essentially the same as 
\cite{KNU10b} rewritten by using weak fan.

In \S3, we consider a relative theory fixing a base $S$. We assume that 
we are given $(H_{(w)})_{w\in \bZ}$, where $H_{(w)}$ is a pure 
log Hodge structure of weight $w$ on an object $S$ of $\cB(\log)$. 
  In this section, we consider the moduli of LMH with the given graded 
quotients $(H_{(w)})_{w\in \bZ}$. 
  In the case where all $H_{(w)}$ are polarizable, we state that 
there exists the  moduli space of LMH whose $w$-th graded quotient for the weight filtration is $H_{(w)}$ for each $w$  (3.3.3).
We expect that the construction works without assuming $H_{(w)}$ are polarizable. In some cases, we can actually show that this is true (see 7.1.4).

  For the relation of results of \S2 and those of \S3, see 3.4.

The proofs of the main results in \S2 and \S3 are given in \S4 and \S5, by using the theory of SL(2)-orbits in Part II.

In \S6, we give the theory of N\'eron models which appear for example when an 
intermediate Jacobian degenerates.

In \S7, we give examples and explain the theories of log abelian varieties and log complex tori in \cite{KKN08} from the point of view of this paper.
We also 
explain the necessity of weak fan and discuss completeness of weak fans.

\vskip20pt

\medskip

\head
\S1. Log mixed Hodge structures
\endhead

\head
\S1.1. The category $\cB(\log)$ and relative log manifolds
\endhead

We review the category $\cB(\log)$ 
 introduced in \cite{KU09}, 
which plays a central role in this paper. We introduce a notion \lq\lq relative log manifold'' (1.1.6), which is a relative version of the notion \lq\lq log manifold'' introduced in \cite{KU09}.
 
 \medskip

{\bf 1.1.1.} Strong topology.

Let $Z$ be a local ringed space over $\bC$, and $S$ a subset of $Z$. 
The {\it strong topology} of $S$ in $Z$ is defined as follows. 
A subset $U$ of $S$ is open for this topology if and only if for any analytic space $A$ (analytic space means a complex analytic space) and any morphism $f : A \to Z$ of local ringed spaces over $\bC$ such that $f(A)\subset Z$, $f^{-1}(U)$ is open in $A$. 
This topology is stronger than or equal to the topology as a subspace of $Z$.

\medskip

{\bf 1.1.2.} 
Let $Z$ be a local ringed space (resp.\  an fs log local ringed space, i.e., a local ringed space endowed with an fs log structure  (\cite{KU09} 2.1.5)) over $\bC$, and $S$ a subset of $Z$. 
The {\it strong subspace of $Z$ defined by $S$} is the subset $S$ with the strong topology in $Z$ (1.1.1) and with the inverse image of $\cO_Z$ (resp.\  the inverse images of $\cO_Z$ and $M_Z$, where $M_Z$ denotes the log structure of $Z$).
  A morphism $\iota:S\to Z$ of local ringed spaces (resp.\ fs log local ringed spaces) over $\bC$ is a {\it strong immersion} if $S$ is isomorphic to a strong subspace of $Z$ via $\iota$. 
  If this is the case, it is easy to see by definition that any morphism $A \to Z$ from an 
analytic space (resp.\ an fs log analytic space) uniquely factors through $S$ if and only if it factors through $S$ set theoretically.

\proclaim{Lemma 1.1.3}
  A composite of two strong immersions is a strong immersion. 
\endproclaim 

\demo{Proof}
  It is enough to show the following: 
  Let $Z$ be a local ringed space over $\bC$ and $T \subset S \subset Z$ are 
two subsets of $Z$. 
  Endow $S$ with the strong topology in $Z$. 
  Then, the strong topology of $T$ in $Z$ coincides with the strong 
topology of $T$ in $S$. 

  First, let $U$ be an open subset of $T$ with respect to the strong topology of $T$ in $Z$.  
  Let $A \to S$ be a morphism whose image is contained in $T$. 
  Then, the image of the composite $A \to S \to Z$ is also contained in 
$T$ and the inverse image of $U$ in $A$ is open.  
  Hence, $U$ is open with respect to the strong topology of $T$ in $S$. 
  Conversely, let $U$ be an open subset of $T$ with respect to the strong topology of $T$ 
in $S$. 
  Let $A \to Z$ be a morphism.
  If its image is contained in $T$, then it is contained in $S$, and 
$A \to Z$ factors through $A \to S$ (1.1.2). 
  Since the image of the last morphism is contained in $T$, 
the inverse image of $U$ in $A$ is open.
  Therefore, $U$ is open with respect to the strong topology of $T$ in $Z$. 
\qed\enddemo
  
{\bf 1.1.4.} The category $\cB(\log)$.

As in \cite{KU09}, in the theory of moduli spaces of log Hodge
structures, we have to
enlarge the category of analytic spaces, for the moduli spaces are often not analytic spaces.

Let $\cA$ be the category of analytic spaces and let $\cA(\log)$ be the
category of fs log analytic spaces (i.e., analytic spaces endowed with an fs log structure). We enlarge $\cA$ and $\cA(\log)$ to $\cB$ and $\cB(\log)$,
respectively, as follows. 
In this paper, we will work mainly in the category $\cB(\log)$, and, 
as we did in \cite{KU09}, we will find the moduli spaces in $\cB(\log)$. 

$\cB$ (resp.\ $\cB(\log)$) is the category of all local ringed spaces $S$ 
over $\bC$ (resp.\ local ringed spaces over $\bC$ endowed with a log structure) having the following property:

\medskip

$S$ is locally isomorphic to a strong subspace of an analytic
space (resp.\ fs log analytic space).

\medskip

 The category $\cB(\log)$ has fiber products (\cite{KU09} 3.5.1).

\medskip

{\bf 1.1.5.} Log manifold (\cite{KU09} 3.5.7).

In \cite{KU09}, we saw that moduli spaces of polarized log Hodge structures were logarithmic manifolds (abbreviated as log manifolds in this paper).

An object $S$ of  $\cB(\log)$ is a {\it log manifold} if,  locally on $S$, there are an fs log analytic space $Z$ which is log smooth over $\bC$ (\cite{KU09} 2.1.11), elements $\omega_1,\ldots, \omega_n$ of $\G(Z, \omega^1_Z)$, and an open immersion from $S$ to 
the strong subspace of $Z$ defined by the subset $\{z\in Z\;|\;\text{$\omega_1, \ldots, \omega_n$ are zero in $\omega^1_z$}\}$. Here $\omega^1_Z$ denotes the sheaf of log differential $1$-forms on $Z$ (\cite{KU09} 2.1.7), and $\omega^1_z$ denotes the space of log differential $1$-forms on the point $z=\Spec(\bC)$ which is endowed with the inverse image of the log structure of $Z$. 

\medskip

\noindent 
{\it Example.} 
$S=(\bC\times \bC \smallsetminus  (\{0\}\times \bC)) \cup \{(0,0)\}$ has a natural structure of a log manifold. 
  See \cite{KU09} 0.4.17.

\medskip

{\bf 1.1.6.} {\it Relative log manifold} over an object $S$ of $\cB(\log)$. 

It is an object $T$ of $\cB(\log)$ over $S$ such that, locally on $S$ and on $T$, there exist a log smooth morphism $Y\to X$ of fs log analytic spaces over $\bC$, a morphism $S\to X$, an open subspace $U$ of $S\times_X Y$, 
elements $\omega_1, \ldots, \omega_n \in \G(U, \omega^1_U)$, and  an open immersion from $T$ to the strong subspace of $U$ defined by the subset 
$\{u\in U\;|\;\text{$\omega_1, \ldots, \omega_n$ are zero in $\omega^1_u$}\}$.

\medskip

\head
\S1.2. Relative monodromy filtrations
\endhead

{\bf 1.2.1.} 
Let $V$ be an abelian group, let $W=(W_w)_{w\in \bZ}$ be an increasing filtration on $V$, and let $N\in \End(V, W)$. Assume that $N$ is nilpotent. 

Then an increasing filtration $M=(M_w)_{w\in \bZ}$ on $V$ is called the 
{\it relative monodromy filtration of $N$ with respect to $W$} if it satisfies the following conditions (1) and (2). 

\medskip

(1) $NM_w\subset M_{w-2}$ for any $w\in \bZ$.  

\medskip

(2) For any $w\in \bZ$ and $m\geq 0$, we have an isomorphism $$N^m: \gr^M_{w+m}\gr^W_w@>\sim>> \gr^M_{w-m}\gr^W_w.$$

\medskip

The relative monodromy filtration of $N$ with respect to $W$ need not exist. If it exists, it is unique (\cite{D80} 1.6.13). If it exists, we denote it   by $M(N, W)$ 
as in \cite{SZ85} (2.5). It exists in the case 
 $W_w=V$ and $W_{w-1}=0$ for some $w$  (\cite{D80} 1.6.1).
 
 Relative monodromy filtration is also called relative weight filtration, or relative monodromy weight filtration. 
 
If $V$ is a vector space over a field $K$ and the $W_w$ are $K$-vector subspaces of $V$, and if $N$ is a $K$-linear map, then the relative monodromy filtration, if it exists, is formed by $K$-subspaces of $V$. If $K$ is of characteristic $0$, 
for an increasing filtration $M$ on $V$ consisting of $K$-vector subspaces, $M$ is the relative monodromy filtration $M(N, W)$ if and only if for any $w\in \bZ$, there is an action of the Lie algebra $\fsl(2,K)$ on $\gr^W_w$ satisfying the following conditions (3) and (4). Let $Y:=\pmatrix -1& 0\\0&1\endpmatrix\in \fsl(2,K)$. 

\medskip

(3) For $k\in \bZ$, $M_k\gr^W_w$ is the sum of $\{v\in \gr^W_w\;|\; Yv=\ell v\}$ for all integers $\ell \leq k-w$. 

\medskip

(4) $N:\gr^W_w\to \gr^W_w$ coincides with the action of $\pmatrix 0&1\\0&0\endpmatrix\in \fsl(2,K)$. 

\medskip

We have the following compatibility of relative monodromy filtrations with various operators (direct sum, tensor product, etc.).  For $j=1,2$, let $V_j$ be a vector space over a field $K$ endowed with  an increasing filtration $W_{\bullet}V_j$ consisting of $K$-vector subspaces, and let $N_j\in \End_K(V_j, W_{\bullet}V_j)$ be nilpotent. 

\medskip

{\bf 1.2.1.1.} 
 Assume that the relative monodromy filtration $M_{\bullet}V_j=M(N_j, W_{\bullet}V_j)$ exists. 
Then, the filtration $M_{\bullet}V_1\oplus M_{\bullet}V_2$ on $V_1\oplus V_2$ is the relative monodromy filtration of $N_1\oplus N_2\in \End(V_1\oplus V_2, \, W_{\bullet}V_1\oplus W_{\bullet}V_2)$. Assume furthermore that  $K$ is of characteristic $0$. Then, the filtration $
M_{\bullet}V_1\otimes M_{\bullet}V_2$ on $V_1\otimes V_2$ 
is the relative monodromy filtration of $N_1\otimes 1 + 1 \otimes N_2\in \End(V_1\otimes V_2, \, W_{\bullet}V_1\otimes W_{\bullet}V_2)$. (Here the $w$-th filter of the filtration $M_{\bullet}V_1 \otimes M_{\bullet}V_2$ is $\ts_{j+k=w} M_jV_1\otimes M_kV_2$, and the filtration $W_{\bullet}V_1\otimes W_{\bullet}V_2$ is defined similarly.) 
The filtration $M$ on $H:=\Hom_K(V_1,V_2)$ defined by $M_w=\{f\in \Hom_K(V_1,V_2)\;|\;f(M_kV_1)\subset M_{k+w}V_2\;
\text{for all}\;k\}$ is the relative monodromy filtration of $N\in \End(H, W)$, where $W$ is the filtration on $H$ defined by $W_wH=\{f\in \Hom_K(V_1,V_2)\;|\;f(W_kV_1)\subset W_{k+w}V_2\;\text{for all}\;k\}$ and $N:H\to H$ is the map $f\mapsto 
N_2f-fN_1$. 

\medskip

{\bf 1.2.1.2.} Assume that the relative monodromy filtration $M$ of $N_1\oplus N_2\in \End(V_1\oplus V_2, \, W_{\bullet}V_1\oplus W_{\bullet}V_2)$ exists. Then the relative monodromy filtration $M_{\bullet}V_j= M(N_j, W_{\bullet}V_j)$ exists for $j=1,2$, and $M= M_{\bullet}V_1\oplus  M_{\bullet}V_2$. 
\medskip

By 1.2.1.1 and 1.2.1.2, we have the compatibility with symmetric powers and exterior powers because these are direct summands of tensor products. 

These 1.2.1.1 and 1.2.1.2 are certainly well known. They are proved as follows.

\demo{Proof of 1.2.1.1}
We give the proof for the tensor product. The proof for $\Hom$ is similar and the proof for the direct sum is easy. For $j=1,2$ and $w\in \bZ$, let $\rho_{j,w}$ be the action of $\fsl(2,K)$ on $\gr^W_wV_j$ satisfying the above conditions (3) and (4) (with $M$ and $W$ in (3) and (4) 
replaced by $M_{\bullet}V_j$ and $W_{\bullet}V_j$, respectively,  and with $N$ in (4) replaced by $N_j$). For $w\in \bZ$, let $\rho_w$ be the representation of $\fsl(2,K)$ on $\bigoplus_{k+\ell=w}\gr^W_kV_1\otimes \gr^W_{\ell}V_2$ defined to be $\bigoplus_{k+\ell=w }\;\rho_{1,k}\otimes \rho_{2,\ell}$. Then $\rho_w$ satisfies the above conditions (3) and (4) (with $M$ and $W$ in (3) and (4) 
replaced by $M_{\bullet}V_1\otimes M_{\bullet}V_2$ and $W_{\bullet}V_1\otimes W_{\bullet}V_2$, respectively,  and with $N$ in (4) replaced by $N_1\otimes 1+1\otimes N_2$). Hence $M_{\bullet}V_1\otimes M_{\bullet}V_2$ is the relative monodromy filtration of $N_1\otimes 1+1\otimes N_2$ with respect to $W_{\bullet}V_1\otimes W_{\bullet}V_2$.
\qed
\enddemo

\demo{Proof of 1.2.1.2} 
By the uniqueness of the relative monodromy filtration, 
$M$ is invariant under the maps $V_1\oplus V_2 @>\sim>> V_1\oplus V_2,\;(x,y)\mapsto (ax,by)$ for all $a, b\in K^\times$. 
If $K$ is not $\bF_2$, 
this shows that $M$ is the direct sum of a filtration $M'$ on $V_1$ and a filtration $M''$ on $V_2$. 
  Also, in the case $K=\bF_2$, too, the last property can be proved by using 
a finite extension of $K$ to get $a, b$ such that $a\neq b$. 
It is easy to see that $M'$ (resp.\ $M''$) is the relative monodromy filtration of $N_1$ (resp.\ $N_2$) with respect to $W_{\bullet}V_1$ (resp.\ $W_{\bullet}V_2$). 
\qed
\enddemo

{\it Remark.} 
If the characteristic of $K$ is $p>0$, 1.2.1.1 need not hold for tensor products and $\Hom$. For example, let $V_1$ be of dimension $p$ and with basis $e_1, \dots, e_p$ such that $N_1(e_1)=0$ and $N_1(e_j)=e_{j-1}$ for $2\leq j\leq p$, let $V_2$ be of dimension $2$ with basis $e'_1, e'_2$ such that $N_2(e'_1)=0$ and $N_2(e'_2)=e'_1$, and let $W_0V_j=V_j$ and $W_{-1}V_j=0$ for $j=1,2$. Then for $M=M_{\bullet}V_1\otimes M_{\bullet}V_2$, $\gr^M_p$ (resp. $\gr^M_{-p}$) is of dimension $1$ and generated by the class of $e_p\otimes e'_2$ (resp. $e_1\otimes e'_1$). But $(N_1\otimes 1 +1 \otimes N_2)^p$ induces the zero map $\gr^M_p\to \gr^M_{-p}$ because $(N_1\otimes 1 +1 \otimes N_2)^p=N_1^p\otimes 1 + 1\otimes N_2^p$ and $N_1^p(e_p)=0$, $N_2^p(e'_2)=0$. 
\medskip

In \S6, in the proof of Theorem 6.2.1, we will use the following fact concerning the relative monodromy filtration.

\medskip

{\bf 1.2.1.3.} Let $(V, W)$ and $N$ be at the beginning of 1.2.1, and assume that the relative  monodromy filtration $M=M(N, W)$ exists. Assume that $W_w=0$ for $w\ll 0$. 
Then we have
$$\Ker(N)\cap W_w\sub M_w\quad \text{for all}\;\;w\in \bZ.$$

\demo{Proof of 1.2.1.3} 
We prove $\Ker(N)\cap W_w\sub M_w + W_k$ ($w, k\in \bZ$) by downward induction on $k$ starting from the trivial case $k=w$. Assume that we have already shown $\Ker(N)\cap W_w\sub M_w+W_k$ with some $k\leq w$. Let $x \in \Ker(N) \cap W_w$, and write $x=y+z$ with $y\in M_w$ and $z\in W_k$. Since $N: \gr^M_j(\gr^W_k)\to \gr^M_{j-2}(\gr^W_k)$ are injective for all $j>k$, we have that the map $N:\gr^W_k/M_j\gr^W_k\to \gr^W_k/M_{j-2}\gr^W_k$ is injective if $j\geq k$. Since $N(z)=-N(y)\in M_{w-2}$ and $w\geq k$, this shows that the class of $z$ in $\gr^W_k$ belongs to $M_w\gr^W_k$. Hence we can write $z= y'+z'$ with $y'\in M_w$ and $z'\in W_{k-1}$. Hence $x=(y+y')+z'\in M_w+W_{k-1}$. 
\qed
\enddemo

\medskip

{\bf 1.2.2.} 
Let $V$ be an $\bR$-vector space endowed with an increasing filtration $W$ (by $\bR$-linear subspaces). Let $\sig$ be a cone in some $ \bR$-vector space (this means that $0\in \sig$, $\sig$ is stable under the addition, and stable under  the multiplication by elements of $\bR_{\geq 0}$), and let $\sig_\bR$ be the $\bR$-linear span of $\sig$. Assume that the cone $\sig$ is finitely generated. Assume that we are given an $\bR$-linear map 
$$h:\sig_\bR\to \End_\bR(V, W)$$
whose image consists of mutually commuting nilpotent operators.

We say that the action of $\sig$ on $V$ via $h$ is 
{\it admissible (with respect to $W$)} if there exists a family $(M(\tau, W))_\tau$ of increasing filtrations $M(\tau, W)$ on $V$ given for each face $\tau$ of $\sig$ satisfying the following conditions (1)--(4). 

\medskip

(1) $M(\sig\cap (-\sig), W)=W$.

\medskip

(2) For any face $\tau$ of $\sig$, any $N\in \sig$ and any $w\in \bZ$, we have $h(N)M(\tau, W)_w\subset M(\tau, W)_w$. 

\medskip

(3) For any face $\tau$ of $\sig$, any $N\in \tau$ and any $w\in \bZ$, we have $$h(N)M(\tau, W)_w\subset M(\tau, W)_{w-2}.$$

\medskip

(4) For any faces $\tau$, $\tau'$ of $\sig$ such that $\tau'$ is the smallest face of $\sig$ containing $\tau$ and $N$, $M(\tau', W)$ is the relative monodromy filtration of $h(N)$ with respect to $M(\tau, W)$. 

\medskip

(We remark that in \cite{KKN08} 2.1.4, the above condition (2) was missed.)

We call the above filtration $M(\tau, W)$ the {\it relative monodromy filtration of $\tau$ with respect to $W$}. The condition (4) (by considering the case where $\tau$ and $\sig\cap(-\sig)$ are taken as $\tau'$ and $\tau$ in (4), respectively) implies that $M(\tau, W)=M(N, W)$ for any $N$ in the interior of $\tau$. 

By 1.2.1.1 and 1.2.1.2, we have the 
 following compatibility of this admissibility  with various operators. For $j=1,2$, Let $V_j$ be a vector space over $\bR$ 
 endowed with  an increasing filtration $W_{\bullet}V_j$, and assume that 
we are given an $\bR$-linear map $\sig_\bR \to \End_\bR(V_j, W_{\bullet}V_j)$ 
 whose image consists of mutually commuting nilpotent operators.

 \medskip

{\bf 1.2.2.1.} 
Assume that for $j=1,2$,  the action of $\sig$ on $V_j$ is admissible with respect to $W_{\bullet}V_j$. 
Then, the diagonal action of $\sig$ on $V_1\oplus V_2$ 
is admissible with respect to the filtration 
 $W_{\bullet}V_1\oplus W_{\bullet}V_2$, and $M(\tau, \, W_{\bullet}V_1\oplus W_{\bullet}V_2)=M(\tau, W_{\bullet}V_1)\oplus M(\tau, W_{\bullet}V_2)$ for any face $\tau$ of $\sig$. 
  The action of $\sig$ on $V_1\otimes V_2$ ($N\in \sig$ acts by sending $x\otimes y$ to $Nx\otimes y+x\otimes Ny$)
is admissible with respect to the filtration 
 $W_{\bullet}V_1\otimes W_{\bullet}V_2$, and $M(\tau, \, W_{\bullet}V_1\otimes W_{\bullet}V_2)=M(\tau, W_{\bullet}V_1)\otimes M(\tau, W_{\bullet}V_2)$. 
 The action of $\sig$ on $H:=\Hom_\bR(V_1,V_2)$ ($N\in \sig$ acts by $f\mapsto Nf-fN$) is admissible with respect to the filtration $W_{\bullet}H$ defined by
 $W_wH=\{f\in \Hom_\bR(V_1,V_2)\;|\;f(W_kV_1)\subset W_{k+w}V_2\;\text{for all}\;k\}$,  and $M(\tau, W_{\bullet}H)$ on $H$ is given by $M(\tau, W_{\bullet}H)_w=\{f\in \Hom_\bR(V_1,V_2)\;|\;f(M(\tau, W)_kV_1)\subset M(\tau, W)_{k+w}V_2\;\text{for all}\;k\}$.

\medskip

{\bf 1.2.2.2.} 
Assume that the diagonal action of $\sig$ on $V_1\oplus V_2$ is admissible with respect to $W_{\bullet}V_1\oplus W_{\bullet}V_2$. Then the action of $\sig$ on $V_j$ is admissible with respect to $W_{\bullet}V_j$, and $M(\tau, \,W_{\bullet}V_1\oplus W_{\bullet}V_2)=
M(\tau, W_{\bullet}V_1)\oplus M(\tau, W_{\bullet}V_2)$.
\medskip

\proclaim{Lemma 1.2.3} 
Let the notation be as in $1.2.2$, let $\sig'$ be a finitely generated cone, and assume that 
we are given an $\bR$-linear map $\sig'_\bR\to \sig_\bR$ which sends $\sig'$ in $\sig$. 

\medskip

{\rm(i)} Assume that the action of $\sig$ on $V$ is admissible with respect to $W$.
Then the action of $\sig'$ on $V$ is also admissible with respect to $W$. For a face $\tau'$ of $\sig'$, $M(\tau', W)=M(\tau,W)$, where $\tau$ denotes the smallest face of $\sig$ which contains the image of $\tau'$. 

\medskip

{\rm(ii)} Assume that the action of $\sig'$ on $V$ is admissible with respect to $W$ and that the map $\sig'\to \sig$ is surjective. Then the action of $\sig$ on $V$ is admissible with respect to $W$. For a face $\tau$ of $\sig$, $M(\tau,W)=M(\tau',W)$, where $\tau'$ denotes the inverse image of $\tau$ in $\sig'$.  
\endproclaim

\demo{Proof} (i) This is essentially 
included in \cite{KKN08} 2.1.5. The proof is easy.

(ii) Denote the given map $\sig'\to \sig$ by $f$.  
Let $\alpha$ and $\beta$ be faces of $\sig$, let $N\in \beta$, and assume that $\beta$ is the smallest face of $\sig$ which contains $\alpha$ and $N$. It is sufficient to prove that $M(f^{-1}(\beta),W)$ is the relative monodromy filtration of $h(N)$ with respect to $M(f^{-1}(\alpha),W)$. Let $x$ be an element of the interior of $f^{-1}(\beta)$. Since $\beta$ is the smallest face of $\sig$ which contains $\alpha$ and $N$, there are $y\in \beta$, $z\in \alpha$, $c\in \bR_{>0}$ such that $f(x)+y = z+cN$. Take $y'\in f^{-1}(\beta)$ and $z'\in f^{-1}(\alpha)$ such that $f(y')=y$ and $f(z')=z$, and let $u=x+y'$. Since $u$ is in the interior of $f^{-1}(\beta)$, $M(f^{-1}(\beta),W)$ is the relative monodromy filtration of $h(f(u))$ with respect to $M(f^{-1}(\alpha),W)$. Since $f(u)= f(z')+cN$, this shows that $M(f^{-1}(\beta),W)$ is the relative monodromy filtration of $h(N)$ with respect to $M(f^{-1}(\alpha),W)$.
\qed
\enddemo

{\bf 1.2.4.} 
Let $S$ be an object of $\cB(\log)$, and let the topological space $S^{\log}$ and the continuous proper map $\tau:S^{\log}\to S$ be as in  \cite{KN99}, \cite{KU09} 2.2. Let $L=(L,W)$ be a locally constant sheaf $L$ of finite dimensional $\bR$-vector spaces over $S^{\log}$ endowed with an increasing filtration $W$ consisting of locally constant $\bR$-subsheaves. 
We say that the local monodromy of $L$ is {\it admissible} if the following two conditions (1) and (2) are satisfied for any $s\in S$ and any $t\in s^{\log}=\tau^{-1}(s)\sub S^{\log}$. 

\medskip

(1) For any $\g \in \pi_1(s^{\log})$, the action of $\g$ on the stalk $L_t$ is unipotent.

\medskip

(2) Let $C(s):= \Hom((M_S/\cO_S^\times)_s, \bR_{\geq 0}^{\add})$, where $ \bR_{\geq 0}^{\add}$ denotes $ \bR_{\geq 0}$ regarded as an additive monoid. Then 
the action of the cone $C(s)$ on $L_t$ via $C(s)_\bR = \bR\ox\pi_1(s^{\log})\to \End_\bR(L_t,W)$, defined by $\g \mapsto \log(\g)$ ($\g\in \pi_1(s^{\log})$), is admissible with respect to $W$. 

\medskip

\proclaim{Lemma 1.2.5}
Admissibility of local monodromy is preserved by pulling back by a morphism $S'\to S$ in $\cB(\log)$. 
\endproclaim

\demo{Proof} 
This follows from 1.2.3 (i). 
\qed
\enddemo

\noindent 
{\bf Remark 1.2.6.} 
There are weaker versions of the notion of  admissibility of  a cone with respect to $W$. 

For example, one version is 
\medskip

(1) The relative monodrmy filtration $M(N, W)$ exists for any $N\in \sig$. 
\medskip

Another version is

\medskip

(2) The relative monodromy filtration $M(N, W)$ exists for any $N\in \sig$, and it depends only on the smallest face of $\sig$ which contains $N$.

\medskip

In the situation with polarized graded quotients explained in 1.3.3 below, it is shown that these weaker admissibilities coincide with the admissibility considered in this subsection \S1.2. See 1.3.2--1.3.4.

In the papers   
\cite{KNU10a}, \cite{KNU10b}, and \cite{KNU10c}, we adopted the admissibility in (2) above. But the difference with the admissibility in this subsection is not essential for this paper.  
  See 2.2.5. 

\head
\S1.3. Log mixed Hodge structures
\endhead

{\bf 1.3.1.} 
Log mixed Hodge structure (LMH, for short) on $S$ (\cite{KKN08} 2.3, \cite{KU09} 2.6). 

Let $S$ be an object in $\cB(\log)$. Recall 
that we have a sheaf of commutative rings $\cO_S^{\log}$ on $S^{\log}$ which is an algebra over the inverse image $\tau^{-1}(\cO_S)$ of $\cO_S$ (\cite{KN99}, \cite{KU09} 2.2).

A {\it pre-log mixed Hodge structure (pre-LMH) on $S$} is a triple $(H_\bZ, W, H_\cO)$, where

$H_\bZ$ is a locally constant sheaf of finitely generated free 
$\bZ$-modules on $S^\loga$,

$W$ is an increasing filtration on $H_\bR$ by 
locally constant rational subsheaves, and

$H_\cO$ is an $\cO_S$-module on $S$ which is locally free of finite rank, 
endowed with a
decreasing filtration $(F^pH_\cO)_p$ and endowed with an isomorphism of $\cO_S^\loga$-modules 
$\cO_S^\loga \otimes_\bZ H_\bZ \simeq \cO_S^\loga \otimes_{\cO_S} H_\cO$ on $S^\loga$,

\noindent 
satisfying the condition that $F^pH_{\cO}$ and 
$H_{\cO}/F^pH_{\cO}$ are locally free of finite rank for any $p \in \bZ$.

Denote $F^p:= \cO_S^\loga \otimes_{\cO_S} F^pH_\cO$ $(p\in\bZ)$ the filtration on $\cO_S^\loga$-module $\cO_S^\loga \otimes_{\cO_S} H_\cO \simeq \cO_S^\loga \otimes_\bZ H_\bZ$.
The filtrations $F^pH_{\cO}$ and $F^p$ determine each other so that we can 
define a pre-LMH as a triple $(H_{\bZ}, W, F)$, where $F$ is a decreasing 
filtration on $\cO^{\log}_S \otimes_{\bZ}H_{\bZ}$, satisfying the corresponding 
conditions.  
In the rest of this paper, we freely use both formulations. 

\medskip

Let $s\in S$ and let $t\in \tau^{-1}(s)\sub S^{\log}$. By a {\it specialization} at $t$, we mean a ring homomorphism $\cO_{S,t}^{\log}\to \bC$ such that the composition
$\cO_{S,s}\to \cO_{S,t}^{\log}\to \bC$ is the evaluation at $s$. Since the evaluation at $s$ is the unique $\bC$-algebra homomorphism  $\cO_{S,s}\to \bC$,  any $\bC$-algebra homomorphism $\cO_{S,t}^{\log} \to \bC$ is a specialization at $t$. 

\medskip

Let $s$ be an fs log point (i.e., an fs log analytic space whose underlying analytic space is a one-point set endowed with the ring $\bC$).

A pre-LMH $(H_\bZ,W,H_\cO)$ on $s$ is a {\it log mixed Hodge structure (LMH) on $s$}, if it satisfies the following conditions (1)--(3).
\smallskip

(1) The local monodromy of $(H_\bR, W)$ is admissible in the sense of 1.2.4.

\smallskip

(2) Let $\nabla=d\ox 1_{H_\bZ}: \cO_s^\loga\ox_\bZ H_\bZ \to \omega_s^{1,\log}\ox_\bZ H_\bZ$.
Then
$$
\nabla F^p \sub \omega_s^{1,\log}\ox_{\cO_s^\loga}F^{p-1} \quad 
\text{for all $p$}.
$$

(3) Let $t \in s^\loga$. 
For a specialization $a: \cO_{s,t}^\loga \to \bC$ at $t$, let $F(a):=\bC \ox_{\cO_{s,t}^\loga}F_t$ be the filtration on $\bC \ox_\bZ H_{\bZ,t}$.
If $a$ is sufficiently twisted in the sense explained below, then, for any face $\tau$ of the cone $C(s)$ in 1.2.4, $(H_{\bZ,t},M(\tau,W),F(a))$ is a mixed Hodge structure in the usual sense.

Here, fixing a finite family $(q_j)_{1\leq j\leq n}$ of elements of $M_s\smallsetminus \bC^\times$ such that 
$(q_j\bmod \bC^\times)_j$ generates $M_s/\bC^\times$ (note that $\bC^\times =\cO_s^\times \subset M_s$),  we say $a$ is sufficiently twisted if 
$\exp(a(\log(q_j)))$ are sufficiently near to $0$.  This notion \lq\lq sufficiently twisted'' is in fact independent of the choice of $(q_j)_j$.\medskip

Let $S$ be an object in $\cB(\log)$ again. 

A pre-LMH on $S$ is a {\it log mixed Hodge structure (LMH)} on $S$ if its pull back to each fs log point $s \in S$ is an LMH on $s$ (\cite{KKN08} 2.3, \cite{KU09} 2.6). 

Direct sums, tensor products, symmetric powers, exterior powers, and inner-Hom of pre-LMH are defined in the evident ways. By using 1.2.2.1 and 1.2.2.2, we see that these operations preserve LMH.

\medskip
For $w\in \bZ$, we call
a  pre-LMH (resp.\ LMH) such that $W_w=H_\bR$ and $W_{w-1}=0$  simply a pre-log Hodge structure (pre-LH) (resp.\ log Hodge structure (LH)) of weight $w$. 
In this case, we often omit the $W$. 

\medskip

 {\bf 1.3.2.} Polarized log Hodge structure (PLH) on $S$ (\cite{KU99} \S5, \cite{KKN08} 2.5, 
\cite{KU09} 2.4).

Let $w$ be an integer.

A {\it pre-polarized log Hodge structure (pre-PLH) of weight $w$} is a triple $(H_{\bZ}, H_{\cO}, \langle\;,\;\rangle)$, where 
\smallskip

$(H_{\bZ}, H_{\cO})$ is a pre-LH of weight $w$, and 
\smallskip

$\langle\;,\;\rangle: H_{\bR} \times H_{\bR} \to \bR$ is a rational non-degenerate $(-1)^w$-symmetric bilinear form, 
\smallskip

\noindent 
satisfying $\lan F^p, F^q \ran=0$ when $p+q>w$. 
\medskip

A pre-PLH $(H_{\bZ}, H_{\cO}, \langle\;,\;\rangle)$ of weight $w$ is a 
{\it polarized log Hodge structure (PLH) of weight $w$} if for any $s\in S$, the pull back of $(H_{\bZ}, H_{\cO}, \langle\;,\;\rangle)$ to the fs log point $s$ has the 
property 1.3.1 (2) and the 
following property: In the notation $t$, $a$ and $F(a)$ in 1.3.1 (3), $(H_{\bZ,t}, \lan\;,\;\ran, F(a))$ is a polarized Hodge structure of weight $w$ in the usual sense for any sufficiently twisted specialization $a$ (1.3.1).

  A PLH of weight $w$ is an LH of weight $w$. 
  This is by results of \cite{CK82} and \cite{Scm73}.

\medskip

{\bf 1.3.3.}  
By results of Cattani-Kaplan, Schmid, and Kashiwara, 
under the existence of polarizations on graded quotients, the definition of LMH becomes simpler as follows. 

 Let $s$ be an fs log point 
and let $(H_{\bZ}, W, H_{\cO})$ be a pre-LMH on $s$. 
Assume that $(H_{\bZ}, W, H_{\cO})$ satisfies the condition 1.3.1 (2). 
Assume moreover that, for each $w \in \bZ$, a rational non-degenerate $(-1)^w$-symmetric bilinear form $\langle\;,\;\rangle_w: H(\gr^W_w)_{\bR} \times H(\gr^W_w)_{\bR} \to \bR$ is given and that the triple $(H(\gr^W_w)_{\bZ}, H(\gr^W_w)_{\cO}, \langle\;,\;\rangle_w)$ is a PLH. 
Let $C(s)$ be the cone in 1.2.4. 
Fix a set of generators $\g_1,\ldots, \g_n$ of $C(s)$. 
By \cite{KU09} 2.3.3, (1) in 1.2.4 is satisfied for any $t \in s^{\log}$ and we can define $N_j := \log(\g_j) \in \End(H_{\bR,t}, W)$ for all $j$. 

\proclaim{Proposition 1.3.4} 
Let the assumption be as in $1.3.3$.  
Then, the following are equivalent.
\smallskip

$(1)$ $(H_{\bZ}, W, H_{\cO})$ is an LMH.
\smallskip

$(2)$ The local monodromy of $(H_{\bR}, W)$ is admissible.
\smallskip

$(3)$ The relative monodromy filtrations $M(N_j,W)$ exist for all $j$. 
\endproclaim

\demo{Proof} The implications (1) $\Rightarrow$ (2) and (2) $\Rightarrow$ (3) are evident.

We prove (2) $\Rightarrow$ (1). Assume (2). 
Since each gr is a PLH, it is also an LH.
Hence 1.3.1 (3) is satisfied for each gr.  By [K86] 5.2.1,
1.3.1 (3) is satisfied, that is, (1) holds.

We prove (3) $\Rightarrow$ (2). Assume (3). 
  Let $F$ be the specialization (1.3.1) of the Hodge filtration with respect to 
some $a: \cO_{s,t}^\loga \to \bC$ at $t$. 
  Then, 
$(H_{\bC,t}; W_{\bC,t}; F, \bar F; N_1,\ldots, N_n)$ is a pre-IMHM in the 
sense of \cite{K86} 4.2. 
  By 4.4.1 of \cite{K86}, it is an IMHM. 
  In particular, 
for any subset $J\sub \{1, \ldots, n\}$, 
$M(\sum_{j\in J}N_j, W)$ exists.  
Let $M(\tau,W):=M(\sum_{j\in J}N_j, W)$, where $\tau$ is the face of 
$C(s)$ generated by  $\g_j$ for $j$ running in $J$.
Then, this family satisfies the condition 1.2.2 (4) by 5.2.5 of \cite{K86}.
Hence, the local monodromy of $(H_{\bR}, W)$ is admissible. 
\qed
\enddemo

{\bf 1.3.5.} 
  Let $S$ be an object of $\cB(\log)$. 
  An {\it LMH with polarized graded quotients on $S$} is a 
quadruple $H=(H_\bZ, W, H_{\cO}, (\langle\;,\;\rangle_w)_w)$
such that $(H_\bZ, W, H_{\cO})$ is an LMH on $S$ and 
$(H(\gr^W_w)_{\bZ}, H(\gr^W_w)_{\cO}, \langle\;,\;\rangle_w)$ is a PLH 
of weight $w$ for any $w$. 

\medskip

\example{Example 1.3.6} 
Let $S=\Delta^{n+r}$, where $\Delta$ is the unit disk. 
Let $S^*=(\Delta^*)^n \times \Delta^r$, where $\Delta^*=\Delta-\{0\}$. 
Let $X$ be an analytic space and $X \to S$ be a projective morphism which is smooth over $S^*$.
Let $X^*=X \times_SS^*$. 
Let $E$ be a divisor on $X$ such that $E \cap X^*$ is relatively normal 
crossing over $S^*$. 
Let $m \in \bZ$. 
Let $h: X^* \smallsetminus E \to S^*$. 
Assume that the local monodromy of $H_{\bR} := R^mh_*(\bR)$ along $E$ is unipotent. 
Then, the natural variation of mixed Hodge structure (VMHS) with polarized graded quotients on $S^*$ underlain by $H_{\bR}$ 
canonically extends to an LMH with polarized graded quotients. 
This was proved by Steenbrink-Zucker \cite{SZ85}, Kashiwara \cite{K86}, 
Saito \cite{Sa90}, and Fujino \cite{F04}. 
See the explanation of the last part of \cite{KNU08} 12.10.
\endexample

\head
\S2. Moduli spaces of log mixed Hodge structures with polarized graded quotients
\endhead

\head
\S2.1. Classifying space $D$ 
\endhead

{\bf 2.1.1. } In this \S2, we fix $\Lambda=(H_0, W, (\langle\;,\;\rangle_w)_w, (h^{p,q})_{p, q})$, where 

\medskip

$H_0$ is a finitely generated free $\bZ$-module,

$W$ is an increasing rational filtration on $H_{0, \bR} = \bR\otimes H_0$,

$\langle\;,\;\rangle_w$ for each $w\in \bZ$ is a rational non-degenerate $\bR$-bilinear form  $\gr^W_w\times \gr^W_w\to \bR$ which is symmetric if $w$ is even and is anti-symmetric if $w$ is odd, 

$h^{p,q}$ are non-negative integers given for each $(p, q)\in \bZ^2$
satisfying the following conditions (1)--(3).

\smallskip

(1) $\sum_{p, q} h^{p,q}= \text{rank}_\bZ(H_0)$,
\smallskip

(2) $\sum_{p+q=w} h^{p,q}= \dim_\bR(\gr^W_w)$ for any $w\in \bZ$,
\smallskip

(3) $h^{p,q}=h^{q,p}$ for any $(p, q)$.

\medskip

{\bf 2.1.2.} 
Let $D$ be the classifying space of gradedly polarized mixed 
Hodge structures in \cite{U84} associated with the data fixed in 2.1.1.
As a set, $D$ consists of all increasing filtrations $F$ on $H_{0,\bC} = \bC \otimes H_0$ such that $(H_0, 
W, (\langle\;,\;\rangle_w)_w, F)$ is a gradedly polarized mixed Hodge structures with $\dim_\bC F^p(\gr^W_{p+q})/F^{p+1}(\gr^W_{p+q}) 
= h^{p,q}$ for all $p, q$. 
Let $\Dc$ be its compact dual.  
Then $D$ and $\Dc$ are complex analytic manifolds, and $D$ is open in $\Dc$. 
Cf.\ \cite{KNU09} 1.5.

\medskip

{\bf 2.1.3.} We recall the  notation used in this series of papers.

As in \cite{KNU09} 1.6, for $A=\bZ, \bQ, \bR, \bC$, let 
$G_A$ be the group of the $A$-automorphisms of $(H_{0,A}, W)$ whose $\gr^W_w$ 
are compatible with $\langle\;,\;\rangle_w$ for all $w$. 
 Then $G_\bC$ acts on $\Dc$ and $G_\bR$ preserves $D$. 

As in \cite{KNU09} 1.7, 
for $A= \bQ, \bR, \bC$, let $\fg_A$ be the associated Lie algebra of $G_A$, which is 
identified with the set of all $A$-endomorphisms 
$N$ of $(H_{0,A}, W)$ whose $\gr^W_w$ satisfies 
$\lan\gr^W_w(N)(x),y\ran_w+\lan x,\gr^W_w(N)(y)\ran_w=0$ for all $w, x, y$. 

As in \cite{KNU09} 1.6 and 1.7, let $G_{A,u}=\{\g\in G_A\;|\;\gr^W(g)=1\}$,
$\fg_{A,u}=\{N\in \fg_A\,|\,\gr^W(N)=0\}$. Then $G_A/G_{A,u}=G_A(\gr^W):=\prod_w G_A(\gr^W_w)$ and $\fg_A/\fg_{A,u}=\fg_A(\gr^W):=\prod_w\fg_A(\gr^W_w)$, where $G_A(\gr^W)$ (resp. $\fg_A(\gr^W_w)$) is \lq\lq the  $G_A$ (resp. $\fg_A$) for $\gr^W_w$''. 

\medskip

\head
\S2.2.  The sets $D_\Sig$ and $D_{\Sig}^\sharp$
\endhead

  Here we define the sets $D_\Sig$ and $D_{\Sig}^\sharp$ in our 
mixed Hodge case.  
See \cite{KU09} for the pure case. 
  Everything is parallel to and generalize the pure case, except that 
we use weak fans instead of fans.

\medskip

{\bf 2.2.1.} {\it Nilpotent cone.}
A subset $\sig$ of $\fg_\bR$ is called a {\it nilpotent cone} if the following conditions (1) and (2) are satisfied.

\medskip

(1) All elements of $\sig$ are nilpotent and $NN'=N'N$ for any $N, N'\in \sig$, as linear maps $H_{0,\bR}\to H_{0,\bR}$.

\medskip

(2) $\sig$ is finitely generated. That is, there is a finite family $(N_j)_j$ of elements of $\sig$ such that
$\sig=\ts_j \bR_{\geq 0}N_j$. 

\medskip

A nilpotent cone is said to be {\it rational} if we can take $N_j\in \fg_\bQ$ in the above condition (2). 
 
For a nilpotent cone $\sig$, let $\sig_\bR$ be the $\bR$-linear span of $\sig$ in $\fg_\bR$, and let $\sig_\bC$ be the $\bC$-linear span of $\sig$ in $\fg_\bC$.
\medskip

We say that $\sig$ is {\it admissible} if the action of $\sig$ on $H_{0,\bR}$ is admissible with respect to $W$ in the sense of 1.2.2. 

\medskip

{\bf 2.2.2.} Let $N_1, \dots, N_n\in \fg_\bR$ be mutually commuting nilpotent elements and let $F\in \Dc$. We say $(N_1, \dots, N_n, F)$ 
{\it generates a nilpotent orbit} if the following (1)--(3) are satisfied.

\medskip

(1)  The action of the cone $\ts_{j=1}^n \bR_{\geq 0}N_j$ on $H_{0,\bR}$ 
is admissible (1.2.2) with respect to $W$.

\medskip

(2) $N_jF^p\sub F^{p-1}$ for $1\leq j\leq n$ and for any $p\in \bZ$.

\medskip

(3) If $y_j \in \bR_{>0}$ and $y_j$ are sufficiently large, we have $\exp(\sum_{j=1}^n iy_jN_j)F\in D$.

\medskip
Note that these conditions depend only on the cone $\sig:=\ts_{j=1}^n \bR_{\geq 0}N_j$ and $F$. We say also that the pair $(\sig, F)$ generates a nilpotent orbit if these conditions are satisfied.

Let $\sig$ be a nilpotent cone. 
 A subset $Z$ of $\Dc$ is called a
{\it $\sigma$-nilpotent orbit} (resp.\ {\it $\sig$-nilpotent $i$-orbit}) if the
 following conditions
(4) and (5)  are satisfied for some $F\in Z$.  

\medskip

(4) $Z=\exp(\sig_\bC)F$ (resp.\ $Z=\exp(i\sig_\bR)F$). 

\medskip

(5) $(\sig, F)$ generates a nilpotent orbit.

\medskip

Such a pair $(\sig, Z)$ is called a {\it nilpotent orbit} (resp.\ {\it nilpotent $i$-orbit}).

Note that if (4) and (5) are satisfied for one $F\in Z$, they are satisfied for all $F\in Z$.

\medskip

{\it Remark} 1.  The notion of nilpotent orbit is essentially  the
same as the notion of LMH on an fs log point.  Cf. [KU09] 2.5
for the pure case.

\medskip

{\it Remark} 2. 
  The notion of nilpotent orbit is closely related to the notion 
infinitesimal mixed Hodge module (IMHM) by Kashiwara (\cite{K86}). 
  A precise relation is as follows.  
(The statement in \cite{KNU10b} 2.1.3 was not accurate.) 

  Let $H_0$ and $W$ be as in 2.1.1. 
  Let $F$ be a decreasing filtration on $H_{0,\bC}$. 
  Let $\bar F$ be its complex conjugation. 
  Let $N_1,\dots, N_n$ be a mutually commuting set of nilpotent 
endomorphisms of $H_{0,\bQ}$. 
  Then, 
$(H_{0,\bC}; W_{\bC}; F, \bar F; N_1,\dots, N_n)$ 
is an IMHM if and only if there exists a 
$\Lambda$ extending $(H_0, W)$ such that 
$N_1, \ldots, N_n \in \fg_{\bR}$ and $F\in \Dc$, and 
that ($N_1,\dots, N_n, F)$ generates a nilpotent orbit. 
  This was seen in the course of the proof of 1.3.4.

\medskip

{\bf 2.2.3.} Weak fan in $\fg_\bQ$ (\cite{KNU10c}). 

A {\it weak fan $\Sig$ in $\fg_{\bQ}$} is a non-empty set of sharp rational nilpotent cones in $\fg_\bR$ satisfying the following conditions (1) and (2). Here a nilpotent cone $\sig$ is said to be {\it sharp} if $\sig\cap (-\sig)=\{0\}$. 

\medskip

(1) If $\sig\in \Sig$, all faces of $\sig$ belong to $\Sig$.

\medskip

(2) Let $\sig, \sig' \in \Sig$, and assume that $\sig$ and $\sig'$ have a common interior point. Assume that 
there is an $F\in \Dc$ such that $(\sig, F)$ and $(\sig', F)$ generate nilpotent orbits. Then $\sig=\sig'$.  

\medskip

A {\it fan} $\Sig$ in $\fg_{\bQ}$ is defined similarly by 
replacing the condition (2) in the above with the condition \lq\lq$\sig \cap \sig'$ is 
a face of $\sig$ for any $\sig, \sig' \in \Sig$''.

For simplicity, we call a weak fan (resp.\ a fan) in $\fg_\bQ$ a weak fan (resp.\ a fan).

An example of a fan is
$$\Xi=\{\bR_{\geq 0}N\;|\; N\in \fg_\bQ, \;N\;\text{is nilpotent}\}.$$

  The next lemma, which was announced in \cite{KNU10c} 1.7, 
implies that a fan in $\fg_{\bQ}$  is a weak fan. 

 \medskip

\proclaim{Lemma 2.2.4}
In the definition $2.2.3$, under the condition $(1)$, 
the condition $(2)$ is equivalent to $(2)'$ below and also 
is equivalent to $(2)''$ below. 

\medskip

$(2)'$ If $\sig, \sig' \in \Sig$ and if 
there is an $F\in \Dc$ such that $(\sig, F)$, $(\sig', F)$, and 
$(\sig \cap \sig', F)$ generate nilpotent orbits, then 
$\sig \cap \sig'$ is a face of $\sig$.

\medskip

$(2)''$ Let $\tau$ be a nilpotent cone. 
Let $(\tau,F)$ $(F \in \Dc)$ generate a nilpotent orbit. 
If the set $A:=\{\sig \in \Sig\,\vert\, \tau \subset \sig, (\sig,F)$ generates 
a nilpotent orbit$\}$ is not empty, $A$ 
has a smallest element, and it is a face of any element of $A$. 
\endproclaim

\demo{Proof}
$(2)''\Rightarrow (2)'$.  Let $\sig, \sig'$ and $F$ be as in $(2)'$.
 Let $\tau=\sig \cap \sig'$.
 Then, $\sig$ and $\sig'$ belong to $A$.
 Let $\sig_0$ be the smallest element of $A$.
 Then, $\sig_0$ is a face of $\sig$ and also is a face of $\sig'$, and hence $\sig_0 \sub \sig \cap
\sig'=\tau$.
 On the other hand, since $\sig_0 \in A$, we have $\tau \sub \sig_0$.
 Thus, $\sig_0 = \tau$.
 Hence, $\sig \cap \sig'=\tau=\sig_0$ is a face of $\sig$.

$(2)'\Rightarrow (2)$. Let $\sig, \sig'$ and $F$ be as in (2).
 Since $\sig \cap \sig'$ contains an interior point of $\sig$,
$(\sig \cap \sig', F)$ generates a nilpotent orbit.
 Hence, $(2)'$ implies that $\sig \cap \sig'$ is a face of $\sig$, but
again, since $\sig \cap \sig'$ contains an interior point of $\sig$,
we have $\sig \cap \sig'=\sig$.
 By symmetry, $\sig \cap \sig'=\sig'$.
 Therefore, $\sig=\sig'$.

$(2)\Rightarrow (2)''$.  Let $(\tau, F)$ be as in $(2)''$.
 For an element $\sig \in A$, we write as
$\sig_0$ the face of $\sig$ spanned by $\tau$.
 Then, $(\sig_0,F)$ generates a nilpotent orbit because both
$(\tau,F)$ and $(\sig,F)$ generate nilpotent orbits.
 Thus, $\sig_0 \in A$.
 It is enough to show $\sig_0 = \sig'_0$ for any $\sig, \sig' \in A$,
because, then, this common $\sig_0$ is the smallest element of $A$ as
soon as $A$ is not empty, and it is indeed a face of any $\sig \in A$.
 But, since an interior point of $\tau$ is a common interior point of
$\sig_0$ and $\sig'_0$, we have certainly $\sig_0 = \sig'_0$ by (2). 
\qed
\enddemo

{\bf 2.2.5.} 
Let $\Sig$ be a weak fan in $\fg_\bQ$.
Let $D_\Sig$ (resp.\ $D^{\sharp}_\Sig$) be the set of all nilpotent orbits (resp.\ nilpotent $i$-orbits) $(\sig, Z)$ with $\sig\in \Sig$. 

\medskip

We have embeddings 
$$D \subset D_\Sig, \quad F \mapsto (\{0\}, \{F\}),$$
$$D \subset D_\Sig^{\sharp}, \quad F \mapsto (\{0\}, \{F\}).$$

We have a canonical surjection
$$D_{\Sig}^{\sharp}\to D_{\Sig}, \quad (\sig, Z)\mapsto (\sig, \exp(\sig_\bC)Z),$$
which is compatible with the above embeddings of $D$.

\medskip

{\it Remark} 1. If $'\Sig\sub \Sig$ denotes the subset of $\Sig$ consisting of all admissible  $\sig\in \Sig$, and if $''\Sig\sub {'\Sig}$ denotes the subset of $\Sig$ consisting of all $\sig\in \Sig$ such that a $\sig$-nilpotent orbit exists, then $'\Sig$ and $''\Sig$ are weak fans, and we have $D_{\Sig}=D_{'\Sig}=D_{''\Sig}$,
$D_{\Sig}^{\sharp}=D_{'\Sig}^{\sharp}=D_{''\Sig}^{\sharp}$.
 
\medskip
{\it Remark} 2.  In our previous papers \cite{KNU10b} and \cite{KNU10c}, we adopted the following formulations (1)--(3) which are slightly different from those in the present paper.   In those papers:

\smallskip
(1)  Nilpotent cone was assumed to be sharp.

\smallskip

(2) In the definition of nilpotent orbit, we adopted the admissibility condition 1.2.6 (2) on the cone instead of the condition 2.2.2 (1).

For sharp cones, the definition of nilpotent
orbit in [KNU10b] and [KNU10c] is equivalent to that in this paper.
Cf. 1.2.6, 1.3.4.  

\smallskip

(3) In the definition of weak fan (resp.\ fan), in \cite{KNU10c} (resp.\ \cite{KNU10b} and\ \cite{KNU10c}), 
we put the condition that any cone in a weak fan (resp.\ fan) should satisfy the admissibility condition 1.2.6 (2). But
 we do not do it in the present paper, for we prefer a simpler definition here. 
 
 \smallskip
 
 These changes of the formulations do not make any difference for the space of nilpotent orbits $D_{\Sig}$ (see Remark 1 above) and for the moduli spaces $\G\bs D_{\Sig}$ of LMH which appear in \S2.4 below.

In the pure case, the definitions of nilpotent cone,  nilpotent orbit, and fan in this paper coincide with those in \cite{KU09}.

\medskip

{\bf 2.2.6.} Compatibility of a weak fan and a subgroup of $G_\bZ$. 

\medskip

Let $\Sig$ be a weak fan in $\fg_\bQ$ and let $\G$ be a subgroup of $G_\bZ$.

We say that $\Sig$ and $\G$ are {\it compatible} if the following condition (1) is satisfied.

\medskip

(1) If $\g\in \G$ and $\sig\in \Sig$, then $\Ad(\g)\sig\in \Sig$. 

\medskip

If $\Sig$ and $\G$ are compatible, 
$\Gamma$ acts on $D_\Sig$ and also on $D_{\Sig}^{\sharp}$ by $(\sig, Z)\mapsto (\Ad(\gamma)\sig,
\gamma Z)$ ($\gamma\in \Gamma$).

\medskip

We say that 
$\Sig$ and $\G$ are {\it strongly compatible} if they are compatible and furthermore the following condition (2) is satisfied.

\medskip

(2) If $\sig\in \Sig$, any element of $\sig$ is a sum of elements of the form
$aN$, where $a\in \bR_{\geq 0}$ and $N\in \sig$ satisfying $\exp(N)\in \G$.

\medskip

{\bf 2.2.7.}  For a sharp rational nilpotent cone $\sig$, let $\text{face}(\sig)$ be the set of all faces of $\sig$. 
It is  
a fan in $\fg_\bQ$. 
Let $$D_\sig:=D_{\text{face}(\sig)}, \quad D_\sig^\sharp:=D_{\text{face}(\sig)}^\sharp.$$ 
 Let $$\G(\sig):=\G \cap \exp(\sig)  \sub  \G(\sig)^{\gp} = 
 \G \cap \exp(\sig_\bR) =\{ab^{-1}\;|\;a, b\in \G(\sig)\}\sub G_\bZ.$$  Then $\G(\sig)$ is a sharp torsion free fs monoid and $\G(\sig)^{\gp}$ is a finitely generated free abelian group. 

  Let $\Sig$ and $\G$ be as in 2.2.6.  
  Assume that they are strongly compatible.
  Then, for any $\sig\in \Sig$, 
the fan $\text{face}(\sig)$ and the group $\G(\sig)^\gp$ are strongly compatible. 

\medskip

{\bf 2.2.8.} We will use the following fact in the proof of Theorem 6.2.1.

\medskip

 Let $\sig\sub \fg_\bR$ be an admissible nilpotent cone. 
Then the adjoint action of $\sig$ on $\fg_\bR$ is admissible.

\medskip

{\it Proof.}
By 1.2.2.1, 
the adjoint action of $\sig$ on $A:=\Hom_\bR(H_{0,\bR},H_{0,\bR})$ is admissible. For a face $\tau$ of $\sig$, let $M(\tau)$ be the relative monodromy filtration of $\tau$ on $A$ for this admissible action. 
 The adjoint action of $\sig$ on
$\fg_\bR(\gr^W_w)$ for $w<0$ is admissible with the relative monodromy filtrations induced by $M(\tau)$ because $\fg_\bR(\gr^W_w)=A(\gr^W_w)$. The adjoint action of $\sig$ on $\fg_\bR(\gr^W_0)$ is admissible with the relative monodromy filtrations induced by $M(\tau)$ by 1.2.2.2 because 
  $\fg_\bR(\gr^W_0)$ is a direct summand of $A(\gr^W_0)$ for the adjoint action of $\sig$. 
It follows that the adjoint action of $\sig$ on $\fg_\bR$ is admissible with the relative monodromy filtrations induced by $M(\tau)$.
\qed

\head
\S2.3.  The sets $E_\sig$ and $E_{\sig}^{\sharp}$
\endhead

Here we define the sets $E_\sig$ and $E_{\sig}^{\sharp}$. 

\medskip

{\bf 2.3.1.} {\it Toric varieties.}
\medskip

Assume that we are given $\Sig$ and $\G$ as in 2.2.6. 
  We assume that they are strongly compatible. 

  Fix $\sig \in \Sig$ in the following.  
  Let $P(\sig)=\Hom(\G(\sig), \bN)$.
Define
$$
\toric_\sig:= \Hom(P(\sig), \bC^{\mult}) \supset \abtoric_\sig:=\Hom(P(\sig), \bR_{\geq 0}^{\mult}).
$$
Here $\bC^{\mult}$ (resp.\ $\bR_{\geq 0}^{\mult}$) denotes the set $\bC$ (resp.\ $\bR_{\geq 0}$) regarded as a multiplicative monoid. 
Let
$$
\torus_\sig:= \G(\sig)^{\gp}\otimes_{\bZ} \bC^\times \supset \abtorus_\sig:=\G(\sig)^{\gp} \otimes_{\bZ} \bR_{>0}.
$$ 
We regard $\torus_\sig$ as an open set of $\toric_\sig$ and $\abtorus_\sig$ as an open set of $\abtoric_\sig$ via the embeddings 
$$
\torus_\sig=\Hom(P(\sig)^{\gp}, \bC^\times)\subset \toric_\sig, \quad \abtorus_\sig=\Hom(P(\sig)^{\gp}, \bR_{>0})\subset \abtoric_\sig.
$$ 
We have a natural action of $\torus_\sig$ on $\toric_\sig$, and a natural action of $\abtorus_\sig$ on $\abtoric_\sig$. 
We have an exact sequence 
$$
0\to \G(\sig)^{\gp}@>\log>> \sig_\bC@>\be>>  \text{torus}_\sig\to 0, 
$$ 
where 
$$
\be(\log(\g)\otimes z)=\g\otimes e^{2\pi iz}\;\; \text{for} \;\;\g\in 
\G(\sig), z\in \bC.
$$

For a face $\tau$ of $\sig$, the surjective homomorphism $P(\sig)\to P(\tau)$ induces an embedding $\toric_{\tau}\to \toric_{\sig}$. Let $0_\tau\in \toric_\tau\sub \toric_\sig$ be the point corresponding to the homomorphism $P(\tau)\to \bC^{\mult}$ which sends $1\in P(\tau)$ to $1$ and all the other elements of $P(\tau)$ to $0$. 
(This point $0_\tau\in \toric_\sig$ is denoted by $1_\tau$ in \cite{KU09}, but we change the notation.)

Any element $q$ of $\toric_\sig$ (resp. $|\toric|_\sig$) is written in the form $q=\be(a)\cdot 0_\tau$ (resp. $q=\be(ib)\cdot 0_\tau$) for $a\in \sig_\bC$ (resp. $b\in \sig_\bR$) and for a face $\tau$ of $\sig$. 
The face $\tau$ is determined by $q$, and $a$ modulo $\tau_\bC+ \log(\G(\sig)^{\gp})$ (resp. $b$ modulo $\tau_\bR$) is determined by $q$. 
\medskip

{\bf 2.3.2.} 
Let 
$$\Ec_\sig
:=\toric_\sig \times \Dc.$$

  We endow $\Ec_\sig$ with the inverse image of the log structure 
of $\toric_\sig$.

  We define a canonical pre-LMH on $\Ec_\sig$. 
  See \cite{KU09} for the pure case. 
  Let $(H_{\sig, \bZ}, W, (\langle\;,\;\rangle_w)_w)$ be the locally constant 
sheaf of $\bZ$-modules on $\toric_\sig^{\log}$ endowed with a weight  filtration and polarizations which is characterized by the property that the stalk at the unit point $1 \in \torus_{\sig} \sub\toric_\sig^{\log}$ is $(H_0, W, (\langle\;,\;\rangle_w)_w)$ and the action of $\pi_1(\toric_\sig^{\log})=\G(\sig)^{\gp}$ is given by $\G(\sig)^{\gp} \sub G_{\bZ}$. 
  We consider its pull-back to $\Ec_\sig$ still denoted by 
the same symbol $(H_{\sig, \bZ}, W, (\langle\;,\;\rangle_w)_w)$.
  Similarly as in the pure case, \cite{KU09} 2.3.7 gives a canonical 
isomorphism $\cO_{\Ec_{\sig}}^{\log}\otimes H_{\sig,\bZ} = 
\cO_{\Ec_{\sig}}^{\log}\otimes H_0$. 
  The $\cO_{\Dc} \otimes H_0$ has the universal Hodge filtration, 
which induces by pulling back a filtration on 
 $\cO_{\Ec_{\sig}}^{\log} \otimes H_0$, and by $\tau_*$ a filtration 
on $H_{\sig,\cO}:=\tau_*(\cO_{\Ec_{\sig}}^{\log} \otimes H_0)$. 
  Then, the triple $H_{\sig}:=(H_{\sig,\bZ}, W, H_{\sig,\cO})$ is a pre-LMH, 
which we call a canonical pre-LMH on $\Ec_\sig$. 
  Note that each graded piece of this pre-LMH 
together with $\langle\;,\;\rangle_w$ is regarded as a 
pre-PLH of weight $w$. 

\medskip

{\bf 2.3.3.} 
We define the subsets 
$\tilde E_\sig$, $E_\sig$, and 
$E^\sharp_\sig$ of $\Ec_\sig$.
  Let $H_{\sig}$ be the canonical pre-LMH on $\Ec_\sig$ defined in 2.3.2. 
  Then, the subset $\tilde E_\sig$ (resp.\ $E_\sig$) is defined as the set of 
points $x \in \Ec_\sig$ such that the pull-back of $H_{\sig}$ on $x$ 
satisfies the Griffiths transversality (resp.\ is an LMH with polarized 
graded quotients). 
  Thus, $E_\sig$ is contained in $\tilde E_\sig$. 

  Consider the subset $\abtoric_\sig$ of 
$\toric_\sig$, let $\Ec^{\sharp}_{\sig}=\abtoric_\sig \times \Dc$, 
and let $E^\sharp_\sig:=E_{\sig} \cap \Ec^{\sharp}_{\sig}$. 

$E_\sig$ and $E_\sig^\sharp$ are characterized as follows (cf.\ \cite{KNU10b} 2.2.3).

For $q \in \toric_\sig$ (resp. $q\in |\toric|_\sig$), write $q= \be(a)\cdot 0_\tau$ (resp. $q=\be(ib)\cdot 0_\tau$), where 
$a\in \sig_\bC$ ($b\in \sig_\bR$), $\tau$ is the face of $\sig$ and $0_\tau$ is as in 2.3.1.
Then, a point $(q, F)$ of $\Ec_\sig= \toric_\sig \times \Dc$ (resp. $\Ec^{\sharp}=|\toric|_\sig\times \Dc$) belongs to  $E_\sig$ (resp. $E_\sig^{\sharp}$) if and only if $(\tau, \exp(a)F)$ (resp. $(\tau, \exp(ib)F)$) generates  a nilpotent orbit.

The proofs of these facts are similar to the pure case (cf.\ \cite{KU09} 3.3.7). 

\medskip

{\bf 2.3.4.} 
In the above notation, we have canonical projections 
$\vf : E_\sig \to \G(\sig)^\gp\bs D_\sig$ and $\vf^\sharp : E^\sharp_\sig \to D^\sharp_\sig$ by
$$
\vf(q, F)=(\tau,\, \exp(\tau_\bC)\exp(a)F) \bmod \G(\sig)^\gp,
$$
$$
\vf^\sharp(q,F)=(\tau,\,\exp(i\tau_\bR)\exp(ib)F).
$$

\medskip

\head
\S2.4. Topology, complex structures, and log structures
\endhead

  Let $\Sig$ and $\G$ be as in 2.2.6.  
  We assume that they are strongly compatible. 

We define a structure of log local ringed space over $\bC$
on $\G\bs D_\Sig$ and a topology on $D_{\Sig}^{\sharp}$.

\medskip

{\bf 2.4.1.} 
As in 2.4.1 of \cite{KNU10b}, for each $\sig \in \Sig$, 
we endow the subsets  $E_\sig$ and $\tilde E_\sig$ of $\Ec_\sig$ in 2.3.3 with the following structures of log local ringed spaces over $\bC$.
The topology is the strong topology in $\Ec_\sig$.
The sheaf $\cO$ of rings and the log structure $M$ are the inverse images of $\cO$ and $M$ of $\Ec_\sig$, respectively.

We endow a topology on $E_{\sig}^{\sharp}$ as a subspace of $E_\sig$.

\medskip

{\bf 2.4.2.}  
We endow $\G\bs D_\Sig$ with the strongest topology for which
the maps $\pi_\sig:E_\sig \overset \vf \to \to\G(\sig)^\gp\bs D_\sig\to\G\bs D_\Sig$ are continuous for all $\sig\in\Sig$.
  Here $\vf$ is as in 2.3.4.
We endow $\G\bs D_\Sig$ with the following sheaf
of rings $\cO_{\G\bs D_\Sig}$ over $\bC$ and the
following log structure $M_{\G\bs D_\Sig}$.
For any open set $U$ of $\G\bs D_\Sig$ and for any
$\sig\in\Sig$, let $U_\sig:=\pi_\sig^{-1}(U)$ and define
$
\cO_{\G\bs D_\Sig}(U)\;
\text{(\resp. $M_{\G\bs D_\Sig}(U)$)}
:=\{\text{map $f:U\to\bC$}\;|\;
f\circ\pi_\sig\in\cO_{E_\sig}(U_\sig)\;
(\resp. \in M_{E_\sig}(U_\sig))\;
(\forall\sig\in\Sig)\}.
$ Here we regard $M_{E_\sig}(U_\sig)$ as a subset of $\cO_{E_\sig}(U_\sig)$ via the structural map $M_{E_\sig}(U_\sig)\to \cO_{E_\sig}(U_\sig)$, which is injective. 

\medskip

{\bf 2.4.3.} 
We introduce on $D_\Sig^\sharp$ the strongest topology for which the maps $E_\sig^\sharp \overset {\vf^\sharp} \to \to D_\sig^\sharp\to D_\Sig^\sharp$ $(\sig\in\Sig)$ are continuous.
Here $\vf^\sharp$ is as in 2.3.4. 
Note that the surjection $D_\Sig^\sharp \to\G\bs D_\Sig$ (cf.\ 2.2.5) becomes continuous.

\head
\S2.5. Properties of $\G\bs D_\Sig$
\endhead

In this \S2.5, let $\Sig$ be a weak fan in $\fg_\bQ$ and let $\G$ be a
subgroup of $G_\bZ$ which is strongly compatible
with $\Sig$ {\rm(2.2.6)}. 

The following Theorem 2.5.1--Theorem 2.5.6 are the mixed Hodge theoretic versions of 4.1.1 Theorem A (i)--(vi) of \cite{KU09}. 

\proclaim{Theorem 2.5.1} 
For $\sig\in\Sig$, $E_\sig$ is open in $\t E_\sig$ in the strong topology of $\t E_\sig$ in $\Ec_\sig$.
Both $\t E_\sig$ and $E_\sig$ are log manifolds.
\endproclaim

We say $\G$ is {\it neat} if for any $\g\in \G$, the subgroup of $\bC^\times$ generated by all eigen values of the action of $\g$ on $H_{0,\bC}$ is torsion free. It is known that there is a neat subgroup of $G_\bZ$ of finite index. 

\proclaim{Theorem 2.5.2} 
If $\G$ is neat, then $\G\bs D_\Sig$ is a log manifold.
\endproclaim

\proclaim{Theorem 2.5.3} 
Let $\sig\in\Sig$ and define the action of $\sig_\bC$ on $E_\sig$ over $\G(\sig)^\gp\bs D_\sig$ by
$$
a\cdot(q,F):=(\be(a)\cdot q,\exp(-a)F)\quad
(a\in\sig_\bC,\,(q,F)\in E_\sig),
$$
where $\be(a)\in \torus_\sig$ is as in  $2.3.1$ and $\be(a)\cdot q$ is
defined by the natural action of $\torus_\sig$ on $\toric_\sig$.
Then, $E_\sig\to\G(\sig)^\gp\bs D_\sig$ is a
$\sig_\bC$-torsor in the category of log manifolds.
That is, locally on the base $\G(\sig)^\gp\bs D_\sig$, $E_\sig$
is isomorphic as a log manifold to the product of $\sig_\bC$
and the base endowed with the evident action of $\sig_\bC$.

\endproclaim

\proclaim{Theorem 2.5.4} 
If $\G$ is neat, then, for any $\sig\in\Sig$, the map
$$
\G(\sig)^\gp\bs D_\sig\to\G\bs D_\Sig
$$
is  locally an isomorphism of log manifolds.
\endproclaim

\proclaim{Theorem 2.5.5} 
The topological space $\G\bs D_\Sig$ is Hausdorff.
\endproclaim

\proclaim{Theorem 2.5.6} 
If $\G$ is neat, then there is a homeomorphism of topological spaces
$$
(\G\bs D_\Sig)^\loga\simeq\G\bs D_\Sig^\sharp,
$$
which is compatible with $\tau:(\G\bs D_\Sig)^\loga\to\G\bs D_\Sig$ and
the projection $\G\bs D_\Sig^\sharp\to\G\bs D_\Sig$ induced by
$D_\Sig^\sharp\to D_\Sig$ in {\rm2.2.5}.
\endproclaim

The proofs of these theorems will be given in \S4--\S5.

\medskip

\head
\S2.6. Moduli  
\endhead

We define the moduli functor of LMH with polarized graded
quotients.

\medskip

{\bf 2.6.1.} 
Fix $\Phi= (\Lambda, \Sig, \G)$, where $\Lambda$ is as in 2.1.1 and $\Sig$ and $\G$ are as in 2.5.

  Let $S$ be an object of $\cB(\log)$. 
  Recall that an LMH with polarized graded quotients on $S$ is a 
quadruple $H=(H_\bZ, W, H_{\cO}, (\langle\;,\;\rangle_w)_w)$
such that $(H_\bZ, W, H_{\cO})$ is an LMH on $S$ and 
$(H(\gr^W_w)_{\bZ}, H(\gr^W_w)_{\cO}, \langle\;,\;\rangle_w)$ is a PLH 
of weight $w$ for any $w$. 

\medskip

{\bf 2.6.2.} 
Let $S$ be an object of $\cB(\log)$. 
By an {\it LMH with polarized graded quotients of type $\Phi$
on $S$}, we mean an LMH with polarized graded quotients  
$H=(H_\bZ, W, H_{\cO}, (\langle\;,\;\rangle_k)_k)$ 
endowed with a global section $\mu$ of the sheaf $\Gamma \bs\Isom((H_\bZ, W, (\langle\;,\;\rangle_k)_k), (H_0, W,
(\langle\;,\;\rangle_k)_k))$ on $S^\loga$ (called a {\it $\G$-level structure}) which satisfies the following conditions (1) and (2). 

\smallskip

(1) $\rank_\bZ(H_\bZ)=\sum_{p,q} h^{p,q}, \quad 
\rank_{\cO_S} (F^pH_\cO(\gr^W_{p+q}))/(F^{p+1}H_\cO(\gr^W_{p+q}))= h^{p, q}\quad \text{for all} \, p, q$.

\smallskip

(2) For any $s\in S$ and $t \in S^\loga$ lying over $s$, if ${\tilde \mu}_t :
(H_{\bZ,t}, W, (\langle\;,\;\rangle_k)_k)\overset{\sim}\to  \to (H_0, W,
(\langle\;,\;\rangle_k)_k)$ is a representative of the germ 
of $\mu$ at $t$, then there exists $\sig\in \Sig$ such that the image of the composite map
$$
\Hom(M_{S,s}/\cO_{S,s}^\times, \bN) \hookrightarrow \pi_1(s^{\log}) \to
\Aut(H_{\bZ,t}, W, (\langle\;,\;\rangle_k)_k)@>{\text{by} \;{\tilde \mu}_t}>> \Aut(H_0, W, (\langle\;,\;\rangle_k)_k)
$$
is contained in $\sig$ and such that the $\exp(\sig_{\bC})$-orbit $Z$ including $\tilde \mu_t(\bC\otimes_{{\cO}^{\log}_{S,t}}F_t )$, which is independent of the choice of a specialization ${\cO}^{\log}_{S,t} \to \bC$ at $t$ (1.3.1), is a $\sig$-nilpotent orbit. 

\medskip

We call an LMH with polarized graded quotients 
of type $\Phi$ on $S$ also an {\it LMH
with polarized graded quotients, with global monodromy in $\Gamma$, and with local monodromy in
$\Sig$}. 

\medskip

{\it Remark.}  There is an error in the  condition (1) in \cite{KNU10b} 3.2.2. 
The correct condition is the above (1). 

\medskip

{\bf 2.6.3.} Moduli functor $\LMH_{\Phi}$.

Let $\LMH_{\Phi} : \cB(\log) \to$ (set) be the contravariant functor defined 
as follows: For an object $S$ of $\cB(\log)$, 
$\LMH_{\Phi}(S)$ is the set of isomorphism classes of 
LMH with polarized graded quotients 
of type $\Phi$ on $S$. 

\medskip

{\bf 2.6.4.} Period maps.

We will have some period maps associated to an LMH with polarized graded quotients of type $\Phi$.
  In the following, assume that $\G$ is neat. 

\medskip

(1) $\LMH_\Phi \to \Map(-, \Gamma \bs D_\Sig)$. 

  Set theoretically, this period map is described as follows. 
  Let $S \in \cB(\log)$. 
  Let $H \in \LMH_{\Phi}(S)$. 
  Then, the image of $s \in S$ by the map corresponding to $H$ 
is $((\sig, Z) \bmod \G) \in \G\bs D_{\Sig}$. 
  Here $\sig$ is the smallest cone of $\Sig$ satisfying 2.6.2 (2), 
which exists by 2.2.4 $(2)''$, and $Z$ is the associated 
$\exp(\sig_{\bC})$-orbit which appeared in 2.6.2 (2). 

\medskip

(2) $\LMH_\Phi\to \Map(-^\loga, \Gamma \bs D_\Sig^{\sharp})$. 

  This is the composite of (1) and 
($S @>f>> \G \bs D_{\Sig})
\mapsto 
(S^\loga @> f^\loga>> (\G \bs D_{\Sig})^\loga \simeq \G \bs 
D^{\sharp}_{\Sig})$ (cf.\ 2.5.6). 
  Set theoretically, this is described as follows. 
  Let $S \in \cB(\log)$. 
  Let $H \in \LMH_{\Phi}(S)$. 
  Then, the image of $t \in S^{\log}$ by the map corresponding to $H$ 
is $((\sig, Z) \bmod \G) \in \G \bs D^{\sharp}_{\Sig}$. 
  Here $\sig$ is the same as in (1), and $Z$ is the 
$\exp(i\sig_{\bR})$-orbit including 
$\tilde \mu_t(\bC\otimes_{{\cO}^{\log}_{S,t}}F_t )$ 
which appeared in 2.6.2 (2). 

\medskip

(3) A variant of (2). Let  $\cB(\log)^{\log}$ be the category of pairs $(S, U)$, where $S$ is an object of $\cB(\log)$ and $U$ is an open set of $S^\loga$.
Let $\LMH_\Phi^\loga$ be the functor on $\cB(\log)^\loga$ defined as follows. 
$\LMH_\Phi^\loga((S,U))$ is the set of isomorphism classes of an 
$H\in \LMH_\Phi(S)$ plus a representative $\tilde \mu: (H_\bZ, W, (\lan\;,\;\ran_k)_k) \simeq (H_0, W, (\lan\;,\;\ran_k)_k)$ of $\mu$ given on $U$. 
Then the map $\LMH_\Phi^\loga((S,U)) \to \Map(U, D^{\sharp}_\Sig)$ is given, whose set-theoretical description is similar to (2). 

\medskip

{\bf 2.6.5.} 
The idea of the definition of the period map  
$$
\LMH_\Phi\to \Map(-, \Gamma\bs D_\Sig)
$$ 
of functors on $\cB(\log)$ in 2.6.4 (1) is the same as in \cite{KU09}: 
  Let $S \in \cB(\log)$. 
  Let $H \in \LMH_{\Phi}(S)$. 
  Then, locally on $S$, 
   the corresponding $S\to \G\bs D_{\Sig}$ 
   comes from a 
morphism $S\to E_\sig$ for some $\sig \in \Sig$.

\medskip

\proclaim{Theorem 2.6.6} 
Let $\Phi=(\Lambda,\Sig,\G)$ be as in $2.6.1$.
Assume that $\G$ is neat. 
Then $\text{\rom{LMH}}_\Phi$ is represented by $\G \bs D_\Sig$ in 
$\cB(\log)$. 
\endproclaim

\medskip

The proof will be given in \S4--\S5.

\bigskip

\head
\S3. Moduli spaces of log mixed Hodge structures with given graded quotients
\endhead

Let $S$ be an object of $\cB(\log)$, and assume that we are given a family
$Q=(H_{(w)})_{w\in \bZ}$, where $H_{(w)}$ is a pure LMH on $S$ 
of weight $w$ for each $w$ and $H_{(w)}=0$ for almost all $w$. Shortly speaking, the subject of \S3 is to construct spaces over $S$ which classify LMH $H$ with $H(\gr^W_w)=H_{(w)}$ for any $w$. 

We assume that all $H_{(w)}$ are polarizable, though we expect that the method of \S3 works (with some suitable modifications) without assuming the polarizability of $H_{(w)}$ (it actually works successfully in some non-polarizable case; cf. 7.1.4).

More precisely, let
$$
\LMH_Q
$$
be the contravariant functor on the category $\cB(\log)/S$ defined as follows. For an object $S'$ of $\cB(\log)$ over $S$, $\LMH_Q(S')$ is the set of all isomorphism classes of an LMH $H$ on $S'$ endowed with an isomorphism $H(\gr^W_w)\simeq H_{(w)}|_{S'}$ for each $w\in \bZ$. Here $H_{(w)}|_{S'}$ denotes the pull-back of $H_{(w)}$ on $S'$.
The functor  
 $\LMH_Q$ is usually not representable. Our subject is to construct relative log manifolds over $S$ which represent big subfunctors of $\LMH_Q$.
   
 We can often construct such a relative log manifold over $S$ by the method of \S2: As is explained in \cite{KNU10a}, we can construct such a 
space as the fiber product of
$S\to \tp_{w\in \bZ}\; \G_w\bs D(\gr^W_w)_{\sig_w}\leftarrow \G\bs D_{\Sig}$, where $S\to \G_w\bs D(\gr^W_w)_{\sig_w}$ is the period map of $H_{(w)}$ (if it exists) and $D_{\Sig}$ is a moduli space in \S2 (see \S3.4 for details). However, this method has a disadvantage that the period map $S\to \G_w\bs D(\gr^W_w)_{\sig_w}$ can be defined in general only after blowing up $S$ as is explained in \cite{KU09} \S4.3.  
(Cf.\ 3.4.3.)

The main purpose of this \S3 is to construct more relative log manifolds over $S$ which represent big subfunctors of $\LMH_Q$, by improving the method of \S2. The difference of \S2 and \S3 lies in that though cones considered in \S2 were inside $\fg_\bR$, the cones considered in \S3 exist outside $\fg_\bR$. 
  In particular, by the method of this \S3, 
  we will obtain in \S6 
  the N\'eron model 
 and the connected N\'eron model, an open subspace of the N\'eron model, 
 both of which represent big subfunctors of $\LMH_Q$ (see Theorem 6.1.1). In this method, we do not need blow up $S$.

\medskip

\head
\S3.1. Relative formulation of cones, the set $D_{S,\Sig}$ of nilpotent orbits
\endhead

Here we give a formulation of cones which is suitable for the study of the 
above functor $\LMH_Q$ of log mixed Hodge structures relative to the prescribed $Q$. 

\medskip

{\bf 3.1.1.} 
Assume that we are given an object $S$ of $\cB(\log)$  and a polarized log Hodge structure  $H_{(w)}$ on $S$  of weight $w$ for each $w\in \bZ$.
Assume that $H_{(w)}=0$ for almost all $w$. 
Write $Q=(H_{(w)})_{w \in \bZ}$.

Then, locally on $S$, we can find the following (1)--(5).

\smallskip

(1) $\Lam=(H_0, W, (\langle\;,\;\rangle_w)_w, (h^{p,q})_{p,q})$ as in 2.1.1.
\smallskip

(2) A sharp fs monoid $P$.
\smallskip

(3) A homomorphism 
$$
\G':=\Hom(P^{\gp},\bZ) \to G'_\bZ:=G_\bZ(\gr^W)=\tp_w G_\bZ(\gr^W_w)
$$
whose image consists of unipotent automorphisms.

\smallskip

((4) and (5) are stated after a preparation.) 
Let 
$$
\sig' := \Hom(P, \bR_{\geq 0}^{\add}), \quad
\toric_{\sig'}:=\Spec(\bC[P])_\an.
$$ 
Let 
$$\sig'_\bR \to \fg'_\bR:=\fg_\bR(\gr^W)$$
be the $\bR$-linear map induced by the logarithm of the homomorphism in (3).

We define the set $E'_{\sig'}$.  
The definition is parallel to the definition of $E_{\sig}$ in 2.3.3.
Let $E'_{\sig'}$ be the set of all points $s\in \toric_{\sig'}\times \Dc(\gr^W)$, where $\Dc(\gr^W)=\prod_w \Dc(\gr^W_w)$,  such that for each $w$, the pull-back of the canonical pre-PLH of weight $w$ on $\toric_{\sig'}\times \Dc(\gr^W_w)$ to the fs log point $s$ is a PLH. 
Here the canonical pre-PLH of weight $w$ on $\toric_{\sig'}\times \Dc(\gr^W_w)$ 
is defined in the same way as in 2.3.2.

In other words, 
$E'_{\sig'}$ is  the set of all pairs $(q', F')\in \toric_{\sig'}\times \Dc(\gr^W)$ satisfying the following condition:
 
If we write $q'=\be(a)\cdot 0_{\tau'}$ with $\tau'$ a face of $\sig'$ and $a\in \sig'_\bC$ (2.3.1), 
then $(\tau'_{\fg'}, \exp(a_{\fg'})F')$ generates a nilpotent orbit, where $(-)_{\fg'}$ denotes the image  under the $\bC$-linear map 
$\sig'_\bC \to \fg'_\bC:=\fg_\bC(\gr^W)$
induced by the above $\bR$-linear map $\sig'_\bR\to \fg'_\bR$. 

We endow $E'_{\sig'}$ with the structure of an object of $\cB(\log)$ by using the embedding $E'_{\sig'}\sub 
\toric_{\sig'} \times \Dc(\gr^W)$. 

\smallskip

(4) A strict morphism $S \to E'_{\sig'}$ in $\cB(\log)$. 

\medskip

Here a morphism $X\to Y$ in $\cB(\log)$ is said to be {\it strict} if via this morphism, the log structure of $X$ coincides with the inverse image of the log structure of $Y$. 

\medskip

(5) An isomorphism of PLH for each $w\in \bZ$ between  $H_{(w)}$ and the pull-back of the canonical PLH of weight $w$ on $E'_{\sig'}$ under the morphism in (4). 

\medskip

{\bf 3.1.2.} 
The local existence of (1)--(5) in 3.1.1 follows from a general theory of PLH in \cite{KU09}. 

In fact, take $s_0\in S$, and take a point  $t_0\in S^{\log}$ lying over $s_0$.
Let $H_0: =\bigoplus_w H_{(w),\bZ,t_0}$ with the evident weight filtration $W$ on $H_{0,\bR}$ and with the bilinear form $\langle\;,\;\rangle_w$ given by the polarization of $H_{(w)}$.  
For each $w$, let $(h^{p,q})_{p,q}$, $p+q=w$, be the Hodge type of $H_{(w)}$ at $s_0$.

Let $P:=(M_S/\cO^\times_S)_{s_0}$. Then $\G':=\Hom(P^{\gp},\bZ)$ is identified with the fundamental group $\pi_1(s_0^{\log})$, and the actions of $\pi_1(s_0^{\log})$ on $H_{(w),\bZ,t_0}$ for $w\in \bZ$ give the homomorphism $\G'\to G'_{\bZ}$. 

On an open neighborhood of $s_0$ in $S$, there is a homomorphism $P\to M_S$ which induces the identity map $P\to (M_S/\cO^\times_S)_{s_0}$ and which is a chart (\cite{KU09} 2.1.5) of the fs log structure $M_S$.  
We replace $S$ by this neighborhood of $s_0$. 
Then this homomorphism induces a strict morphism $S\to \toric_{\sig'}$. 
By \cite{KU09} 8.2.1--8.2.3, we obtain a strict morphism $S\to E'_{\sig'}$ over $\toric_{\sig'}$ and the isomorphism (5) in 3.1.1. 

\medskip

{\bf 3.1.3.} 
Assume that we are given $S$, $H_{(w)}$, $Q=(H_{(w)})_{w \in \bZ}$, and (1)--(5) of 3.1.1 and fix them. 

We will fix a splitting $H_0=\bigoplus_w H_0(\gr^W_w)$ of the filtration $W$ on $H_0$. We will denote $H_0(\gr^W_w)$ also by $H_{0,(w)}$. 

Let $\mu$ be the $\G'$-level structure on $H_\bZ(\gr^W)
:=\bigoplus_w H_{(w),\bZ}$, i.e., the class of local isomorphisms $H_\bZ(\gr^W)\simeq H_0$ on $S^{\log}$ which is global mod $\G'$, defined by the isomorphism (5) in 3.1.1.

We will denote $E'_{\sig'}$ often by $S_0$.

In \S3, the notions nilpotent cone, nilpotent orbit, fan, and weak fan are different from those in \S2. 
Recall $\sig'=\Hom(P, \bR_{\geq 0}^{\add})$ together with the linear map 
$\sig'_\bR \to \fg'_\bR:=\fg_{\bR}(\gr^W)$. In \S3,  nilpotent cones, fans, and weak fans appear in the fiber product $\sig'\times_{\fg_\bR'} \fg_\bR$, and nilpotent orbits appear in the product $\sig'_\bC\times \Dc$, as explained below.

\medskip

{\bf 3.1.4.} Nilpotent cone.

In \S3, 
a {\it nilpotent cone} means a
finitely generated cone $\sig$ in the fiber product 
$$
\sig' \times_{\fg'_\bR} \fg_\bR
$$ 
satisfying the following (1) and (2). 

\smallskip

(1) For any $N\in \sig$, the $\fg_\bR$-component $N_\fg\in \fg_\bR$ of $N$ is nilpotent. 

\smallskip

(2) For any $N_1, N_2\in \sig$, $N_{1,\fg}N_{2,\fg}=N_{2,\fg}N_{1,\fg}$.

\smallskip

We say a nilpotent cone $\sig$ is {\it admissible} if the following condition (3) is satisfied. 

\smallskip

(3) The action of $\sig$ on $H_{0,\bR}$ via $\sig\to\fg_\bR,\;N\mapsto N_\fg$, 
is admissible with respect to $W$ in the sense of 1.2.2. 

\medskip

{\bf 3.1.5.} Nilpotent orbit. 

Let $\sig$ be a nilpotent cone. A 
{\it $\sig$-nilpotent orbit} is a subset $Z$ of $\sig'_\bC \times \Dc$ satisfying the following conditions (1) and (2).

\smallskip

(1) If $(a,F)\in Z$ ($a\in \sig'_\bC, F\in \Dc$), then $Z=\{(a+b', \exp(b_{\fg})F)\;|\; b\in \sig_\bC\}$ (here $b'$ denotes the image of $b$ in $\sig'_\bC$).

\smallskip

(2) If $(a,F)\in Z$, then $(\sig_{\fg}, F)$ generates a nilpotent orbit.

\smallskip

Note that for $F\in \Dc$ and $a \in \sig'_\bC$, $\{(a+b', \exp(b_{\fg})F)\;|\;b\in \sig_\bC\}$ is a $\sig$-nilpotent orbit if and only if the above condition (2) is satisfied, which depends on $F$ but not on $a$. 
If this condition (2) is satisfied, 
we say $F$ {\it generates a $\sig$-nilpotent orbit} (or $(\sig, F)$, or further, $(N_1, \dots, N_n, F)$ {\it generates a nilpotent orbit} 
when $\sig$ is generated by $N_1,\ldots, N_n$).

\medskip

{\bf 3.1.6.} Fan and weak fan.

In \S3, a {\it fan} is a non-empty set of sharp rational nilpotent cones in $\sig'\times_{\fg'_\bR}\fg_\bR$ which is closed under the operations of taking a face and taking the intersection. 

A {\it weak fan} is a non-empty set $\Sig$ of sharp rational nilpotent cones in $\sig'\times_{\fg'_\bR}\fg_\bR$ which is closed under the operation of taking a face and satisfies the following condition. 

\smallskip

(1) Let $\sig, \tau \in \Sig$, and assume that $\sig$ and $\tau$ have a common interior point. 
Assume that there is an $F\in \Dc$ such that $(\sig, F)$ and $(\tau, F)$ generate nilpotent orbits in the sense of 3.1.5. 
Then $\sig=\tau$.  

\medskip

A fan is a weak fan. (This is proved in the same way as 2.2.4.)

An example of a fan is
$$\Xi= \{\bR_{\geq 0}N\;|\; N\in \sig'\times_{\fg'_\bR} \fg_\bR, N\;\text{is rational, $N_{\fg}$ is nilpotent}\}.$$

\medskip

{\bf 3.1.7.} Recall 
$$
\G' = \Hom(P^{\gp}, \bZ).
$$ 
We denote the given homomorphism $\G'\to G'_\bZ$ (3.1.1 (3)) by $a\mapsto a_{G'}$. 
We regard $G'_\bZ$ as a subgroup of $G_\bZ$ by the fixed splitting of $W$ of $H_0$. 

As is easily seen, there is a one-to-one correspondence between

\medskip

(1) A subgroup $\G_u$ of $G_{\bZ,u}$ such that $\g_{G'} \G_u\g_{G'}^{-1}=\G_u$ for all $\g\in \G'$. 

\medskip

(2) A subgroup $\G$ of $\G'\times_{G'_{\bZ}} G_\bZ$ 
containing $\G'=\{(a, a_{G'})\;|\; a\in \G'\}$. 
%such that the projection $\G\to \G'$ is surjective.

\medskip

The correspondence is:

From (1) to (2). $\G:=\{(a, a_{G'}b)\;|\; a\in \G', b \in \G_u\}$. 
  
From (2) to (1). $\G_u:=\Ker(\G\to \G')$.

In what follows, we will denote the projection $\G'\times_{G'_\bZ} G_\bZ\to G_\bZ$ by $a\mapsto a_{G}$. 

\medskip

{\bf 3.1.8.} 
Let $\Sig$ be a weak fan, and let $\G$ be a subgroup of $\G' \times_{G'_{\bZ}}
G_\bZ$ satisfying (2) in 3.1.7. 

We say that $\Sig$ and $\G$ are {\it compatible} if $\Ad(\g)(\sig)\in \Sig$ for any $\g\in \G$ and $\sig\in \Sig$. 
Here $\Ad(\g)$ is $(x,y)\mapsto (x, \Ad(\g_G)y)$ $(x\in \sig'$, $y\in \fg_\bR$).

We say $\Sig$ and $\G$ are {\it strongly compatible} 
if they are compatible and, for any $\sig\in \Sig$, any element of $\sig$ 
is a sum of elements of the form $a\log(\g)$ with 
$a \in \bR_{\geq 0}$ and $\g\in \G(\sig)$. 
Here $\log$ is the map $\G \to \G' \times \fg_\bR$, 
where $\G\to \G'$ is the projection and $\G \to \fg_\bR$ is $\g\mapsto \log(\g_G)$, and 
$\G(\sig)$ is defined by
$$
\G(\sig):=\{\g\in \G\;|\;\log(\g)\in \sig\}.
$$

\medskip

{\bf 3.1.9.} $\G$-level structure on $H$. 

Let $\G$ be as in 3.1.8. 
Let $\G_{G}$ be the image of $\G$ in $G_{\bZ}$. 

Let $S' \in \cB(\log)/S$. 
Let $H$ be an LMH on $S'$. 
Assume that an isomorphism between $H(\gr^W_w)$ and $H_{(w)}\vert_{S'}$ 
is given for each $w$. 

Then, there is a one-to-one correspondence between the following (1) and (2).

\medskip

(1) The class $(\bmod\ \G_{G})$ of a collection of local isomorphisms  $(H_\bZ, W)\simeq (H_0,W)$ which is compatible with the given $\G'$-level structure on $H_\bZ(\gr^W)$ and which is global $\bmod \G_G$.

\medskip

(2) The class $(\bmod\ \G_{G})$ of a collection of local splittings $H_\bZ(\gr^W) \simeq H_\bZ$ of the filtration $W$ on the local system $H_\bZ$ on $S^{\log}$ which is global mod $\G_{G}$.  

\medskip

Such a class is called a {\it $\G$-level structure} on $H_\bZ$. 

\medskip

{\bf 3.1.10.} 
Assume that $\Sig$ and $\G$ as in 3.1.8 are given.
Assume that $\Sig$ is strongly compatible with $\G$. 
Define
$$
D_{S, \Sig}= \{(s, \sig, Z)\;|\;s\in S, \sig\in \Sig, Z\sub \sig'_\bC\times \Dc\;\text{satisfying (1) and (2) below}\}.
$$
\smallskip

(1) $Z$ is a $\sig$-nilpotent orbit (3.1.5). 
\smallskip

(2) For any $(a, F) \in Z$, the image of $s$ under $S\to E'_{\sig'}$ coincides with $(\be(a)\cdot 0_\sig(\gr^W),\,  \exp(-a_{\fg'})F(\gr^W))$. 
Here $0_{\sig}(\gr^W)$ denotes the image of 
$0_{\sig}\in \toric_{\sig}$ in $\toric_{\sig'}$. 
(That is, if  $\alpha$ denotes the smallest face of $\sig'$ which contains the image of $\sig\to \sig'$, then $0_{\sig}(\gr^W)=0_{\alpha}$.)

\medskip Note that under the condition (1), 
the element $(\be(a)\cdot 0_\sig(\gr^W),\,  \exp(-a_{\fg'})F(\gr^W))$ of $E'_{\sig'}$ in the condition (2) is independent of the choice of $(a, F)\in Z$. 

If $'\Sig\sub \Sig$ denotes the subset of $\Sig$ consisting of all admissible  $\sig\in \Sig$, and if $''\Sig\sub {}'\Sig$ denotes the subset of $\Sig$ consisting of all $\sig\in \Sig$ such that a $\sig$-nilpotent orbit exists, then $'\Sig$ and $''\Sig$ are weak fans, and we have $D_{S,\Sig}=D_{S,'\Sig}=D_{S,''\Sig}$.  

\vskip20pt

\head
\S3.2. Toroidal partial compactifications $\G \bs D_{S,\Sig}$
\endhead

\medskip

Notation is as in \S3.1. 

Assume that 
$\Sig$ is a weak fan and is strongly compatible with $\G$. 
In this \S3.2, we endow  the set $\G \bs D_{S,\Sig}$ with a structure of an object of $\cB(\log)$. Here $\g\in \G$ acts on $D_{S,\Sig}$ as
$(s, \sig,Z)\mapsto (s, \Ad(\g)(\sig), \g Z)$, where $\g Z$ is defined by the action of $\g$ on $\sig'_\bC\times \Dc$ given by 
$(a,F)\mapsto (a-\g', \g_{G}F)$ with $\g'$ being 
the image of $\g$ in $\G'$ (3.1.7). 

\medskip

{\bf 3.2.1.} 
 Let $\sig \in \Sig$. 
Let
$$
\toric_\sig= \Spec(\bC[ \G(\sig){\spcheck}])_\an,
$$
where $\G(\sig){\spcheck}=\Hom(\G(\sig), \bN)$. 

Let $S$ be an object of $\cB(\log)$ and let
$$
E_{S, \sig}=\{(s,z, q, F)\;|\; s\in S, z\in \sig'_{\bC}, q\in \toric_{\sig}, F\in \Dc\;\text{satisfying the following (1) and (2)}\}.
$$

(1) If we write $q=\be(b)\cdot 0_\tau$ with $\tau$ a face of $\sig$ and $b\in \sig_\bC$ (2.3.1), then $\exp(b_{\fg})F$ generates a $\tau$-nilpotent orbit (3.1.5). 

\smallskip

(2) The image of $s$ under $S \to E'_{\sig'}$ coincides with $(\be(z)q(\gr^W), \exp(-z_{\fg'})F(\gr^W))$, 
where $q(\gr^W)$ is the image of $q$ in $\toric_{\sig'}$. 

\medskip

{\bf 3.2.2.} 
Note that, for $\sig\in \Sig$, we have a canonical map
$$
E_{S, \sig}\to \G(\sig)^{\gp}\bs D_{S,\sig},\;\;(s, z, q, F)\mapsto
\text{class}(s, \tau, Z).
$$
Here $(\tau, Z)$ is as follows.  
Write $q=\be(b)\cdot 0_{\tau}$ (2.3.1). 
This determines $\tau$.  
Further,  $Z$ is the unique $\tau$-nilpotent orbit containing $(z+b', \exp(b_{\fg})F)$. $D_{S,\sig}$ denotes $D_{S,\text{face}(\sig)}$. 

\proclaim{Proposition 3.2.3} 
Set-theoretically, $E_{S, \sig}$ is a $\sig_\bC$-torsor over $\G(\sig)^{\gp}\bs D_{S,\sig}$ with respect to the 
following action of $\sig_\bC$ on $E_{S, \sig}$. 
For $a\in \sig_\bC$, $a$ acts on $E_{S, \sig}$ as $(s,z, q, F)\mapsto (s,z-a', \;\be(a)q,\; \exp(-a_{\fg})F)$, where $a'$ denotes the image of $a$ in $\sig'_\bC$. 
\endproclaim

This is proved easily.

\medskip

{\bf 3.2.4.} 
Let $S_0=E'_{\sig'}$ as in 3.1.3. 
We consider the case that the morphism $S\to S_0$ is an isomorphism (we will express this situation as the case $S=S_0$). 
We endow $E_{S_0,\sig}$ with the structure of an object of $\cB(\log)$ by using the embedding $$E_{S_0,\sig}\sub \sig'_\bC \times \toric_\sig\times \Dc,\;\;(s,z,q, F)\mapsto (z,q, F).$$
That is, the topology of $E_{S_0,\sig}$ is the strong topology in $\sig'_\bC \times \toric_\sig\times \Dc$, and $\cO$ and the log structure $M$ of $E_{S_0,\sig}$ are the inverse images of $\cO$ and $M$ of $\sig'_\bC \times \toric_\sig\times \Dc$, respectively.
\medskip

{\bf 3.2.5.} 
We endow $\G\bs D_{S_0,\Sig}$ with the structure of a log local ringed space over $\bC$, as follows.

We define the topology of $\G \bs D_{S_0,\Sig}$ as follows. 
A subset $U$ of $\G\bs D_{S_0,\Sig}$ is open if and only if, for any $\sig\in \Sig$, the inverse image $U_\sig$ of $U$ in $E_{S_0,\sig}$ is open.

We define the sheaf of rings $\cO$ on $\G \bs D_{S_0,\Sig}$, as a subsheaf of the sheaf of $\bC$-valued functions,  as follows. 
For an open set $U$ of $\G \bs D_{S_0,\Sig}$ and for a map $f: U\to \bC$, $f$ belongs to $\cO$  if and only if, for any $\sig \in \Sig$, the pull-back of $f$ on the inverse image of $U$ in $E_{S_0,\sig}$ belongs to the $\cO$ of $E_{S_0,\sig}$.

We define the log structure of $\G\bs D_{S_0, \Sig}$ as a subsheaf of  $\cO$ of $\G\bs D_{S_0,\Sig}$, as follows. 
For an open set $U$ of $\G\bs D_{S_0, \Sig}$ and for an $f\in\cO(U)$, $f$ belongs to the log structure if and only if, for any $\sig \in \Sig$, the pull-back of $f$ on the inverse image of $U$ in $E_{S_0,\sig}$ belongs to the log structure of $E_{S_0,\sig}$. Here we regard the log structure of $E_{S_0,\sig}$ as a subsheaf of $\cO$ of $E_{S_0,\sig}$.
\medskip

\proclaim{Proposition 3.2.6}  {\rm (i)}
$E_{S_0,\sig}$ and $\G\bs D_{S_0,\Sig}$  are objects of $\cB(\log)$. 

\medskip

{\rm (ii)} $E_{S_0,\sig}$ is a $\sig_\bC$-torsor over $\G(\sig)^{\gp}\bs D_{S_0,\sig}$ in $\cB(\log)$. 
\endproclaim

This will be proved in \S5. 

\medskip

{\bf 3.2.7.} 
Note that, as a set, $\G \bs D_{S,\Sig}$ is the fiber product of 
$S\to S_0\leftarrow \G\bs D_{S_0,\Sig}$. 

We endow  $\G\bs D_{S,\Sig}$ with the structure as an object of $\cB(\log)$ by regarding it as the fiber product of $S\to S_0\leftarrow \G\bs D_{S_0,\Sig}$
in $\cB(\log)$. 
The following theorem will be proved in \S5. 

\proclaim{Theorem 3.2.8} 
{\rm(i)}  $\G\bs D_{S,\Sig}$ is a relative log manifold over $S$. 
 
 \medskip
 
{\rm (ii)} For $\sig\in \Sig$, the morphism $\G(\sig)^{\gp}\bs D_{S,\sig}\to \G\bs D_{S,\Sig}$ is locally an isomorphism in $\cB(\log)$.
 
  \medskip

{\rm (iii)} If $S$ is Hausdorff, $\G \bs D_{S, \Sig}$ is Hausdorff. 
 
 \endproclaim

{\bf 3.2.9.} 
We define $D_{S,\Sig}^{\sharp}$ in the evident way, replacing nilpotent orbit by nilpotent $i$-orbit which is defined in the evident way. 
Then $\G\bs D_{S,\Sig}^{\sharp}$ is identified with $(\G \bs D_{S, \Sig})^{\log}$ via the natural map.

\vskip20pt

\head
\S3.3. Moduli 
\endhead

\medskip

{\bf 3.3.1.} 
Let $\Sig$ be a weak fan in $\sig' \times_{\fg'_\bR} \fg_\bR$, where $\sig' = \Hom(P,\bR_{\geq 0}^{\add})$. 

Assume that we are given a subgroup $\G_u$ of $G_{\bZ,u}$ which satisfies the condition in 3.1.7 (1), and let $\G$ be the corresponding group as in 3.1.7.
Assume that $\Sig$ and $\G$ are strongly compatible. 
\medskip

\medskip

{\bf 3.3.2.} 
Let $\LMH_{Q,\G}^{(\Sig)}$ be the following functor
on the category $\cB(\log)/S$. 
For any object $S'$ of $\cB(\log)/S$, $\LMH_{Q,\G}^{(\Sig)}(S')$ is the set of all isomorphism classes of the following $H$.

\smallskip

$H$ is an LMH on $S'$ endowed with a $\G$-level structure $\mu$, whose $\gr^W$ is identified with  $(H_{(w)})_w$ endowed with the given $\G'$-level structure (3.1.3), 
satisfying the following condition. 

\smallskip

For any $s\in S'$, any $t\in (S')^{\log}$ lying over $s$,  any representative $\mu_t : H_{\bZ,t}\simeq H_0$ of $\mu$, and for any specialization $a$ at $t$ (1.3.1), there exists $\sig\in \Sig$ such that the image  of $\Hom((M_{S'}/\cO_{S'}^\times)_s, \bN)\to \sig' \times \fg_\bR$ is contained in $\sig$ and $\mu_t(F(a))$ generates a $\sig$-nilpotent orbit (3.1.5). 
Here  $F$ denotes the Hodge filtration of $H$. 

\medskip

We have the canonical morphism $\LMH_{Q,\G}^{(\Sig)}\to \LMH_Q$.

\medskip

In the case $\G_u=G_{\bZ,u}$, $\LMH_{Q,\G}^{(\Sig)}$ is the following functor.

For any object $S'$ of $\cB(\log)/S$, $\LMH_{Q,\G}^{(\Sig)}(S')$ is the set of all isomorphism classes of the following $H$.

\smallskip

$H$ is an LMH on $S'$, whose $\gr^W$ is identified with $(H_{(w)})_w$, satisfying the following condition. 

\smallskip

For any $s\in S'$, any $t\in (S')^{\log}$ lying over $s$,  any isomorphism $\mu_t : H_{\bZ,t}\simeq H_0$ whose $\gr^W$ belongs to the given $\G'$-level structure (3.1.3), and for any specialization $a$ at $t$ (1.3.1), there exists $\sig \in \Sig$ such that the image  of $\Hom((M_{S'}/\cO_{S'}^\times)_s, \bN)\to \sig'\times \fg_\bR$ is contained in  $\sig$ and  $\mu_t(F(a))$ generates a $\sig$-nilpotent orbit (3.1.5). 
Here $F$ denotes the Hodge filtration of $H$. 

\medskip

From this, we see that in the case $\G_u=G_{\bZ,u}$, $\LMH_{Q,\G}^{(\Sig)}\to \LMH_Q$ is injective and hence $\LMH_{Q,\G}^{(\Sig)}$ is regarded as a subfunctor of $\LMH_Q$. We will denote this subfunctor by $\LMH_Q^{(\Sig)}$. 
(This injectivity is explained also from the fact that, in the case $\G_u=G_{\bZ,u}$, $H$ with the given $\gr^W=(H_{(w)})_w$ is automatically endowed with a $\G$-level structure.)

In the case $\G_u=G_{\bZ,u}$, we denote $\G\bs D_{S,\Sig}$ by
$$
J_{S,\Sig}, \quad \text{or simply by} \;\;J_{\Sig}.
$$

\proclaim{Theorem 3.3.3} 
As an object of $\cB(\log)$ over $S$, $\G\bs D_{S,\Sig}$  represents {\rm{LMH}}${}_{Q,\G}^{(\Sig)}$.
\endproclaim

This will be proved in \S5. 

\bigskip

\head
\S3.4. Relation with \S2
\endhead

{\bf 3.4.1.}  
Let the situation be as in \S2. 

Let $\Sig$ be a weak fan, and let $\G$ be a neat subgroup of $G_{\bZ}$ which is strongly compatible with $\Sig$. 

  Let $\sig_w$ be a sharp rational nilpotent cone in  $\fg_\bR(\gr^W_w)$ for each $w$. 
Assume that, for any $\sig\in \Sig$ and for any $w$, the image of $\sig$ in $\fg_\bR(\gr^W_w)$ is contained in  $\sig_w$.

  Let $\G_w$ be a neat subgroup of $G_{\bZ}(\gr^W_w)$ for each $w$. Assume that for each $w$, the image of $\G$ in $G_\bZ(\gr^W_w)$ is contained in $\G_w$.  
  
Assume that the fan  $\text{face}(\sig_w)$ and $\G_w$ are strongly compatible. 

Let $S$ be an object of $\cB(\log)$. 
Assume that we are given a morphism $S \to \G_w\bs D(\gr^W_w)_{\sig_w}$ for each $w$. 
As in the papers \cite{KNU10a} and \cite{KNU10c}, we consider the fiber product of
$$
S \to \tp_w \G_w\bs D(\gr^W_w)_{\sig_w}\leftarrow \G\bs D_{\Sig}.
$$

Assume that $S\to \prod_w \G_w\bs D(\gr^W_w)_{\sig_w}$ is strict. Assume also that the homomorphism $\G\to \prod_w \G_w$ is surjective.  Then we have the following relation of this fiber product  to the space constructed in this section.

Let $H_{(w)}$ be the pull-back on $S$ of the universal polarized log Hodge structure of weight $w$ on $\G_w\bs D(\gr^W_w)_{\sig_w}$, and let $Q=(H_{(w)})_w$.

Locally on $S$, the morphism $S\to \G_w\bs D(\gr^W_w)_{\sig_w}$  factors as 
$S\to E(\gr^W_w)_{\sig_w} \to \G_w\bs D(\gr^W_w)_{\sig_w}$. We have the situation of 3.1.1 by taking $H_{0,(w)}= H_0(\gr^W_w)$ and $P$ the dual fs monoid of $\prod_w \G_w(\sig_w)$. 
We have $\sig'=\prod_w \sig_w$, $\G'=\prod_w \G_w(\sig_w)^{\gp}$,  and $\Sig$ can be regarded as a weak fan in the sense of 3.1.6 by identifying $\sig \in \Sig$ with $\{(N_{\sig'}, N)\;|\; N \in \sig\}$, where $N_{\sig'}$ is the image of $N$ in $\sig'$.

 Let $\tilde \G\sub \G$ be the inverse image of $\G'\sub \prod_w \G_w$. 

\proclaim{Proposition 3.4.2} 
We have a canonical isomorphism 
$$
\tilde \G\bs D_{S,\Sig}\simeq S \times_{\prod_w \G_w\bs D(\gr^W_w)_{\sig_w}} \G\bs D_{\Sig}
$$
in the category $\cB(\log)/S$.
\endproclaim
 
{\it Proof.} 
We have a canonical map $D_{S,\Sig}\to D_{\Sig}$ which sends $(s, \sig, Z)$ to $(\sig, \exp(\sig_\bC)F)$, where $F$ is any element of $\Dc$ such that $(a,F)\in Z$ for some $a\in \sig'_\bC$. 

The induced map 
 $\tilde \G\bs D_{S, \Sig} \to Y:=S \times_{\prod_w \G_w\bs D(\gr^W_w)_{\sig_w} } \G\bs D_{\Sig}$ is bijective. In fact, it is easy to see the surjectivity, and the injectivity is reduced to the following Claim 1. 

\medskip

{\bf Claim 1.} Assume
 $(s, \sig, Z_j) \in D_{S,\Sig}$ and $(a_j, F)\in Z_j$ for $j=1,2$, with $s\in S$, $\sig\in \Sig$,  
 $Z_j\sub \sig'_\bC\times \Dc$, 
  $a_j\in \sig'_\bC=\prod_w \sig_{w,\bC}$, $F\in \Dc$. 
 Then $a_1=a_2$ (hence $Z_1=Z_2$).
 
 \medskip
 
  We prove Claim 1. By the definition of $D_{S,\Sig}$, for $j=1,2$, the image of $s$ in $S_0=E'_{\sig'}$ is 
 $(\be(a_j)\cdot 0_\sig(\gr^W), \,\exp(-a_j)F(\gr^W))$. Hence $\exp(-a_1)F(\gr^W)=\exp(-a_2)F(\gr^W)$. Since $(\sig', F(\gr^W))$ generates a nilpotent orbit, 
 this implies $a_1=a_2$ by \cite{KU09} Proposition 7.2.9 (i). 
 
 We will prove
 
 \medskip
 
{\bf Claim 2.} 
In the case $S= S_0 = 
\prod_w E(\gr^W_w)_{\sig_w}$, the above bijection ${\tilde \G}\bs D_{S,\Sig}\to Y$ is in fact an isomorphism $\tilde \G \bs D_{S_0,\Sig}@>\sim>> Y$ in $\cB(\log)$.

\medskip By 
the base change by $S\to S_0$, the isomorphism in Claim 2 induces an isomorphism  $\tilde \G \bs D_{S,\Sig}@>\sim>> Y$ in $\cB(\log)/S$ in general.

We prove Claim 2. By Theorem 2.5.4 and Theorem 3.2.8 (ii), it is sufficient to prove that for each $\sig\in \Sig$, the morphism of $\cB(\log)$ from 
$\tilde \G(\sig)^{\gp}\bs D_{S_0,\sig}$ to the fiber product of
$$S_0\to \prod_w \G_w(\sig_w)^{\gp}\bs D(\gr^W_w)_{\sig_w}\leftarrow \G(\sig)^{\gp}\bs D_{\sig}$$
is locally an isomorphism. By Theorem 2.5.3, 
$S_0=E'_{\sig'}\to 
\prod_w \G_w(\sig_w)^{\gp}\bs D(\gr^W_w)_{\sig_w}$ is a $\sig'_\bC$-torsor and $E_\sig$ in \S2 is a $\sig_\bC$-torsor over $\G(\sig)^{\gp}\bs D_{\sig}$. 
Furthermore, $E_{S_0,\sig}$ is a $\sig_\bC$-torsor over $\tilde \G(\sig)^{\gp}\bs D_{S_0,\sig}$ by Proposition 3.2.6 (ii). 
Hence 
Claim 2 is reduced to the evident fact that the canonical morphism $E_{S_0,\sig}\to \sig'_\bC\times E_{\sig}$ is an isomorphism. 
\qed
 
\medskip
 
{\bf 3.4.3.}  
Thus, to construct the moduli space of LMH with a given graded quotients of the weight filtration, there are two methods. 
One is the method in \cite{KNU10a} and \cite{KNU10c} to take the fiber product as in the above 3.4.1, and the other is the method of this \S3. 
As advantages of the method of  \S3, we have:
\medskip

(1) As is shown in \S6, the N\'eron model and the connected N\'eron model 
can be constructed by the method of \S3 even in the case of higher dimensional base $S$. In the method of \S2, to obtain such models, 
the period map 
$S\to \G'\bs D_{\Sig'}$ of the $\gr^W$ is necessary to define the fiber 
product. 
  But, since such a period map of the $\gr^W$ 
can exist in general only after some blowing-up (see \S4.3 of \cite{KU09}), 
we may have to blow up the base $S$ in the method of \S2.

\medskip
(2) The case with no polarization on $\gr^W$ can be sometimes treated. (See 7.1.4  for example.)

\bigskip

\head
\S4. Associated $\SL(2)$-orbits
\endhead

In the next \S5, we will prove main results in \S2--\S3.  
This \S4 is a preparation for \S5. 
  We review some necessary definitions and results on the space 
$D_{\SL(2)}$ of $\SL(2)$-orbits in \cite{KNU.p}, introduce 
the most important map (called the CKS map) $D^{\sharp}_{\Sig,\val} \to 
D_{\SL(2)}$ in the fundamental diagram in Introduction, and prove that 
it is continuous, which will be crucial in \S5.

\medskip

Let the situation and terminologies (weak fan, nilpotent cone, etc.) be as in \S2 (not as in \S3).

\head
\S4.1.  $D_{\SL(2)}$ 
\endhead

  We review from \cite{KNU.p} 
the definition of SL(2)-orbit, the space $D_{\SL(2)}$ 
and the $\SL(2)$-orbit associated to an $(n+1)$-tuple 
$(N_1, \ldots, N_n, F)$ which generates a nilpotent orbit. 

\medskip

{\bf 4.1.1.} $\SL(2)$-orbit (\cite{KNU.p} 2.3).

Assume first that we are in the pure case (originally
considered by [CKS86]), that is, 
$W_w=H_{0,\bR}$ and $W_{w-1}=0$ for some $w$. Then an {\it $\SL(2)$-orbit in $n$ variables} is a pair 
$(\rho, \vf)$, where $\rho: \SL(2,\bC)^n\to G_\bC$ is 
a homomorphism of algebraic groups defined over 
$\bR$ and $\vf: {\bold P}^1(\bC)^n\to \Dc$ is a holomorphic map, satisfying the following three conditions (\cite{KU02}, \cite{KU09}).
\smallskip

(1) $\vf(gz)=\rho(g)\vf(z)\quad \text{for any}\;g\in \SL(2,\bC)^n,
\;z\in \bP^1(\bC)^n$.
\smallskip

(2) $\vf(\fh^n)\sub D$.
\smallskip

(3) $\rho_*(\fil^p_z (\fsl(2,\bC)^{\op n})) \sub \fil^p_{\vf(z)}(\fg_{\bC}) \quad \text{for any $z \in \bP^1(\bC)^n$ and any $p\in \bZ$}$.
\smallskip

\noindent
Here in (2), $\fh \subset \bP^1(\bC)$ is the upper-half plane.
In (3), $\rho_*$ denotes the Lie algebra homomorphism induced by $\rho$,
and $\fil_z(\fsl(2,\bC)^{\op n}))$, $\fil_{\vf(z)}(\fg_{\bC})$ are the filtrations
induced by the Hodge filtrations at $z$, $\vf(z)$, respectively.

Now we consider the general mixed Hodge situation. A {\it non-degenerate $\SL(2)$-orbit of rank $n$} is a 
  pair $((\rho_w, \vf_w)_{w \in \bZ}, \br)$, where $(\rho_w, \vf_w)$ is 
an $\SL(2)$-orbit in $n$ variables for $\gr^W_w$ (that is,  an $\SL(2)$-orbit for the pure case) for each $w \in \bZ$ and $\br$ is an element of $D$ 
 satisfying the following conditions (4)--(6).

\smallskip

(4) $\br(\gr^W_w)=\vf_w(\bold i)$ for any $w\in \bZ$, where $\bi=(i, \dots, i)\in {\bold P}^1(\bC)^n$. 
\smallskip

(5) If $2\leq j \leq n$, there exists $w\in \bZ$ such that the $j$-th component of $\rho_w$ is a non-trivial homomorphism.

\smallskip

(6) If $\br \in D_{\spl}$ and $n \geq 1$, there exists $w\in \bZ$ such that the 1-st component of $\rho_w$ is a non-trivial homomorphism. Here, as in \cite{KNU09},  
$D_{\spl}$ denotes the subset of $D$ consisting of $\bR$-split mixed Hodge structures.

\smallskip

  Let $((\rho_w, \vf_w)_w, \br)$ be a non-degenerate $\SL(2)$-orbit 
of rank $n$. 

Then, the associated homomorphism of algebraic groups over $\bR$
$$
\tau : \bG_{m,\bR}^n \to \Aut_{\bR}(H_{0, \bR}, W)
$$
and the associated set (or family) of weight filtrations 
are defined as follows. 

Recall that the canonical splitting (over $\bR$) of the weight filtration of a mixed Hodge structure is defined as in \cite{KNU08} \S1, as is reviewed in \cite{KNU09} \S4,  \cite{KNU.p} \S1.2. It was originally defined in \cite{CKS86}. By associating the canonical splitting $\spl_W(F)$ of $W$ to a mixed Hodge structure $F\in D$, we have a continuous map 
$$\spl_W:D\to \spl(W),$$ 
where $\spl(W)$ denotes the set of all splittings of $W$. 

Let $s_\br:=\spl_W(\br): \gr^W \overset \sim \to \to H_{0,\bR}$ be the canonical splitting of $W$ associated to $\br$.
Then  
$$
\tau(t_1,\dots, t_n)= s_\br\circ\Big(\tsize\bigoplus_{w\in \bZ}\Big(\tsize\prod_{j=1}^n t_j\Big)^w\rho_w(g_1, \dots, g_n)\;\text{on}\;\gr^W_w\Big)\circ s_\br^{-1}
$$ 
$$
\text{with}\quad g_j=\pmatrix 1/ \tsize\prod_{k=j}^n t_k& 0 \\
0 & \tsize\prod_{k=j}^n t_k\endpmatrix.
$$
For $1\leq j\leq n$, let $\tau_j: \bG_{m,\bR} \to \Aut_{\bR}(H_{0,\bR}, W)$ be the $j$-th component of $\tau$. 

For $1\leq j\leq n$, let $W^{(j)}$ be the increasing filtration on $H_{0,\bR}$ defined by
$$W^{(j)}_w=\tsize\bigoplus_{k\leq w} \;\{v\in H_{0,\bR}\;|\; \tau_j(a)v=a^k\;(\forall\; a\in \bR^\times)\}.$$
The associated set (resp. family) of weight filtrations is defined to be the set
$\{W^{(1)}, \cdots, W^{(n)}\}$ (resp. the family $(W^{(j)})_{1\leq j\leq n}$). 

\medskip

{\bf 4.1.2.} $D_{\SL(2)}$ (\cite{KNU.p} 2.5).

Two non-degenerate $\SL(2)$-orbits $p=((\rho_w, \vf_w)_w, \br)$ and $p'=
((\rho'_w, \vf'_w)_w, \br')$ of rank $n$ are said to be {\it equivalent} if there is a $t \in \bR^n_{>0}$ such that
$$
\rho'_w=\Int(\gr^W_w(\tau(t))) \circ \rho_w, \quad \vf'_w=\gr^W_w(\tau(t))\circ \vf_w\quad(\forall\; w\in \bZ), \quad \br'=\tau(t)\br.
$$
This is actually an equivalence relation.

The associated homomorphism of algebraic groups 
$
\tau : \bG_{m,\bR}^n \to \Aut_{\bR}(H_{0, \bR}, W),
$
and the associated set (resp. family) of weight filtrations depend only on the equivalence class.

\smallskip

Let $D_{\SL(2)}$ be the set of all equivalence classes $p$ of non-degenerate 
$\SL(2)$-orbits of various ranks satisfying the following condition (1).  

\smallskip

\noindent
(1) If $W'$ is a member of set of weight filtrations associated to $p$, then the $\bR$-spaces $W'_k\gr^W_w$ are rational (i.e. defined over $\bQ$)  for any $k$ and $w$.

\medskip

In \cite{KNU.p}, we defined two topologies of $D_{\SL(2)}$, which we call the stronger topology and the weaker topology, respectively, and denoted by $D_{\SL(2)}^I$ and $D_{\SL(2)}^{II}$ the set $D_{\SL(2)}$ endowed with the stronger topology and the weaker topology, respectively. The identity map of sets $D_{\SL(2)}^I\to D_{\SL(2)}^{II}$ is continuous in these topologies. In \cite{KNU.p}, we defined some kind of real analytic structures on $D_{\SL(2)}^I$ and on $D_{\SL(2)}^{II}$, but we do not use them in the present paper.

A finite set $\Psi$ of increasing filtrations on $H_{0,\bR}$ is called an admissible set of weight filtrations if $\Psi$ is the associated set of weight filtrations of some $p\in D_{\SL(2)}$. For an admissible set $\Psi$ of weight filtrations, let $D_{\SL(2)}^I(\Psi)$  be the subset of $D_{\SL(2)}$ consisting of all points $p$ of $D_{\SL(2)}$ such that the set of weight filtrations associated to $p$ is a subset of $\Psi$. Then $\{D_{\SL(2)}^I(\Psi)\}_{\Psi}$ is an open covering of $D_{\SL(2)}^I$.

\medskip

{\bf 4.1.3.} Associated $\SL(2)$-orbits (\cite{KNU.p} 2.4.2).
  
  Let $(N_1, \dots, N_n, F)$ generate a nilpotent orbit (2.2.2). 
  Then, we associate to it a non-degenerate $\SL(2)$-orbit 
$((\rho'_w, \vf'_w)_w, \br_1)$ as follows. 

First, by \cite{CKS86}, for each $w\in \bZ$, we have the $\SL(2)$-orbit $(\rho_w, \varphi_w)$ in $n$ variables for $\gr^W_w$ associated to $(\gr^W_w(N_1), \dots, \gr^W_w(N_n), F(\gr^W_w))$, which generates a nilpotent orbit for $\gr^W_w$. 
Let $k=\min(\{j\;|\;1\leq j \leq n,\; N_j\not=0\} \cup \{n+1\})$. 
Let 
$$
J'=\{j\;|\;1\leq j\leq n, \;\text{the $j$-th component of $\rho_w$ is non-trivial for some $w\in \bZ$}\}.
$$ 
Let $J=J'=\emptyset$ if $k=n+1$, and let $J=J'\cup\{k\}$ if otherwise.
  Let $J=\{a(1), \dots, a(r)\}$ with $a(1)<\dots < a(r)$. 
Then $(\rho_w', \vf_w')$ is an $\SL(2)$-orbit on $\gr^W_w$ characterized by 
$$
\rho'_w(g_{a(1)}, \dots, g_{a(r)}):=\rho_w(g_1, \dots, g_n), 
\quad 
\vf'_w(z_{a(1)}, \dots, z_{a(r)}):=\vf_w(z_1, \dots, z_n).
$$

  Next, if $y_j\in \bR_{>0}$ and $y_j/y_{j+1}\to \infty$ $(1\leq j\leq n$, $y_{n+1}$ means $1)$, the canonical splitting $\spl_W(\exp(\ts_{j=1}^n iy_jN_j)F)$ of $W$ (4.1.1) associated to $\exp(\ts_{j=1}^n iy_jN_j)F$ converges in $\spl(W)$. 
Let $s\in \spl(W)$ be the limit. 

Let $\tau: \bG_{m,\bR}^n \to \Aut_\bR(H_{0,\bR}, W)$ be the homomorphism of algebraic groups defined by 
$$
\tau(t_1,\dots, t_n)=s\circ \Big(\tsize\bigoplus_{w\in \bZ}\big(\big(\tsize\prod_{j=1}^n t_j\big)^w\rho_w(g_1, \dots, g_n) \;\text{on}\; \gr^W_w \big)\Big)\circ s^{-1}, 
$$
where $g_j$ is as in $4.1.1$.
Then, as $y_j>0$, $y_1=\cdots=y_k$, 
$y_j/y_{j+1}\to \infty $ $(k\leq j\leq n$, $y_{n+1}$ means $1)$,
$$
\tau\left(\sqrt{\frac{y_2}{y_1}},\dots, \sqrt{\frac{y_{n+1}}{y_n}}\right)^{-1}\exp(\ts_{j=1}^n iy_jN_j)F 
$$
converges in $D$. 
Let $\br_1\in D$ be the limit (cf.\ \cite{KNU.p} 2.4.2 (ii)). 

\medskip

\head
\S4.2. Valuative spaces
\endhead

{\bf 4.2.1.} 
$D_{\val}$. 

Let $D_{\val}$ (resp. $D_{\val}^{\sharp}$) be the set of all triples $(A, V, Z)$, where

\smallskip

$A$ is a $\bQ$-subspace of $\fg_\bQ$ consisting of mutually commutative nilpotent elements. 
Let $A_\bR = \bR\otimes_\bQ A$ and $A_\bC = \bC\otimes_\bQ A$.

\smallskip

$V$ is an (additive)  submonoid of $A^*=\Hom_\bQ(A, \bQ)$ satisfying the following conditions: $V\cap (-V)=\{0\}$, $V\cup (-V)=A^*$.

\smallskip

$Z$ is a subset of $\Dc$ such that $Z=\exp(A_\bC)F$ (resp.\ $Z=\exp(iA_\bR)F$) for any $F\in Z$ and that there exists a finitely generated rational subcone $\tau$ of $A_\bR$ satisfying the following condition (1).

\smallskip

(1) $Z$ is a $\tau$-nilpotent orbit (resp. $i$-orbit), and 
$V$ contains any element $h$ of $A^*$ such that $h: A \to \bQ$ sends $A\cap \tau$ to $\bQ_{\geq 0}$.
\smallskip

We have a canonical map
$$
D_{\val}^{\sharp}\to D_{\val},\quad
(A,V,Z) \mapsto (A,V, \exp(A_\bC)Z),
$$ 
which is surjective.
\medskip

{\bf 4.2.2.} 
$D_{\Sig,\val}$ and $D_{\Sig,\val}^{\sharp}$. 

Let $\Sig$ be a weak fan (2.2.3) in $\fg_\bQ$. 

Let $D_{\Sig,\val}$ (resp.\ $D_{\Sig,\val}^{\sharp}$) be the subset of $D_{\val}$ (resp.\ $D_{\val}^{\sharp}$) consisting of all $(A, V, Z) \in D_{\val}$ (resp.\ $D_{\val}^{\sharp}$) such that there exists $\sig\in \Sig$ satisfying the following conditions (1) and (2). 

\smallskip

(1) $\exp(\sig_\bC)Z$ is a $\sig$-nilpotent orbit (resp.\ $\exp(i\sig_\bR)Z$ is a $\sig$-nilpotent $i$-orbit).
  
\smallskip

(2) $V$ contains any element $h$ of $A^*$ such that $h:A\to \bQ$ sends $A\cap\sig$ to $\bQ_{\geq 0}$.  

\smallskip

Note that (2) implies $A\sub \sig_\bR$.

The canonical map $D_{\val}^{\sharp}\to D_{\val}$ (4.2.1) induces a canonical map $D_{\Sig,\val}^{\sharp}\to D_{\Sig,\val}$ which is surjective. 

\medskip

Our next subject is to define canonical maps
$D_{\Sig,\val} \to D_\Sig$ and $D_{\Sig,\val}^{\sharp}\to D_{\Sig}^\sharp$ for a weak fan $\Sig$. 

\proclaim{Lemma 4.2.3}
Let $\Sig$ be a weak fan and let $(A, V, Z)\in D_{\Sig,\val}$. Let $S$ be the set of all $\sig\in \Sig$ satisfying the conditions $(1)$ and $(2)$ in $4.2.2$. 
Then $S$ has a smallest element.
\endproclaim

{\it Proof.}
Let $T$ be the set of all rational nilpotent cones $\sig$ which satisfy (1) and (2) in 4.2.2. 
Hence $S=T\cap \Sig$. 
Let $T_1$ be the set of all finitely generated subcones of $A_\bR$ satisfying the condition 4.2.1 (1). 

Then we have
\medskip

{\bf Claim 1.} 
If $\sig_1\in T_1$ and $\sig_2\in T$, then $\sig_1\cap \sig_2\in T_1$.

\medskip

We prove Claim 1.
By $(\sig_1 \cap \sig_2)^\v = \sig_1^\v + \sig_2^\v$ (\cite{O88} A.1 (2)), 
4.2.2 (2) is satisfied for $\sig:=\sig_1 \cap \sig_2$.  
Here $(\;\;)^\v$ denotes the dual of a cone.
Then, $\dim \sig=\dim \sig_1$. 
On the other hand, $Z$ generates a $\sig_1$-nilpotent orbit by 
4.2.1 (1) for $\sig_1$. 
  Thus $Z$ also generates a $\sig$-nilpotent orbit. 
  Hence $\sig \in T_1$. 

\medskip

We now prove 4.2.3. 
Let $\sig_1$ and $\sig_2$ be minimal elements of $S$. We prove $\sig_1=\sig_2$. 
Since $(A,V,Z) \in D_{\val}$, the set $T_1$ is not empty. Let $\tau$ be an element of $T_1$, and let $\tau'=\tau\cap \sig_1\cap \sig_2$. Then $\tau'\in T_1$ by Claim 1. Since $\Sig$ is a weak fan, by 2.2.4 (2)$''$, the set
$\{\alpha\in S\;|\;\tau'\sub \alpha\}$ has a smallest element $\sig$. Since $\sig_1$ and $\sig_2$ belong to this set, $\sig\sub \sig_1$ and $\sig\sub \sig_2$. Since $\sig_1$ and $\sig_2$ are minimal in $S$, we have $\sig_1=\sig=\sig_2$. 
\qed

\medskip

{\bf 4.2.4.} 
Let $\Sig$ be a weak fan. We define the maps $D_{\Sig,\val}\to D_{\Sig}$ and 
$D_{\Sig,\val}^{\sharp}\to D_{\Sig}^{\sharp}$ as
$$
(A, V, Z) \mapsto (\sig, \exp(\sig_\bC)Z),  \quad 
(A,V,Z) \mapsto  (\sig, \exp(i\sig_\bR)Z),
$$
respectively, where $\sig$ is the smallest element of the set $S$ in 4.2.3. 

\medskip

{\bf 4.2.5.} 
The sets $E_{\sig,\val}$ and $E_{\sig,\val}^{\sharp}$.

Let $E_\sig$ and $E_\sig^\sharp$ be as in 2.3.3, and let  $\toric_{\sig,\val}$ over $\toric_\sig$ and  $\abtoric_{\sig,\val}$ be as in \cite{KU09} 5.3.6.
Define
$$
E_{\sig,\val}:=\toric_{\sig,\val} \times_{\toric_{\sig}}E_{\sig}, \quad
E_{\sig,\val}^{\sharp}=\abtoric_{\sig,\val} \times_{\abtoric_{\sig}} 
E_{\sig}^{\sharp}.
$$

Let $E_{\sig,\val}\to E_{\sig}$ and $E_{\sig,\val}^{\sharp}\to E_{\sig}^{\sharp}$ be the canonical projections. 

\medskip

{\bf 4.2.6.} 
Analogously as in \cite{KU09} 5.3.7, we have the projections $E_{\sig,\val} \to \G(\sig)^\gp\bs D_{\sig,\val}$ and $E^\sharp_{\sig,\val} \to D^\sharp_{\sig,\val}$. 

\medskip

{\bf 4.2.7.}
By 2.4 and \cite{KU09} 3.6.20, we define on $E_{\sig,\val}$ and $\G \bs D_{\Sig,\val}$ structures of log local ringed spaces over $\bC$, and  on $E_{\sig,\val}^{\sharp}$ and $D_{\Sig,\val}^{\sharp}$ topologies, analogously as in \cite{KU09} 5.3.6--5.3.8.
\medskip

\head
\S4.3. The map $D_{\val}^{\sharp} \to D_{\SL(2)}$.
\endhead
\medskip

We prove two theorems 4.3.1 and 4.3.2, which follow from propositions 4.3.3 and 4.3.4.

\proclaim{Theorem 4.3.1}
Let $p=(A,V,Z)\in D^\sharp_\val$.
\medskip

{\rm(i)} There exists a family $(N_j)_{1\le j\le n}$ of elements of $A_\bR = \bR \otimes_{\bQ} A$ satisfying the following two
conditions$:$
\smallskip

\noindent
{\rm(1)}
If $F\in Z$, $(N_1,\dots,N_n,F)$ generates a nilpotent orbit.

\smallskip

\noindent
{\rm(2)} Via $(N_j)_{1\le j\le n}:A^*=\Hom_\bQ(A,\bQ)\to\bR^n$, $V$ coincides with the set of all elements of $A^*$ whose images in $\bR^n$ are $\ge0$ with respect to the lexicographic order.
\medskip

{\rm(ii)} Take $(N_j)_{1\le j\le n}$ as in {\rm(i)}, let $F\in Z$ and define $\psi(p)\in D_{\SL(2)}$ to be the class of the $\SL(2)$-orbit associated to
$(N_1,\dots,N_n,F)$.
Then $\psi(p)$ is independent of the choices of $(N_j)_{1\le j\le n}$ and $F$.
\endproclaim

By Theorem 4.3.1, we obtain a map
$$
\psi:D^\sharp_\val\to D_{\SL(2)}.\tag3
$$
We call this map the {\it CKS map}, like in the pure case in \cite{KU09}. (CKS stands for Cattani-Kaplan-Schmid from whose work \cite{CKS86} the map $\psi$ in the pure case was defined in \cite{KU09}).

\proclaim{Theorem 4.3.2}
Let $\Sig$ be a weak fan in $\fg_\bQ$.
Then, $\psi:D_{\Sig,\val}^\sharp\to D_{\SL(2)}^I$ is continuous. 
\endproclaim

Thus $\psi:D_{\Sig,\val}^\sharp\to D_{\SL(2)}^I$ is the unique continuous extension of the identity map of $D$.

The $C^\infty$-property of the CKS map will be discussed in a later part of this series.

As in \cite{KU09}, 
these theorems follow from Proposition 4.3.3 below, which is an analog of 
\cite{KU09} 6.4.1.

\medskip

\proclaim{Proposition 4.3.3}
Let $(N_s)_{s\in S}$ be a finite family of mutually commuting
nilpotent elements of $\fg_\bR$, let $F\in\Dc$, and assume
that $((N_s)_{s\in S},F)$ generates a nilpotent orbit.
Let $a_s\in\bR_{>0}$ for $s\in S$.
Assume that $S$ is the disjoint union of nonempty subsets
$S_j$ $(1\le j\le n)$.
Denote $S_{\le j}:=\bigsqcup_{k \le j}S_k$ and $S_{\ge j}:=\bigsqcup_{k \ge j}S_k$. 
For $1\le j\le n$, let $\Dc_j$ be the subset of $\Dc$ consisting
of all $F'\in\Dc$ such that $((N_s)_{s\in S_{\le j}},F')$ generates
a nilpotent orbit. 
Let $L$ be a directed ordered set, let $F_\lam\in\Dc$
$(\lam\in L)$, $y_{\lam,s}\in\bR_{>0}$ $(\lam\in L,s\in S)$,
and assume that the following five conditions are satisfied.
\medskip

\noindent
{\rm(1)} $F_\lam$ converges to $F$.
\smallskip

\noindent
{\rm(2)} $y_{\lam,s}\to\infty$ for any $s\in S$.
\smallskip

\noindent
{\rm(3)} If $1\le j<n$, $s\in S_{\le j}$ and
$t\in S_{\ge j+1}$, then
$\tfrac{y_{\lam,s}}{y_{\lam,t}}\to\infty$.
\smallskip

\noindent
{\rm(4)} If $1\le j\le n$ and $s,t\in S_j$, then
$\tfrac{y_{\lam,s}}{y_{\lam,t}}\to
\tfrac{a_s}{a_t}$.
\smallskip

\noindent
{\rm(5)} For $1\le j\le n$ and $e\ge0$,
there exist
$F^*_\lam\in\Dc$ $(\lam\in L)$ and
$y^*_{\lam,t}\in\bR_{>0}$
$(\lam\in L,t\in S_{\ge j+1})$ such that
$$
\align
&\exp\big(\ts_{t\in S_{\ge j+1}}
iy_{\lam,t}^*N_t\big)F_\lam^*\in\Dc_j\quad
(\lam:\text{sufficiently large}),\\
&\,y_{\lam,s}^e d(F_\lam,F^*_\lam)\to0\quad
(\forall s\in S_j),\\
&\,y_{\lam,s}^e|y_{\lam,t}-y_{\lam,t}^*|\to0
\quad(\forall s\in S_j,\forall t\in S_{\ge j+1}).
\endalign
$$
Here $d$ is a metric on a neighborhood of $F$ in
$\Dc$ which is compatible with the analytic
structure of $\Dc$.
\medskip

For each $1\le j\le n$, take $c_j\in S_j$ and denote $N_j:=\ts_{s\in S_j}
\tfrac{a_s}{a_{c_j}}N_s$.

Then, $(N_1, \cdots, N_n, F)$ generates a nilpotent orbit, and 
we have the convergence
$$\exp(\ts_{s\in S}iy_{\lam,s}N_s)F_\lam\to \operatorname{class}(\psi(N_1, \cdots, N_n,
F)) \quad \text{in}\;\;D_{\SL(2)}^I.
$$
\endproclaim

By \cite{KNU.p} 3.2.12 (i) which characterizes the topology of $D_{\SL(2)}^I$, for the proof of Proposition 4.3.3, it is sufficient to prove the following

\medskip

\proclaim{Proposition 4.3.4}
Let the assumption be as in $4.3.3$. Let  
 $\br_1 \in D$ be the point and $\tau$ be the homomorphism which are associated to $(N_1,\dots,N_n,F)$ as in $4.1.3$.

Then we have the following convergences {\rm(1)} and {\rm(2)}.
\smallskip

{\rm(1)}
$\tau\big(\sqrt{\frac{y_{\lam,c_1}}{y_{\lam,c_2}}},\dots,
      \sqrt{\frac{y_{\lam,c_n}}{y_{\lam,c_{n+1}}}}\big)
\exp(\ts_{s\in S}iy_{\lam,s}N_s)F_\lam \to \br_1$ in $D$, 
where $y_{\lam,c_{n+1}}=1$.

\smallskip

{\rm(2)} $\spl_W(\exp(\ts_{s\in S}iy_{\lam,s}N_s)F_\lam) \to \spl_W(\br_1)$ in $\spl(W)$.

\smallskip

Here, for $F'\in D$, $\spl_W(F')\in \spl(W)$ denotes the canonical splitting of $W$ associated to $F'$ as in $4.1.1$. 
\endproclaim

(1) of Proposition 4.3.4 is the case $j=k$ ($k$ is as in 4.1.3) 
of Lemma 4.3.6 below by $\br_1=\exp(iN_k)\hat F_{(k)}$ (\cite{KNU.p} 2.4.8, 
Claim). 

\medskip

{\bf 4.3.5.}
Let the notation be as in 4.3.3.
Before stating 
the lemma, we recall the definition of $\hat F_{(j)}$ and $\tau_j$ associated to $(N_1,\dots, N_n, F)$ (cf.\ \cite{KNU08} 10.1.1,  \cite{KNU.p} 2.4.6).

Let $W_0 = W$ and, for $1\leq j\leq n$, let $W^{(j)}=M(N_1+\dots+N_j, W)$.

For $0 \leq j \leq n$, we define an $\bR$-split mixed Hodge structure $(W^{(j)}, \hat F_{(j)})$ and the associated splitting $s^{(j)}$ of $W^{(j)}$ inductively starting from $j=n$ and ending at $j=0$. 
First, $(W^{(n)}, F)$ is a mixed Hodge structure as is proved 
by Deligne (see \cite{K86} 5.2.1). 
Let $(W^{(n)}, \hat F_{(n)})$ be the $\bR$-split mixed Hodge
structure associated to the mixed Hodge structure $(W^{(n)}, F)$.
Then $(W^{(n-1)}, \exp(iN_n)\hat F_{(n)})$ is a mixed Hodge structure. 
Let $(W^{(n-1)}, \hat F_{(n-1)})$ be the $\bR$-split mixed Hodge structure associated to $(W^{(n-1)}, \exp(iN_n)\hat F_{(n)})$. 
Then $(W^{(n-2)}, \exp(iN_{n-1})\hat F_{(n-1)})$ is a mixed Hodge structure. 
This process continues. 
In this way we define $\hat F_{(j)}$ downward-inductively as the $\bR$-split mixed Hodge structure associated to the mixed Hodge structure $(W^{(j)}, \exp(iN_{j+1})\hat F_{(j+1)})$, and define $s^{(j)}: \gr^{W^{(j)}} \simeq H_{0,\bR}$ to be the splitting of $W^{(j)}$ associated to $\hat F_{(j)}$. 
The splitting $s$ in 4.1.3 is nothing but $s^{(0)}$ (\cite{KNU08} 10.1.2). 
The $j$-th component
$$
\tau_j : \bG_{m,\bR} \to \Aut_{\bR}(H_{0,\bR},W)\quad (0 \leq j \leq n)
$$ 
of the homomorphism $\tau$ associated to $(N_1,\dots,N_n, F)$ (4.1.3) is characterized as follows. For $a\in \bR^\times$ and $w\in \bZ$, $\tau_j(a)$ acts on $s^{(j)}(\gr^{W^{(j)}}_w)$ as the multiplication by $a^w$. 

\medskip

\proclaim{Lemma 4.3.6}
Let the notation be as in $4.3.3$.
Let $1\le j\le n$.
Then
$$
\Big(\tp_{j\le k\le n}\tau_k\Big(
\tsqrt{\tfrac{y_{\lam,c_k}}{y_{\lam,c_{k+1}}}}
\Big)\Big)
\exp(\ts_{s\in S_{\ge j}}iy_{\lam,s}N_s)
F_\lam \to\exp(iN_j)\hat F_{(j)}
\quad\text{in}\;\Dc,
$$
where $y_{\lam,c_{n+1}}:=1$.
$($Recall $N_j:=\ts_{s\in S_j}\tfrac{a_s}{a_{c_j}}N_s.)$
\endproclaim

{\it Proof.} This is a mixed Hodge version of \cite{KU09} Lemma 6.4.2, and proved by the same method in the proof of that lemma. At the place in the proof of that lemma where we used \cite{KU09} Lemma 6.1.10, we use its mixed Hodge version \cite{KNU08} 10.3. 
\qed

\medskip
{\bf 4.3.7.} 
We prove Proposition 4.3.4 (2).

We prove the following assertion $(A_j)$ by downward induction on $j$.
(Note that $(A_0)$ is what we want to prove.) 

\smallskip

$(A_j)$ : Proposition 4.3.4 (2) is true in the case  $\exp\big(\ts_{t\in S_{\ge j+1}} iy_{\lam,t}N_t\big)F_\lam\in\Dc_j$ for all $\lam$.   
Here $\check D_0=D$. 

\medskip

{\it Proof.} 
Let $0\leq j\leq n$.

Let $p_\lam= \exp(\ts_{s\in S}iy_{\lam,s}N_s)F_\lam$, 
$t_{\lam,j}= \tp_{j< k\le n}\tau_k\Big(
\tsqrt{\tfrac{y_{\lam,c_k}}{y_{\lam,c_{k+1}}}}\Big)$, and $p_{\lam,j}= t_{\lam,j}p_\lam$. 
Then, by \cite{KNU08} 10.3, 
$$
p_{\lam,j}= t_{\lam,j}p_\lam=\exp\big(\ts_{s\in S_{\leq j}}i
\tfrac{y_{\lam,s}}{y_{\lam,c_{j+1}}}N_s\big)U_{\lam,j},
$$
$$
\text{where}\quad 
U_{\lam,j}:=t_{\lam,j}\exp(\ts_{t\in S_{\ge j+1}}iy_{\lam,t}N_t)F_\lam.
$$
  Assume $\exp\big(\ts_{t\in S_{\ge j+1}} iy_{\lam,t}N_t\big)F_\lam\in\Dc_j$.
  Then 
$(N'_1,\ldots,N'_j, U_{\lam,j})$ generates a nilpotent orbit, 
where $N'_k=\sum_{s \in S_k}\frac{y_{\lam,s}}{y_{\lam,c_k}}N_s$ 
($1 \le k \le j$).
  Let $s_{\lam}$ be the associated limit splitting.
  By \cite{KNU08} 0.5 (2) and 10.8 (1), there is a convergent power series 
$u_{\lam}$ whose constant term is $1$ and 
whose coefficients depend on $U_{\lam,j}$ and $y_{\lam,s}/y_{\lam,c_k}$ ($1\leq k\leq j$, $s\in S_k$) real analytically 
such that 
$\spl_W(p_{\lam,j}) =u_{\lam}\Bigl(\frac{y_{\lam,c_2}}{y_{\lam,c_1}}, \dots, 
\frac{y_{\lam,c_{j+1}}}{y_{\lam,c_j}}\Bigr)s_{\lam}$.  
  Since $s_{\lam}$ also depends real analytically on $U_{\lam,j}$, 
we have 

\smallskip

\noindent
(1) $\spl_W(p_{\lam,j})$ converges to $\spl_W(\br_1)$.

\medskip

  This already showed ($A_n$). 

  Next, assume $j < n$ 
 and that $(A_{j+1})$ is true. We prove $(A_j)$ is true. 

  Choose a sufficiently big $e>0$ depending on $\tau_{j+1},\ldots,\tau_n$. 

  Take $F^*_{\lam}$ and $y^*_{\lam,t}$ ($t \in S_{\ge j+2}$) by 
the assumption 4.3.3 (5). 
  Let $y^*_{\lam,t}=y_{\lam,t}$ ($t \in S_{\le j+1}$). 
  Define $p_{\lam}^*$, $t_{\lam,j}^*$ and $p_{\lam,j}^*$ similarly 
as $p_{\lam}$, $t_{\lam,j}$ and $p_{\lam,j}$, respectively. 

Then we have 

\smallskip

\noindent
(2) $\spl_W(p^*_{\lam,j})$ converges to $\spl_W(\br_1)$, and
$y_{\lam, c_{j+1}}^ed(\spl_W(p_{\lam,j}), \spl_W(p^*_{\lam,j})) \to 0$.

\smallskip

By inductive hypothesis on $j$, we have

\smallskip

\noindent
(3) $\spl_W(p^*_{\lam})= t_{\lam,j}^{*-1} \spl_W(p^*_{\lam, j})
\gr^W(t^*_{\lam,j})$ converges to $\spl_W(\br_1)$.

\smallskip

By (1)--(3), we have $\spl_W(p_\lam)=t_{\lam,j}^{-1} \spl_W(p_{\lam,j}) \gr^W(t_{\lam,j})$ also converges to $\spl_W(\br_1)$.\qed

\bigskip

\head
\S5. Proofs of the main results
\endhead

In this section, we prove the results in \S2 and \S3. 

The results in \S3 are reduced to the ones in \S2.  
We explain this reduction in the last subsection 5.6 in this section. 

The proofs 5.1--5.5 
for the results in \S2 are the mixed Hodge theoretic and 
weak fan versions of 
\cite{KU09} \S7, and arguments are parallel to those in \cite{KU09}. 
We replace the results and lemmas used there with the 
corresponding ones here (e.g., the Hausdorffness of $D_{\SL(2)}$ in 
the pure case (\cite{KU02} 3.14 (ii)) there is replaced by the 
corresponding result in the mixed Hodge situation (\cite{KNU.p})), and 
make obvious modifications (e.g., fans are replaced by weak fans, 
the associated filtration $W(-)$ is replaced by $M(-,W)$, 
$\tilde \rho$ is replaced by $\tau$,...). 
For the readers' conveniences, we include several important rewritten 
arguments from those of \cite{KU09} \S7, even in the case they are more or less 
repetitions of loc.\ cit.

\head
\S5.1.  Proof of Theorem 2.5.1
\endhead

We prove Theorem 2.5.1.

\medskip

{\bf 5.1.1.} 
The fact that $\tilde  E_\sig$ is a log manifold is proved by the similar arguments in the proof of the pure case  \cite{KU09} 3.5.10. 

\medskip

{\bf 5.1.2.}
We prove that $E_\sig$ is open in $\tilde E_\sig$ for the strong topology of $\tilde E_\sig$ in $\Ec_\sig$. 

For a face $\tau$ of $\sig$, let $U(\tau)$ be the open set of $\toric_\sig$ defined to be $\Spec(\bC[P])_\an$ where $P\subset \Hom(\G(\sig)^{\gp},\bZ)$ is the inverse image of $\G(\tau)^{\spcheck}\sub \Hom(\G(\tau)^{\gp},\bZ)$ under the canonical map $\Hom(\G(\sig)^{\gp}, \bZ)\to \Hom(\G(\tau)^{\gp},\bZ)$. Let $V(\tau)$ be the open set of $\tilde E_\sig$ defined to be $\tilde E_\sig\cap (U(\tau) \times \Dc)$. Let $V$ be the union of $V(\tau)$ where $\tau$ ranges over all faces of $\sig$ which are admissible. Then $V$ is open in $\tilde E_\sig$, and $E_\sig\sub V$. 
Concerning the canonical projection $V \to \prod_w \tilde E_{\sig(\gr^W_w)}$, $E_\sig$ is the inverse image of $\prod_w 
E_{\sig(\gr^W_w)}$ in $V$. By \cite{KU09} 4.1.1 Theorem A (i), each  $E_{\sig(\gr^W_w)}$ is open in $\tilde E_{\sig(\gr^W_w)}$ in the strong topology of $\tilde E_{\sig(\gr^W_w)}$ in $\Ec_{\sig(\gr^W_w)}$. 
Hence  $E_\sig$ is open in $\tilde E_\sig$ in the strong topology of $\tilde E_\sig$ in $\Ec_\sig$.

 This shows that
  $E_{\sig}$ is  a log manifold.

  \medskip

\head
\S5.2. Study of $D_\sig^{\sharp}$
\endhead

Let $\Sig$ be a weak fan, $\G$ a subgroup of $G_\bZ$, and assume that they are strongly compatible. Let $\sig\in \Sig$. 

We prove some results on $D_\sig^{\sharp}$.

\medskip

{\bf5.2.1.}
The actions of $\sig_\bC$ 
$$
\sig_\bC \times E_\sig \to E_\sig,\quad
(a, (q,F)) \mapsto (\be(a)q,\exp(-a)F),
$$ 
where $a \in\sig_\bC$, $q\in\toric_\sig$,
$F\in\Dc$ and $(q,F)\in E_\sig$, 
and 
$$
\sig_\bC \times E_{\sig,\val}\to E_{\sig,\val},\quad
(a, (q,F)) \mapsto (\be(a)q,\exp(-a)F),
$$ 
where 
$a \in\sig_\bC$, $q\in\toric_{\sig,\val}$,
$F\in\Dc$ and $(q,F)\in E_{\sig,\val}$, 
are proved to be continuous by the same argument in the 
pure case \cite{KU09} 7.2.1.
  They induce the actions of $i\sig_\bR$ on
$E_\sig^\sharp$ and $E_{\sig,\val}^\sharp$.

The quotient spaces for these actions are identified as
$$
\align
&\sig_\bC\bs E_\sig=\G(\sig)^\gp\bs D_\sig,\quad
\sig_\bC\bs E_{\sig,\val}
=\G(\sig)^\gp\bs D_{\sig,\val},\\
&i\sig_\bR\bs E_\sig^\sharp=D_\sig^\sharp,\quad
i\sig_\bR\bs E_{\sig,\val}^\sharp
=D_{\sig,\val}^\sharp.
\endalign
$$

\proclaim{Proposition 5.2.2}
{\rm(i)} The action of $\sig_\bC$ on $E_\sig$
$($resp.\ $E_{\sig,\val})$ is free.

\medskip

{\rm(ii)} The action of $i\sig_\bR$ on $E^\sharp_\sig$
$($resp.\ $E^\sharp_{\sig,\val})$ is free.
\endproclaim

Here and in the following, free action always means
set-theoretically free action. 
That is, an action of a group $H$ on a set $S$ is free if and only if $hs\neq s$  for any $h\in H-\{1\}$ and $s\in S$.

Before proving 5.2.2, we state a lemma, whose variant in the pure case is \cite{KU09} 6.1.7 (2).

\proclaim{Lemma 5.2.2.1} 
Let $(H_\bR, W,F)$ be an $\bR$-mixed Hodge structure, let
$$
\delta(F) \in L_\bR^{-1,-1}(F(\gr^W))\sub \End_{\bR}(H(\gr^W)_{\bR})
$$ be as in \cite{KU09} \S$4$, \cite{KNU.p} $1.2.2$, and  
let $s' \in \spl(W)$ be the splitting of $W$ characterized by 
$F=s'(\exp(i\delta(F))F(\gr^W))$ $($loc.\ cit.$)$. 

  Let $N \in \End_{\bR}(H_{\bR})$ and assume that 
$(s')^{-1}Ns' \in L_{\bR}^{-1,-1}(F(\gr^W))$ 
 and that 
$(s')^{-1}Ns'$ commutes with $\delta(F)$. 
  Let $z \in \bC$. 
  Then, we have $\delta(\exp(zN)F) = \delta(F) + \Im(z) (s')^{-1}Ns'$. 
\endproclaim

\demo{Proof}
  Write $z=x+iy$ with $x$ and $y$ real. 
  Let $\delta=\delta(F)$. 
  Apply $\exp(zN)$ to $F=s'(\exp(i\delta)F(\gr^W))$.
  Then, we have $\exp(zN)F = (\exp(xN)s')\exp(iy(s')^{-1}Ns')\exp(i\delta)
F(\gr^W)$. 
  By the commutativity assumption, $\exp(iy(s')^{-1}Ns')\exp(i\delta)
= \exp(i(y(s')^{-1}Ns'+\delta))$, and the conclusion follows. 
\qed
\enddemo

\demo{Proof of 5.2.2}
We show here that the action of $\sig_\bC$ on
$E_\sig$ is free.
The rest is proved in the similar way.

Assume  $a\cdot(q,F)=(q,F)$ $(a\in\sig_\bC$,
$(q,F)\in E_\sig)$, that is, $\be(a)\cdot q=q$ and $\exp(-a)F=F$. 

By the pure case considered in \cite{KU09} 7.2.9, 
the image of $a$ in $\fg_\bC(\gr^W)$ is zero. 
  The rest of the proof of 5.2.2 is parallel to the argument in \cite{KU09} 7.2.9 as follows. 
  
By $\be(a)\cdot q=q$, we have
$a=b+c$ ($b\in\sig(q)_\bC\cap \fg_{\bC,u}$, $c\in\log(\G(\sig)^\gp\cap G_{\bZ,u})$), where 
$\sig(q)$ is the unique face of $\sig$ such that $q \in 
\bold e(\sig_{\bC})\cdot 0_{\sig(q)}$. 
On the other hand, taking
$\delta(\;\;)=\delta(M(\sig(q), W),\;\;)$ of $\exp(-a)F=F$, we
have $\delta(F)=\delta(\exp(-a)F)=-\Im(b)+\delta(F)$ by 5.2.2.1.
It follows $\Im(b)=0$ and hence $\exp(\Re(b)+c)F=F$.
Take an element $y$ of the interior of
$\sig(q)$ such that $F':=\exp(iy)F\in D$.
Then $\exp(\Re(b)+c)F'=F'$.
Since $\Re(b)+c\in \fg_{\bR,u}$, 
we have $\exp(\Re(b)+c)\spl_W(F')=\spl_W(F')$. 
This proves  $\Re(b)+c=0$ and hence $a=0$. 
\qed
\enddemo

\medskip

\proclaim{Proposition 5.2.3}
Let $\sig$ and $\sig'$ be sharp rational nilpotent cones.
Let $\a\in E^\sharp_{\sig,\val}$ and
$\a'\in E^\sharp_{\sig',\val}$.
Assume that $(\be(iy_\lam),F_\lam)\in
E^\sharp_{\sig,\val,\triv}$ $(y_\lam\in\sig_\bR$,
$F_\lam\in\Dc)$ $($resp.\ $(\be(iy'_\lam),F'_\lam)\in
E^\sharp_{\sig',\val,\triv}$ $(y'_\lam\in\sig'_\bR$,
$F'_\lam\in\Dc))$ converges to $\a$ $($resp.\ $\a')$ in
the strong topology (here $(-)_{\triv}$ means the part where the log structure of $E_{\sig,\val}$ or of $E_{\sig',\val}$ becomes trivial, respectively), and that
$$
\exp(iy_\lam)F_\lam=\exp(iy'_\lam)F'_\lam\;\;
\text{in}\;\;D.
$$
Then we have$:$
\medskip

{\rm(i)} The images of $\a$ and $\a'$ in
$D^\sharp_\val$ coincide.
\medskip

{\rm(ii)} $y_\lam-y'_\lam$ converges in $\fg_\bR$.
\endproclaim

{\it Proof.}
  First of all, we remark that 
the proof below is just a revision of that of \cite{KU09} 
7.2.10.  
  Notice that in the course of the proof of Claim A in the mixed Hodge case below, the map
$\b$ is 
undefined on $D_{\spl}$ when $W^{(1)}=W$.
  But, in this case, $\br_1 \not\in D_{\spl}$, so $\b$ can be still 
applied. 

  Since the composite map $E^\sharp_{\sig,\val}$
(resp.\ $E^\sharp_{\sig',\val})\to D^\sharp_\val
@>\psi>>D_{\SL(2)}^I$ is continuous by 4.2.7 
and Theorem 4.3.2, the image of $\a$ (resp.\ $\a'$)
under this composite map is the limit of
$\exp(iy_\lam)F_\lam$
(resp.\ $\exp(iy'_\lam)F'_\lam$) in $D_{\SL(2)}^I$.
Since $\exp(iy_\lam)F_\lam=\exp(iy'_\lam)F'_\lam$
and since $D_{\SL(2)}^I$ is Hausdorff (\cite{KNU.p} 3.5.17 (i)),
these images coincide, say $p\in D_{\SL(2)}^I$.
Let $m$ be the rank of $p$ and
let $\Psi=(W^{(k)})_{1\le k\le m}$ be the associated family
of weight filtrations (4.1.1, 4.1.2).
Let $(A,V,Z)$ (resp.\ $(A',V',Z')$) be the image of $\a$
(resp.\ $\a'$) in $D^\sharp_\val$.
Take an excellent basis $(N_s)_{s\in S}$ (resp.\
$(N'_s)_{s\in S'}$) for $(A,V,Z)$ (resp.\ $(A',V',Z')$) in the mixed Hodge case 
(the definition of an excellent basis \cite{KU09} 6.3.8 
obviously extends to the mixed Hodge case, and its existence is seen similarly 
to \cite{KU09} 6.3.9), and let $(a_s)_{s\in S}$ and $(S_j)_{1\le j\le n}$ (resp.\ $(a'_s)_{s\in S'}$ and $(S'_j)_{1\le j\le n'}$) be as in (the mixed Hodge version of) \cite{KU09} 6.3.3.
For each $l$ with $1\le l\le m$, let $f(l)$ (resp.\ 
$f'(l)$) be the smallest integer such that $1\le f(l)\le n$
(resp.\ $1\le f'(l)\le n'$), and that $W^{(l)}=M(\ts_{s\in
S_{\le f(l)}}\bR_{\ge0}N_s, W)$ (resp.\ $=M(\ts_{s\in
S'_{\le f'(l)}}\bR_{\ge0}N'_s, W)$). 
%and that $W^{(l)}\ne
%W(\ts_{s\in S_{\le f(l)-1}}(\bR_{\ge0})N_s)$ (resp.\
%$\ne W(\ts_{s\in S'_{\le f'(l)-1}}(\bR_{\ge0})N'_s)$).
Let $g(l)$ (resp.\ $g'(l)$) be any element of $S_{f(l)}$
(resp.\ $S'_{f'(l)}$).

Take an $\bR$-subspace $B$ of $\sig_\bR$ (resp.\ $B'$
of $\sig'_\bR$) such that $\sig_\bR=A_\bR\op B$
(resp.\ $\sig'_\bR=A'_\bR\op B'$).
Write
$$
y_\lam=\ts_{s\in S}y_{\lam,s}N_s+b_\lam,\quad
y'_\lam=\ts_{s\in S'}y'_{\lam,s}N'_s+b'_\lam
$$
with $y_{\lam,s}$, $y'_{\lam,s}\in\bR$, $b_\lam\in B$,
$b'_\lam\in B'$.
We may assume $y_{\lam,s}$, $y'_{\lam,s}>0$.
Then, by the mixed Hodge version of \cite{KU09} 6.4.11, proved similarly, 
$y_{\lam,s}$ $(s\in S)$,
$y'_{\lam,s}$ $(s\in S')$,
$\tfrac{y_{\lam,s}}{y_{\lam,t}}$ $(s\in S_{\le j}$,
$t\in S_{\ge j+1})$ and
$\tfrac{y'_{\lam,s}}{y'_{\lam,t}}$ ($s\in S'_{\le j}$,
$t\in S'_{\ge j+1}$) tend to $\infty$,
$\tfrac{y_{\lam,s}}{y_{\lam,t}}$ $(s,t\in S_j)$ tends to
$\tfrac{a_s}{a_t}$, $\tfrac{y'_{\lam,s}}{y'_{\lam,t}}$
$(s,t\in S'_j)$ tends to $\tfrac{a'_s}{a'_t}$, and
$b_\lam$ and $b'_\lam$ converge.

\proclaim{Claim A}
$\tfrac{y'_{\lam,g'(l)}}{y_{\lam,g(l)}}$ converges to an
element of $\bR_{>0}$ for $1\le l\le m$.
\endproclaim

We prove this claim.
The assumption of Proposition 4.3.3 is satisfied if we
take $(N_s)_{s\in S}$, $a_s$, $S_j$, $y_{\lam,s}$ as
above the same ones in 4.3.3, the above
$\exp(b_\lam)F_\lam$ as $F_\lam$ in 4.3.3, the limit
of the above $\exp(b_\lam)F_\lam$ as $F$ in 4.3.3.
Let $\tau: (\bR^\times)^m\to \Aut_\bR(H_{0,\bR}, W)$ be the homomorphism associated to $p$ (4.1.1, 4.1.2). 
Then, by 4.3.4 (1), we have
$$
\tau\Big(
\sqrt{\tfrac{y_{\lam,g(1)}}{y_{\lam,g(2)}}},\dots,
\sqrt{\tfrac{y_{\lam,g(m)}}{y_{\lam,g(m+1)}}}\Big)
\exp(iy_\lam)F_\lam\to\br_1\quad
\text{in $D$}.
$$
See 4.1.3 for $\br_1$. 
Similarly,  for some $t=(t_1,\dots,t_m)\in\bR_{>0}^m$, we have
$$
\tau\bigg(
\sqrt{\tfrac{y'_{\lam,g'(1)}}{y'_{\lam,g'(2)}}},\dots,
\sqrt{\tfrac{y'_{\lam,g'(m)}}{y'_{\lam,g'(m+1)}}}\bigg)
\exp(iy'_\lam)F'_\lam\to\tau(t)\cdot\br_1\quad
\text{in $D$}
$$
(4.1.2).
Here $y_{\lam,g(m+1)}=y'_{\lam,g'(m+1)}=1$.
Recall that $\Psi=\{W^{(1)}, \dots, W^{(m)}\}$ is the set of weight filtrations associated to $p$. 
In the case $W\notin \Psi$ (resp.\ $W\in \Psi$), take a distance to $\Psi$-boundary 
$$
\b:D\to\bR_{>0}^m \quad 
(\text{resp.} \;\; D_{\text{nspl}}:=D\smallsetminus D_{\spl}\to \bR_{>0}^m)
$$
in \cite{KNU.p} 3.2.5.
So $\b(\tau(t)x)=t\b(x)$ $(x\in D$ (resp.\ $D_{\text{nspl}})$, $t\in\bR_{>0}^m)$.
Note that, in the case $W\in \Psi$, $\br_1\in D_{\text{nspl}}$. 
Applying $\b$ to the above convergences, taking their
ratios, and using $\exp(iy_\lam)F_\lam
=\exp(iy'_\lam)F'_\lam$, we have
$$
\tfrac{y'_{\lam,g'(l)}}{y_{\lam,g(l)}}
\to t_l^2\quad (1\le l\le m).
$$
Claim A is proved.

Next, we prove the following Claims B$_l$ and C$_l$
$(1\le l\le m+1)$ by induction on $l$.

\proclaim{Claim B$_l$}
$y_{\lam,g(l)}^{-1}\big(\ts_{s\in S_{\le f(l)-1}}
y_{\lam,s}N_s -\ts_{s\in S'_{\le f'(l)-1}}
y'_{\lam,s}N'_s\big)$ converges.
\endproclaim

\proclaim{Claim C$_l$}
$\ts_{s\in S_{\le f(l)-1}}\bQ N_s
=\ts_{s\in S'_{\le f'(l)-1}}\bQ N'_s$.
\endproclaim

\noindent
Here $y_{\lam,g(m+1)}:=1$, $S_{\le f(m+1)-1}:=S$,
$S'_{\le f'(m+1)-1}:=S'$.

Note that Proposition 5.2.3 follows from Claims
B$_{m+1}$ and C$_{m+1}$.
In fact, $A=A'$ follows from Claim C$_{m+1}$, $V=V'$
follows from Claim B$_{m+1}$, and $Z=Z'$ follows from the
facts that the limit of $\exp(b_\lam)F_\lam$ (resp.\
$\exp(b'_\lam)F_\lam'$) is an element of $Z$ (resp.\
$Z'$) and that these limits coincide by Claim B$_{m+1}$,
Claim C$_{m+1}$ and the assumption $\exp(iy_\lam)F_\lam
=\exp(iy'_\lam)F'_\lam$.

We prove these claims.
First, Claims B$_1$ and C$_1$ are trivial by definition. 
Assume $l>1$.
By the hypothesis Claim C$_{l-1}$ of induction, $N_s$
$(s\in S_{\le f(l-1)-1})$ and elements of $\sig'$ are commutative.
Hence, by the formula $\exp(x_1+x_2)=
\exp(x_1)\exp(x_2)$ if $x_1x_2=x_2x_1$ and by the
assumption $\exp(iy_\lam)F_\lam=\exp(iy'_\lam)
F'_\lam$, we have
$$
\exp\big(iy_\lam-\ts_{s\in S_{\le f(l-1)-1}}
iy_{\lam,s}N_s\big)F_\lam
=\exp\big(iy'_\lam-\ts_{s\in S_{\le f(l-1)-1}}
iy_{\lam,s}N_s\big)F'_\lam.
$$
Applying $\tp_{k=l}^m\tau_k\Big(
\tsqrt{\tfrac{y_{\lam,g(k)}}{y_{\lam,g(k+1)}}}\Big)$
to this and using \cite{KNU08} 10.3, we obtain
$$
\align
&\exp\big(\ts_{f(l-1)\le j<f(l)}\ts_{s\in S_j}
i\tfrac{y_{\lam,s}}{y_{\lam,g(l)}}N_s\big)\tag1\\
&\hskip80pt\cdot
\tp_{k\ge l}\tau_k
\Big(\tsqrt{\tfrac{y_{\lam,g(k)}}
{y_{\lam,g(k+1)}}}\Big)
\exp\big(\ts_{s\in S_{\ge
f(l)}}iy_{\lam,s}N_s+ib_\lam\big)F_\lam\\
=\;&\exp\big(\ts_{s\in S'_{\le f'(l)-1}}
i\tfrac{y'_{\lam,s}}{y_{\lam,g(l)}}N'_s -\ts_{s\in
S_{\le f(l-1)-1}}
i\tfrac{y_{\lam,s}}{y_{\lam,g(l)}}N_s\big)\\
&\hskip80pt\cdot
\tp_{k\ge l}\tau_k\Big(
\tsqrt{\tfrac{y_{\lam,g(k)}}{y_{\lam,g(k+1)}}}\Big)
\exp\big(\ts_{s\in S'_{\ge f'(l)}}
iy'_{\lam,s}N'_s+ib'_\lam\big)F'_\lam.
\endalign
$$
By Lemma 4.3.6,
$$
\align
&\tp_{k\ge l}\tau_k
\Big(\tsqrt{\tfrac{y_{\lam,g(k)}}{y_{\lam,g(k+1)}}}
\Big)\exp\big(\ts_{s\in S_{\ge f(l)}}
iy_{\lam,s}N_s+ib_\lam\big)F_\lam\;\;\text{and}\tag2\\
&\tp_{k\ge l}\tau_k\Big(
\tsqrt{\tfrac{y_{\lam,g(k)}}{y_{\lam,g(k+1)}}}\Big)
\exp\big(\ts_{s\in S'_{\ge f'(l)}}
iy'_{\lam,s}N'_s+ib'_\lam\big)F_\lam'\;\;\text{converge in $\Dc$.}
\endalign
$$
Let $d$ (resp.\ $d'$) be a metric on a neighborhood in $\Dc$ of
the limit of $F_\lam$ (resp.\ $F'_\lam$) which is
compatible with the analytic structure.
Let $e\ge1$ be an integer.
Then, since $(\be(iy_\lam),F_\lam)$ (resp.\
$(\be(iy'_\lam),F'_\lam)$) converges to $\a$ (resp.\
$\a'$) in the strong topology, there exist, by \cite{KU09} 3.1.6, $y_\lam^*
=\ts_{s\in S}y_{\lam,s}^*N_s+b_\lam^*\in\sig_\bR$
$(b_\lam^*\in B_\bR)$, $y_\lam^{\p*}
=\ts_{s\in S'}y_{\lam,s}^{\p*}N'_s+b_\lam^{\p*}
\in\sig'_\bR$ $(b_\lam^{\p*}\in B'_\bR)$ and
$F_\lam^*,F_\lam^{\p*}\in\Dc$ having the following
three properties:
\medskip

\noindent
(3) $y_{\lam,g(l-1)}^e(y_\lam-y_\lam^*)$,
$y_{\lam,g(l-1)}^e(y'_\lam-y_\lam^{\p*})$,
$y_{\lam,g(l-1)}^ed(F_\lam,F_\lam^*)$, and
$y_{\lam,g(l-1)}^ed(F'_\lam,F_\lam^{\p*})$
\newline
converge to $0$.
\medskip

\noindent
(4) $y_{\lam,s}=y_{\lam,s}^*$ $(s\in S_{\le f(l-1)})$
and $y'_{\lam,s}=y^{\p*}_{\lam,s}$ $(s\in S'_{\le
f'(l-1)})$.
\medskip

\noindent
(5) $\big((N_s)_{s\in S_{\le f(l)-1}},\,
\exp\big(\ts_{s\in S_{\ge f(l)}}
iy_{\lam,s}^*N_s+b_\lam^*\big)F_\lam^*\big)$ and
\newline
$\big((N'_s)_{s\in S'_{\le f'(l)-1}},\,
\exp\big(\ts_{s\in S'_{\ge f'(l)}}
iy_{\lam,s}^{\p*}N'_s+b_\lam^{\p*}\big)
F_\lam^{\p*}\big)$ generate nilpotent orbits.
\medskip

\noindent
Take $e$ sufficiently large.
Then, by the hypothesis Claim B$_{l-1}$ of induction
and by (1) and (2), 
the difference between 
$\delta(W^{(l-1)},\;\;)$ of
$$
\align
&\exp\Big(\ts_{f(l-1)\le j<f(l)}\ts_{s \in S_j}
i\tfrac{y_{\lam,s}}{y_{\lam,g(l)}}N_s\Big)\\
&\hskip96pt
\cdot\tp_{k\ge l}\tau_k\Big(
\tsqrt{\tfrac{y_{\lam,g(k)}}{y_{\lam,g(k+1)}}}\Big)
\exp\big(\ts_{s\in S_{\ge f(l)}}
iy_{\lam,s}^*N_s+ib_\lam^*\big)F_\lam^*
\endalign
$$
and $\delta(W^{(l-1)},\;\;)$ of
$$
\align
&\exp\Big(\ts_{s\in S'_{\le f'(l)-1}}
i\tfrac{y'_{\lam,s}}{y_{\lam,g(l)}}N'_s
-\ts_{s\in S_{\le f(l-1)-1}}
i\tfrac{y_{\lam,s}}{y_{\lam,g(l)}}N_s\Big)\\
&\hskip96pt
\cdot\tp_{k\ge l}\tau_k\Big(
\tsqrt{\tfrac{y_{\lam,g(k)}}{y_{\lam,g(k+1)}}}\Big)
\exp\big(\ts_{s\in S'_{\ge f'(l)}}
iy_{\lam,s}^{\p*}N'_s+ib_\lam^{\p*}\big)F_\lam^{\p*}
\endalign
$$
converges to $0$. By 5.2.2.1, 
the former is equal to
$$
\ts_{f(l-1)\le j<f(l)} \ts_{s \in S_j}
i\tfrac{y_{\lam,s}}{y_{\lam,g(l)}}N_s
+(\text{a term which converges}),
$$
and the latter is equal to
$$
\ts_{s\in S'_{\le f'(l)-1}}
i\tfrac{y'_{\lam,s}}{y_{\lam,g(l)}}N'_s
-\ts_{s\in S_{\le f(l-1)-1}}
i\tfrac{y_{\lam,s}}{y_{\lam,g(l)}}N_s
+(\text{a term which converges}).
$$
This proves Claim B$_l$.
Claim C$_l$ follows from Claim B$_l$ easily. 
\qed

\medskip

We recall the notion of proper action. 
In this series of papers, a continuous map of topological spaces is said to be 
proper if it is proper in the sense of \cite{B66} and is separated. 

\medskip

\definition{Definition 5.2.4
{\rm(\cite{B66} Ch.3 \S4 no.1 Definition 1)}}
Let $H$ be a Hausdorff topological group acting continuously on a
topological space $X$.
$H$ is said to act {\rm properly} on $X$ if the map
$$
H\x X\to X\x X,\;\;(h,x)\mapsto(x,hx),
$$
is proper.
\enddefinition

As is explained above, the
meaning of properness of a continuous map is slightly different
from that in \cite{B66}. However, for a continuous action of
a Hausdorff topological group $H$ on a topological
space $X$, the map $H\x X\to X\x X,(h,x)\mapsto (x, hx)$
is always separated, and hence this map is proper
in our sense if and only if it is proper
in the sense of \cite{B66}.

The following 5.2.4.1--5.2.4.7 hold. Recall (5.2.2) that when we say an action is free, it means free set-theoretically.  

\medskip

{\bf 5.2.4.1.}
If a Hausdorff topological group $H$ acts properly on a
topological space $X$, then the quotient space
$H\bs X$ is Hausdorff.

\medskip
{\bf 5.2.4.2.}
If a discrete group $H$ acts properly and freely on a
Hausdorff space $X$, then the projection $X\to H\bs X$
is a local homeomorphism.
\medskip

{\bf 5.2.4.3.}
Let $H$ be a Hausdorff topological group acting continuously on
topological spaces $X$ and $X'$.
Let $\psi:X\to X'$ be an equivariant continuous map.

\medskip

{\rm(i)} If $\psi$ is proper and surjective and if $H$
acts properly on $X$, then $H$ acts properly on $X'$.
\medskip

{\rm(ii)} If $H$ acts properly on $X'$ and if $X$ is
Hausdorff, then $H$ acts properly on $X$.
\medskip

{\bf 5.2.4.4.}
Assume that a Hausdorff topological group $H$ acts on a
topological space $X$ continuously and freely.
Let $X'$ be a dense subset of $X$.
Then, the following two conditions $(1)$ and $(2)$ are
equivalent.
\medskip

\noindent
{\rm(1)} The action of $H$ on $X$ is proper.
\medskip

\noindent
{\rm(2)} Let $x,y\in X$, $L$ be a directed set,
$(x_\lam)_{\lam\in L}$ be a family of elements of
$X'$ and $(h_\lam)_{\lam\in L}$ be a family of
elements of $H$, such that $(x_\lam)_\lam$ $($resp.\
$(h_\lam x_\lam)_\lam)$ converges to $x$ $($resp.\
$y)$.
Then
$(h_\lam)_\lam$ converges to an element $h$ of $H$
and $y=hx$.
\medskip

If $X$ is Hausdorff, these equivalent conditions are
also equivalent to the following condition {\rm(3)}.
\medskip

\noindent
{\rm(3)} Let $L$ be a directed set,
$(x_\lam)_{\lam\in L}$ be a family of elements of
$X'$ and $(h_\lam)_{\lam\in L}$ be a family of
elements of $H$, such that $(x_\lam)_\lam$ and
$(h_\lam x_\lam)_\lam$ converge in $X$.
Then $(h_\lam)_\lam$ converges in $H$.
\medskip

{\bf 5.2.4.5.} Let $H$ be a Hausdorff topological group acting continuously on Hausdorff topological spaces $X_1$ and $X_2$. Let $H_1$ be a closed normal subgroup of $H$, and assume that the action of $H$ on $X_2$ factors through $H_2:=H/H_1$. Assume that for $j=1, 2$, the action of $H_j$ on $X_j$ is proper and free. Assume further that there are a neighborhood $U$ of $1$ in $H_2$ and a continuous map $s: U\to H$ such that the composition $U@>s>> H\to H_2$ is the inclusion map. Then the diagonal action of $H$ on $X_1\times X_2$ is proper and free. 
\medskip

{\bf 5.2.4.6.} Let $H$ be a topological group acting continuously on a Hausdorff topological space $X$. Assume that $Y:= H\bs X$
is Hausdorff and assume that $X$ is an $H$-torsor over $Y$ in the category of topological spaces. Then the action of $H$ on $X$ is proper.

\medskip

The following 5.2.4.7 is not related to a group action, but we put it here for 5.2.4.1 and 5.2.4.3 (i) imply it in the case $H=\{1\}$. 

\medskip

{\bf 5.2.4.7.} Let $f: X\to Y$ be a continuous map between topological spaces and  assume that $f$ is proper and surjective. Assume 
$X$ is Hausdorff. Then $Y$ is Hausdorff.

\medskip

For 5.2.4.1, 5.2.4.2, 5.2.4.3, see \cite{B66} Ch.3 \S4 no.2 Proposition 3, ibid. Ch.3 \S4 no.4 Corollary, ibid. Ch.3 \S2 no.2 Proposition 5, respectively. 
5.2.4.4 is \cite{KU09} Lemma 7.2.7. See \cite{KU09} for the proof of it.

5.2.4.5 is proved as follows. The freeness is clear. We prove the properness. Let $L$ be a directed set, let $(x_\lam, y_\lam)_{\lam\in L}$ be a family of elements of $X_1\times X_2$, let $(h_\lam)_{\lam\in L}$ be a family of elements of $H$. Assume that $x_\lam$, $y_\lam$, $h_\lam x_\lam$, $h_\lam y_\lam$ converge. By 5.2.4.4, it is sufficient to prove that $h_\lam$ converges. Let $\bar h_\lam$ be the image of $h_\lam$ in $H_2$. Since the action of $H_2$ on $X_2$ is proper, $\bar h_\lam$ converges by 5.2.4.4. Let $h\bmod H_1$ ($h\in H$) be the limit. By replacing $h_\lam$ by $h^{-1}h_\lam$, we may assume $\bar h_\lam\to 1$. By replacing $h_\lam$ by $s(\bar h_\lam)^{-1}h_\lam$, we may assume $h_\lam\in H_1$. Since the action of $H_1$ on $X_1$ is proper, $h_\lam$ converges by 5.2.4.4.

5.2.4.6 is proved as follows. Let $L$ be a directed set, let $(x_\lam)_{\lam\in L}$ be a family of elements of $X$, and let $(h_\lam)_{\lam\in L}$ be a family of elements of $H$. Assume that $x_\lam$ converges to $x\in X$ and $h_\lam x_\lam$ converges to $y\in X$. By 5.2.4.4, it is sufficient to prove that $h_\lam$ converges. Since $Y=H\bs X$ is Hausdorff, the images of $x$ and $y$ in $Y$ coincide. Let $z\in Y$ be their image. By replacing $Y$ by a sufficiently small neighborhood of $z$, we may assume that $X=H\times Y$ with the evident action of $H$. Then the convergence of $(h_\lam)_\lam$ is clear.
\medskip

\proclaim{Proposition 5.2.5}
The action of $i\sig_\bR$ on $E^\sharp_\sig$
$($resp.\ $E^\sharp_{\sig,\val})$ is proper.
\endproclaim

\medskip
{\it Proof.} 
The proof is exactly the same as in the pure case \cite{KU09} 7.2.11, that is, 
it reduces to 5.2.3 (ii) by 5.2.4.3 (i) and 5.2.4.4. 
\qed
\medskip

By 5.2.5 and  5.2.4.1, we have the following result. 

\proclaim{Corollary 5.2.6}
The spaces $D_\sig^\sharp$ and $D_{\sig,\val}^\sharp$ are Hausdorff.
\endproclaim

This corollary will be generalized in \S5.3 
by replacing $\sig$ with $\Sig$ (5.3.3).

\medskip

The case (b) of the following Lemma will be used in \S5.4 later. 

\proclaim{Lemma 5.2.7}
Let $\cC$ be either one of the following two categories$:$
\roster
\item"{\rm(a)}" the category of topological spaces.
\smallskip
\item"{\rm(b)}" the category of log manifolds $(1.1.5)$.
\endroster
In the case $(a)$ $(resp.\ (b))$, let $H$ be a topological
group $($resp.\ a complex analytic group$)$, $X$ an object of
$\cC$ and assume that we have an action $H\x X\to X$ in $\cC$.
$($In the case {\rm(b)}, we regard $H$ as having the
trivial log structure.$)$
Assume this action is proper $(5.2.4)$ topologically and is free
set-theoretically.
Assume moreover the following condition {\rm(1)} is
satisfied.
\medskip

\noindent
{\rm(1)}
For any point $x\in X$, there exist an object $S$ of
$\cC$, a morphism $\iota:S\to X$ of $\cC$ whose
image contains $x$ and an open neighborhood $U$ of
$1$ in $H$ such that $U\x S\to X$,
$(h,s)\mapsto h\iota(s)$, induces an isomorphism in
$\cC$ from $U\x S$ onto an open set of $X$.
\medskip

Then$:$
\medskip

{\rm(i)} In the case {\rm(b)}, the quotient topological
space $H\bs X$ has a unique structure of an object of
$\cC$ such that, for an open set $V$ of $H\bs X$,
$\cO_{H\bs X}(V)$ $($resp.\ $M_{H\bs X}(V))$ is the set
of all functions on $V$ whose pull-backs to the inverse
image $V'$ of $V$ in $X$ belong to $\cO_X(V')$
$($resp.\ $M_X(V'))$.
\medskip

{\rm(ii)} $X\to H\bs X$ is an $H$-torsor in the
category $\cC$.
\endproclaim

{\it Proof.} This is \cite{KU09} Lemma 7.3.3. See \cite{KU09} for the proof.

\medskip

\proclaim{Proposition 5.2.8} 
$E_\sig^{\sharp} \to i\sig_\bR\bs E_\sig^{\sharp}=D_\sig^{\sharp}$ and  $E_{\sig,\val}^{\sharp} \to i\sig_\bR\bs E_{\sig,\val}^{\sharp}=D_{\sig,\val}^{\sharp}$ are $i\sig_\bR$-torsors in the category of topological spaces.
\endproclaim

{\it Proof.} This is proved by using the arguments in \cite{KU09} 7.3.5 for the pure case. We apply Lemma 5.2.7 by taking $H=i\sig_\bR$, $X=E_\sig^{\sharp}$, and $\cC$ to be the category of topological spaces as in 5.2.7 (a). By Proposition 5.2.5, Proposition 5.2.2 and Lemma 5.2.7, it is sufficient to prove that the condition 5.2.7 (1) is satisfied. Let $x=(q,F)\in E_\sig^{\sharp}$. Let $S_1\sub E_\sig$ be the log manifold containing $x$ constructed in the same way as the log manifold denoted by $S$ in the pure case \cite{KU09} 7.3.5. (In the argument to construct this log manifold, we replace $W(\sig(q))[-w]$ there by $M(\sig(q), W)$ here.)  Let $S_2=E^{\sharp}_\sig \cap S_1$, let $U$ be a sufficiently small neighborhood of $0$ in $\sig_\bR$, and let $S= \{(q', \exp(a)F')\,|\,(q', F')\in S_2, a\in U\}$. Then $S$ has the desired property. 
\qed
\medskip

\head
\S5.3. Study of $D_{\Sig}^{\sharp}$ 
\endhead

Let $\Sig$ be a weak fan. We prove some results on $D_{\Sig}^{\sharp}$.

\medskip

\proclaim{Theorem 5.3.1} For $\sig\in \Sig$, the inclusion maps 
$D_\sig^{\sharp}\to D_\Sig^{\sharp}$ and $D_{\sig,\val}^{\sharp}\to D_{\Sig,\val}^{\sharp}$ are open maps.
\endproclaim

{\it Proof.}
We prove that $D_\sig^{\sharp}\to D_\Sig^{\sharp}$ is an open map. 
For $\sig,\tau\in\Sig$, $D_\sig^\sharp\cap
D_\tau^\sharp$ is the union of $D_{\alpha}^{\sharp}$, where $\alpha$ ranges over all elements of $\Sig$ having the following property:  $\alpha$ is a face of $\sig$ and a face of $\tau$ at the same time. 
Hence, by the definition of the topology of
$D_\Sig^\sharp$, it is sufficient to prove the
following (1).
\medskip

\noindent
(1) If $\sig\in\Sig$ and $\tau$ is a face of $\sig$,
the inclusion $D_\tau^\sharp\to D_\sig^\sharp$ is a
continuous open map.
\medskip

We prove (1).
Let the open set $U(\tau)$ of $\toric_\sig$ be as in 5.1.2,  and let $|U|(\tau)$ be the open set $U(\tau)\cap \abtoric_\sig$ of $\abtoric_\sig$. Then $\abtoric_\tau\sub |U|(\tau)$ as subsets of $\abtoric_\sig$. Let 
$|\t U|(\tau)\sub E_\sig^\sharp$ be the inverse image
of $|U|(\tau)$ under $E_\sig^\sharp\to\abtoric_\sig$.
 We have  commutative diagrams of topological spaces
$$
\CD
E_\tau^\sharp@>\sub>> |\t U|(\tau) @>\sub>>
E_\sig^\sharp\\@VVV@VVV@VVV\\
\abtoric_\tau@>\sub>>|U|(\tau)@>\sub>>\;
\abtoric_\sig,
\endCD\qquad
\CD
E_\tau^\sharp@>\sub>> E_\sig^\sharp\\
@VVV@VVV \\ D_\tau^\sharp@>\sub>>\;D_\sig^\sharp.
\endCD
$$
Furthermore, the inverse image of $D_\tau^\sharp$
under $E_\sig^\sharp\to D_\sig^\sharp$ coincides with
$|\t U|(\tau)$.
The second diagram shows that
$D_\tau^\sharp\to D_\sig^\sharp$ is continuous.
Let $B$ be an $\bR$-vector subspace of $\sig_\bR$
such that $\tau_\bR\op B=\sig_\bR$.
Then we have a homeomorphism
$$
\abtoric_\tau\x B\tra|U|(\tau),\;\;
(a,b)\mapsto\be(ib)a.
$$
   From this, we see that the projection $|\t U|(\tau)\to
D_\tau^\sharp$ factors through $|\t U|(\tau)\to
E_\tau^\sharp$ which sends $(\be(ib)a,F)$ to
$(a,\exp(-ib)F)$.
This shows that a subset $U$ of $D_\tau^\sharp$ is
open in $D_\sig^\sharp$, that is, its inverse image in
$E_\sig^\sharp$ is open, if and only if it is open in
$D_\tau^\sharp$, that is, its inverse image in
$E_\tau^\sharp$ is open.
This completes the proof of (1).

The val version is proved similarly.
\qed
\medskip

\proclaim{Proposition 5.3.2} The map
$D_{\Sig,\val}^{\sharp}\to D_\Sig^{\sharp}$ is proper. 
\endproclaim

{\it Proof.}
By Theorem 5.3.1, we are reduced to the case
$\Sig=\{\text{face of}\;\sig\}$.
In this case, by Proposition 5.2.8, we are reduced to
the fact that $E^\sharp_{\sig,\val}\to E^\sharp_\sig$ is
proper, and hence to the fact that $\abtoric_{\sig,\val}\to \abtoric_\sig$ is 
proper (cf.\ 4.2.5).
\qed
\medskip

\proclaim{Proposition 5.3.3}
 $D_{\Sig}^{\sharp}$ and $D_{\Sig,\val}^{\sharp}$ are Hausdorff spaces.
 \endproclaim
 
 {\it Proof.}
The statement for $D^\sharp_\Sig$ follows from that
for $D_{\Sig,\val}^\sharp$, since $D^\sharp_{\Sig,\val}
\to D^\sharp_\Sig$ is proper by Proposition 5.3.2 and
is surjective.

We prove that $D_{\Sig,\val}^\sharp$ is Hausdorff.
By Theorem 5.3.1, it is sufficient to prove the
following (1).
\medskip

\noindent
(1) Let $\sig$ and $\sig'$ be rational nilpotent cones
and let $\b\in D_{\sig,\val}^\sharp$ and $\b'\in
D_{\sig',\val}^\sharp$.
Assume $x_\lam\in D$ converges to $\b$ in
$D_{\sig,\val}^\sharp$ and to $\b'$ in
$D_{\sig',\val}^\sharp$.
Then $\b=\b'$ in $D_\val^\sharp$.
\medskip

We prove this.
By Proposition 5.2.8, there exist an open neighborhood
$U$ of $\b$ in $D_{\sig,\val}^\sharp$ (resp.\ $U'$ of
$\b'$ in $D_{\sig',\val}^\sharp$) and a continuous
section $s_\sig:U\to E_{\sig,\val}^\sharp$ (resp.\
$s_{\sig'}:U'\to E_{\sig',\val}^\sharp$) of the
projection $E_{\sig,\val}^\sharp\to
D_{\sig,\val}^\sharp$ (resp.\ $E_{\sig',\val}^\sharp\to
D_{\sig',\val}^\sharp$).
Write
$$
s_\sig(x_\lam)=(\be(iy_\lam),F_\lam),\quad
s_{\sig'}(x_\lam)=(\be(iy'_\lam),F'_\lam).
$$
Write $\a=s_\sig(\b)$, $\a'=s_{\sig'}(\b')$.
Then the assumption of Proposition 5.2.3 is satisfied.
Hence we have $\b=\b'$ by Proposition 5.2.3 (i).
\qed
\medskip

\proclaim{Lemma 5.3.4}
Let $V$ be a vector space over a field $K$ and let $N, h: V\to V$ be $K$-linear maps such that $Nh=hN$. 
Let $I$ be an increasing filtration on $V$ with $I_k = 0$ for $k\ll 0$ and $I_k = V$ for $k\gg0$.
Assume that $NI_k\sub I_k$ and $hI_k\sub I_{k-1}$ for all $k$, and that the relative monodromy filtration $M$ of $N$ relative to $I$ exists.
Then $\gr^M(h)=0$.
\endproclaim

{\it Proof.} 
We prove the following statement ($A_j$) for $j\geq 1$, by induction on $j$.
\medskip

\noindent
($A_j$) $hI_k\gr^M_l\sub I_{k-j}\gr^M_l$ for any $k, l\in \bZ$. 

\medskip

First, ($A_1$) holds by the assumption $hI_k\sub I_{k-1}$. Let $j\geq 1$, and assume 
that ($A_j$) holds. 
We prove that ($A_{j+1}$) holds. 
It is sufficient to prove that the map $\gr^I_k\gr^M_l\to \gr^I_{k-j}\gr^M_l$ induced by $h$ is the zero map for any $k, l\in \bZ$. Note $\gr^I_k\gr^M_l= (I_k\cap M_l)/(I_{k-1}\cap M_l + I_k \cap M_{l-1})=\gr^M_l\gr^I_k$. We prove that the map
$\gr^M_l\gr^I_k \to \gr^M_l\gr^I_{k-j}$ induced by $h$ is the zero map.
If $l\geq k$, we denote the kernel of $N^{(l-k)+1}: \gr^M_l\gr^I_k \to \gr^M_{2k-l-2}\gr^I_k$ by $P(\gr^M_l\gr^I_k)$. 
By \cite{D80} 1.6, $\gr^M_l\gr^I_k$ is the sum of the images of 
$N^s: P(\gr^M_{l'}\gr^I_k)\to \gr^M_l\gr^I_k$, where $(s, l')$ ranges over all pairs of integers such that $l'\geq k$, $s\geq 0$, and $l=l'-2s$. 
Hence it is sufficient to prove that the map $P(\gr^M_l\gr^I_k)\to \gr^M_l\gr^I_{k-j}$ induced by $h$ is the zero map for any $k, l\in \bZ$ such that $l\geq k$. Let $x\in \gr^M_l\gr^I_{k-j}$ be an element of the image of $P(\gr^M_l\gr^I_k)$ under this map.  
By the definition of $P(\gr^M_l\gr^I_k)$ and by $Nh=hN$, the map $N^{l-k+1}: \gr^M_l\gr^I_{k-j} \to \gr^M_{l-2(l-k+1)}\gr^I_{k-j}$ kills $x$. 
On the other hand, since $M$ is the monodromy filtration of $N$ relative to $I$, the map $N^{l-(k-j)}: \gr^M_l\gr^I_{k-j} \to \gr^M_{l-2(l-(k-j))}\gr^I_{k-j}$ is an isomorphism. Since $l-(k-j) \geq l-k+1$, we have $x=0$. 
\qed

\medskip

Let $\G$ be a subgroup of $G_\bZ$ which is strongly compatible with $\Sig$.

\proclaim{Theorem 5.3.5}
Assume that $\G$ is neat. 

\medskip

{\rm(i)} Let $p\in D_\Sig^\sharp$, $\g\in \G$,  and assume $\g p=p$. Then $\g=1$.

\medskip

{\rm(ii)} Let $p\in D_\Sig$, $\g\in \G$, 
and assume $\g p=p$.
Then $\g\in\G(\sig)^\gp$.
\endproclaim

{\it Proof.} 
Let $(\sig,Z)\in D_\Sig$ and $\g\in\G$.
Assume $\g(\sig,Z)=(\sig,Z)$, that is, $\Ad(\g)(\sig)=\sig$, $\g Z=Z$. We prove first 

\medskip

\noindent
(1) $\g N=N\g$ for any $N\in \sig_\bC$.

\medskip

Since $\g\G(\sig)\g^{-1}=\G(\sig)$,
 we have an automorphism $\Int(\g):y\mapsto
\g y\g^{-1}$ $(y\in\G(\sig))$ of the sharp fs monoid
$\G(\sig)$.
Since the automorphism group of a sharp fs monoid is finite (see for example, \cite{KU09} 7.4.4; 
an alternative proof: see the generators of the 1-faces), 
this automorphism of $\G(\sig)$ is of finite order.
Since $\sig_\bC$ is generated over $\bC$ by
$\log(\G(\sig))$, the $\bC$-linear map
$\Ad(\g):\sig_\bC\to\sig_\bC$, $y\mapsto\g y\g^{-1}$,
is of finite order.
On the other hand, any eigenvalue $a$ of this
$\bC$-linear map is equal to $bc^{-1}$ for some
eigenvalues $b$, $c$  of the $\bC$-linear map
$\g:H_{0,\bC}\to H_{0,\bC}$, and hence the neat
property of $\G$ shows $a=1$.
Thus we have $\g y\g^{-1}=y$ for any $y\in\sig_\bC$.

Let $M:=M(\sig, W)$.
By (1), we have $\g M=M$.

Take $F\in Z$.
Since $\g Z=Z$, there exists $N\in \sig_\bC$ such that
$$
\g F=\exp(N)F.\tag2
$$

Recall that the pure Hodge theoretic version of 5.3.5 is proved in \cite{KU09} 7.4.5, and its proof shows $\log(\gr^W(\g))=\gr^W(N)$. 
Hence $\gr^W(\g)$ is in the image of $\G(\sig)^\gp$, 
and we may assume $\gr^W(\g)=1$.

Since $\exp(N)$ acts on $\gr^M:=(\gr^M)_\bC$ trivially,
$\gr^M(\g)F(\gr^M)=F(\gr^M)$ follows.

\proclaim{Claim 1}
$\gr^M(\g)=1$.
\endproclaim

This follows from (1) and Lemma 5.3.4 applied to $h=\g-1$.

\medskip

\proclaim{Claim 2}
$\log(\g)M_k\sub M_{k-2}\quad(\forall k\in\bZ)$.
\endproclaim

In fact, by (2), we have
$$
(\log(\g)-N)F^p\sub F^p\quad(\forall p\in\bZ).\tag3
$$
Since $N(M_k)\sub M_{k-2}$, (3) shows that
the map $\gr^M_k\to \gr^M_{k-1}$ induced by $\log(\g)$, which we denote as $\gr^M_{-1}(\log(\g)):\gr^M_k\to\gr^M_{k-1}$, satisfies
$$
\gr^M_{-1}(\log(\g))F^p(\gr^M_k)\sub F^p(\gr^M_{k-1})\quad
(\forall k,\forall p).\tag4
$$
By taking the complex conjugation of (4), we have
$$
\gr^M_{-1}(\log(\g))\Fb^p(\gr^M_k)\sub\Fb^p(\gr^M_{k-1})
\quad(\forall k,\forall p).\tag5
$$
Since
$$
\align
&\gr^M_k=\tOp_{p+q=k}
F^p(\gr^M_k)\cap\Fb^q(\gr^M_k)\quad\text{and}\\
&F^p(\gr^M_{k-1})\cap\Fb^q(\gr^M_{k-1})=0\quad
\text{for $p+q=k>k-1$},
\endalign
$$
(4) and (5) show that the map $\gr^M_{-1}(\log(\g)):\gr^M_k\to\gr^M_{k-1}$ is the zero map.
This proves Claim 2.
\medskip

By (3), $\log(\g)F^p\sub F^{p-1}$, $\log(\g)\Fb^p\sub\Fb^{p-1}$ $(\forall p\in\bZ)$.
These and Claim 2 show, by \cite{KU09} Lemma 6.1.8 (iv),
$$
\log(\g),N\in L^{-1,-1}(M,F).
$$
Hence, by \cite{KU09} Lemma 6.1.8 (iii), $\g F=\exp(N)F$ proves
$\log(\g)=N$.
This proves $\g\in\G(\sig)^\gp$.

If $\g(\sig,Z')=(\sig,Z')$ for some $(\sig,Z')\in
D_\Sig^\sharp$ then, for $Z:=\exp(\sig_\bR)Z'$, we
have $(\sig,Z)\in D_\Sig$ and $\g(\sig,Z)=(\sig,Z)$.
In the above argument, we take $N\in i\sig_\bR$, and
have $\log(\g)=N$.
Since $\log(\g)$ is real and $N$ is purely imaginary,
this shows that $\log(\g)=0$.
Hence $\g=1$.
\qed
\medskip

\proclaim{Theorem 5.3.6}
{\rm(i)} The actions of $\G$ on $D^\sharp_\Sig$ and
$D^\sharp_{\Sig,\val}$ are proper.
In particular, the quotient spaces
$\G\bs D^\sharp_\Sig$ and
$\G\bs D^\sharp_{\Sig,\val}$ are Hausdorff.
\smallskip

{\rm(ii)} Assume $\G$ is neat.
Then the canonical maps $D^\sharp_\Sig\to
\G\bs D^\sharp_\Sig$ and $D^\sharp_{\Sig,\val}\to
\G\bs D^\sharp_{\Sig,\val}$ are local homeomorphisms.
\endproclaim

{\it Proof.} 
We first prove (i).
By \cite{KNU.p} 3.5.17, the action of $\G$ on $D_{\SL(2)}^I$ is proper. 
By 5.2.4.3 (ii), the statement for $D^\sharp_{\Sig,\val}$ follows from
this, from the continuity of $\psi:D^\sharp_{\Sig,\val} \to D_{\SL(2)}^I$, and from Proposition 5.3.3.
By 5.2.4.3 (i), the statement for $D^\sharp_\Sig$ follows from the above result, since $D^\sharp_{\Sig,\val}\to D^\sharp_\Sig$ is proper by
Proposition 5.3.2 and is surjective.

Next, by 5.2.4.2, (ii) follows from (i) and Theorem 5.3.5.
\qed
\medskip

\head 
\S5.4. Study of $\G(\sig)^{\gp}\bs D_{\sig}$
\endhead

Let $\Sig$ be a weak fan, let $\sig\in \Sig$, and let $\G$ be a subgroup of $G_\bZ$ which is strongly compatible with $\Sig$. 

In this subsection, we prove Theorem 2.5.3 and the local cases (the cases where  $\Sig=\text{face}(\sig)$) of Theorems 2.5.2, 2.5.5, 2.5.6. 

\proclaim{Lemma 5.4.1} 
Let $\sig\in \Sig$ and assume that a $\sig$-nilpotent orbit exists.  
Let $\cI$ be the set of all admissible sets $\Psi$ of weight filtrations on $H_{0,\bR}$ $(4.1.2)$ such that $M(\sig,W)\in \Psi$. 
Then, the image of the CKS map
$D_{\sig,\val}^{\sharp}\to D_{\SL(2)}^I$ is contained in the open set $\bigcup_{\Psi \in \cI} D_{\SL(2)}^I(\Psi)$ $(4.1.2)$ of $D_{\SL(2)}^I$. 
\endproclaim

{\it Proof.}
Let $\alpha =(A, V, Z)\in D_{\sig,\val}^{\sharp}$. 
Then there exist $N_1,\dots, N_n\in A\cap \sig$ which generate $A$ over $\bQ$ such that, for the injective homomorphism $(N_j)_{1\leq j\leq n}: \Hom_{\bQ}(A, \bQ)\to \bR^n$, $V$ coincides with the inverse image of the subset of $\bR^n$ consisting of  all elements  which are $\geq 0$ for the lexicographic order. 
Take elements $N_{n+1},\dots, N_m$ of $\fg_\bQ\cap \sig$ such that $N_1,\dots, N_m$ generate the cone $\sig$. 
By assumption, there is an 
$F\in \Dc$ which generates a $\sig$-nilpotent orbit. Let $\Psi:=\{M(N_1+\dots+N_j, W)\;|\;1\leq j \leq m\}\supset \Psi' := \{M(N_1+\dots+N_j, W)\;|\;1\leq j 
\leq n\}$. 
Then $\Psi$ coincides with the set of weight filtrations associated to the class of the $\SL(2)$-orbit associated to $((N_j)_{1\leq j\leq m}, F)$, and $\Psi'$ coincides with the set of weight filtrations associated to the image $p$ of $\alpha$ in $D_{\SL(2)}$. 
We have $p\in D_{\SL(2)}^I(\Psi)$ and $M(\sig, W)=M(N_1+\dots+N_m, W)\in \Psi$. 
 \qed

\medskip

The following lemma will be improved as Proposition 5.4.6 later.

\proclaim{Lemma 5.4.2} Let $\sig\in \Sig$, let $\tau$ be a face of $\sig$ such that a $\tau$-nilpotent orbit exists, let $U(\tau)$ be the open set of $\toric_\sig$ defined in $5.1.2$, and let $\tilde U(\tau)$ $($resp. \
$\tilde U(\tau)_{\val})$ be the inverse image of $U(\tau)$ under $E_{\sig}\to\toric_\sig$ $($resp.\ $E_{\sig,\val}\to \toric_\sig)$. Then the action $(5.2.1)$ of $\sig_\bC$  on $\tilde U(\tau)$ $($resp.\
$\tilde U(\tau)_{\val})$ is proper. 

\endproclaim

\medskip

{\it Proof.} Since $\tilde U(\tau)_{\val}\to \tilde U(\tau)$ is proper and surjective, it is sufficient to consider 
 $\tilde U(\tau)_\val$ (5.2.4.3 (i)).
 
Let $W'$ be the filtration $M(\tau, W)(\gr^W)$ on $\gr^W$ induced by $M(\tau, W)$. We consider five continuous maps
$$f_1: \tilde U(\tau)_{\val} \to E_{\sig,\val}^{\sharp},\quad 
f_2: \tilde U(\tau)_{\val} \to D_{\tau,\val}^{\sharp},$$
$$f_3: \tilde U(\tau)_{\val} \to \spl(W),\quad
f_4: \tilde U(\tau)_{\val}\to \spl(W'),$$
$$f_5: \tilde U(\tau)_{\val}\to \sig_\bR/(\tau_\bR+\log(\G(\sig)^{\gp})).$$
Here $f_1$ is induced by $|-|: E_{\sig,\val}\to E_{\sig,\val}^{\sharp},\;(q,F)\mapsto (|q|, F)$, where for $q: \G(\sig)^{\spcheck}\to \bC^{\mult}$, $|q|$ denotes the composition of $q$ and $\bC\to \bR_{\geq 0},\;a\mapsto |a|$. The map $f_2$ is induced from $f_1$ via the canonical map $E_{\sig,\val}^{\sharp}\to D_{\sig,\val}^{\sharp}$. The map
$f_3$ is induced from $f_2$, the CKS map $D_{\tau,\val}^{\sharp}\to D_{\SL(2)}^I$, and the canonical map $\spl_W: D_{\SL(2)}^I\to \spl(W)$ (\cite{KNU.p} \S3.2). The map $f_4$ is induced from the map $f_2$, the CKS map $D_{\tau,\val}^{\sharp}\to \bigcup_{\Psi\in \cI} D_{\SL(2)}^I(\Psi)$ where $\cI$ denotes the set of all admissible sets of weight filtrations $\Psi$ such that $M(\tau, W) \in \Psi$ (here we apply 5.4.1 replacing $\sig$ there by $\tau$), and the Borel-Serre splitting 
 $\spl_{W'}^{\BS}:\bigcup_{\Psi \in \cI} D_{\SL(2)}^I(\Psi)\to \spl(W')$ (\cite{KNU.p} \S3.2). The map $f_5$ is the composition  
 $$\tilde U(\tau)_{\val}\to U(\tau)\to (\G(\sig)^{\gp}/\G(\tau)^{\gp})\otimes \bC^\times
 \to \sig_\bR/(\tau_\bR+\log(\G(\sig)^{\gp}))$$ where the second arrow is defined by the fact $$\Ker(\Hom(\G(\sig)^{\gp},\bZ) \to \Hom(\G(\tau)^{\gp}, \bZ))\sub P,$$
 where $P$ is the fs monoid which appeared in the definition of $U(\tau)=\Spec(\bC[P])_{\an}$ (5.1.2), and the last arrow is the homomorphism $\g \otimes e^{2\pi iz}\mapsto \text{Re}(z)\log(\g)$ ($\g\in \G(\sig)^{\gp}$, $z\in \bC$). 
 
 These maps $f_j$ have the following compatibilities with the action of $\sig_\bC$ on $\tilde U(\tau)$:
 
 \smallskip
 For any $a\in \sig_\bC$ and any $x\in \tilde U(\tau)_{\val}$,
 $f_j(a\cdot x)= a\cdot f_j(x)$, 
 
 \smallskip
\noindent 
  where $\sig_\bC$ acts on the target space of $f_j$ as follows. 
 For $x=(q, F)\in E_{\sig,\val}^{\sharp}$, $a\cdot x = (\be(i\Im(a))\cdot q,\, \exp(-a)F)$. Note that in the case $a\in i\sig_\bR$, this action coincides with the original action of $a$ on $E_{\sig,\val}^{\sharp}$. For $x=(A,V, Z)\in D_{\tau,\val}^{\sharp}$, $a\cdot x=(A,V, \exp(-\text{Re}(a))Z)$. For $s\in \spl(W)$, $a\cdot s= \exp(-\text{Re}(a))\circ s\circ\gr^W(\exp(\text{Re}(a)))$. 
 For $s\in \spl(W')$, $a\cdot s= \exp(-\text{Re}(a))\circ s\circ\gr^{W'}(\exp(\text{Re}(a)))$. For $x\in \sig_\bR/(\tau_\bR+\log(\G(\sig)^{\gp}))$, $a\cdot x= x+\text{Re}(a)$.
  
  Let $H=\sig_\bC$, and define closed subgroups  $H(j)$ ($0\leq j\leq 5$) of $H$ 
  such that $0=H(0) \sub H(1)\sub \dots \sub H(5)=H$ as follows. 
  $H(1)=i\sig_\bR$, $H(2)= H(1)+\log(\G(\sig)^{\gp})$, $H(3)=H(2)+\tau_{\bR,u}$ where $\tau_{\bR,u}=\tau_\bR\cap \fg_{\bR,u}$, $H(4) = H(3)+\tau_\bR$.
  Define the spaces $X_j$ ($1\leq j\leq 5$) with actions of $H$ as follows: 
  $X_1=E_{\sig,\val}^{\sharp}$, $X_2= D_{\tau, \val}^{\sharp}$,
  $X_3$ is the quotient space of $\spl(W)$ under the action of $\G(\sig)^{\gp}$ given by $s\mapsto \g \circ s \circ \gr^W(\g)^{-1}$ ($s\in \spl(W)$, $\g\in \G(\sig)^{\gp}$), $X_4$ is the quotient space of $\spl(W')$ under the action of $\G(\sig)^{\gp}$ given by $s\mapsto \g \circ s \circ \gr^{W'}(\g)^{-1}$ ($s\in \spl(W')$, $\g\in \G(\sig)^{\gp}$), and  
    $X_5= \sig_\bR/(\tau_\bR+\log(\G(\sig)^{\gp}))$. Then for $1\leq j\leq 5$, the action of $H(j)$ on $X_j$ factors through $H_j:=H(j)/H(j-1)$, and the action of $H_j$ on $X_j$ is proper and free. In fact, for $j=1$, this is by Proposition 5.2.5 and Proposition 5.2.2 (ii). For $j=2$, this is by Theorem 5.3.6 (i) and Theorem 5.3.5 (i). For $j=3,4,5$, this is easily seen. 
  Hence by 5.2.4.5, the action of $H$ on $X_1\times X_2\times X_3\times X_4\times X_5$ is proper. By 5.2.4.3 (ii), this proves that the action of $\sig_\bC$ on $\tilde U(\tau)_{\val}$ is proper.\qed

\medskip

\proclaim{Proposition 5.4.3}
 The log local ringed space $\G(\sig)^{\gp}\bs D_\sig$ over $\bC$ is a log manifold.
\endproclaim

We prove this proposition together with Theorem 2.5.3 (that is, $E_\sig$ is a $\sig_\bC$-torsor over $\G(\sig)^{\gp}\bs D_\sig$).
\medskip

{\it Proof.} Since $E_\sig=\bigcup_{\tau} \tilde U(\tau)$ where $\tau$ ranges over all faces of $\sig$ such that a $\tau$-nilpotent orbit exists, it is sufficient to prove that for such $\tau$, $\G(\sig)^{\gp}\bs D_\tau= \sig_\bC\bs \tilde U(\tau)$ is a log manifold and $\tilde U(\tau)$ is a $\sig_\bC$-torsor over $\G(\sig)^{\gp}\bs D_\tau$. 

For such $\tau$, we apply Lemma 5.2.7 by taking $H=\sig_\bC$,
$X=\tilde U(\tau)$ and $\cC$ to be the category of
log manifolds as in 5.2.7 (b). 
By Lemma 5.4.2, Proposition 5.2.2 (i) and Lemma
5.2.7, it is sufficient to prove that the condition 5.2.7
(1) is satisfied. Let $x=(q,F)\in \tilde U(\tau)$. Let $S_1\sub E_\sig$ be the log manifold which contains $x$ 
constructed in the same way as the log manifold $S$ in the pure case \cite{KU09} 7.3.5. 
 (In the argument to construct this log manifold, we replace $W(\sig(q))[-w]$ there by $M(\sig(q), W)$ here.)  Let $S=S_1\cap \tilde U(\tau)$. Then 
 $S$ has the desired property.
\qed
\medskip

Since $\G(\sig)^{\gp}\bs D_\sig$ belongs to $\cB(\log)$ by 5.4.3, the topological space $(\G(\sig)^{\gp}\bs D_\sig)^{\log}$ is defined.

\proclaim{Proposition 5.4.4} 
We have a canonical homeomorphism 
 $$\G(\sig)^{\gp}\bs D_\sig^{\sharp}\simeq (\G(\sig^{\gp})\bs D_\sig)^{\log}$$ over $\G(\sig)^{\gp}\bs D_\sig$.
\endproclaim

{\it Proof.}
We define a canonical homeomorphism
$(\G(\sig)^\gp\bs D_\sig)^\loga\simeq\G(\sig)^\gp\bs
D_\sig^\sharp$ over $\G(\sig)^{\gp}\bs D_\sig$.
We have
$$
(\toric_\sig)^\loga
\simeq\Hom(\G(\sig)^\v,\bS^1\x\bR_{\ge0}^\mult)
\simeq(\bS^1\ox_\bZ\G(\sig)^\gp)\x\abtoric_\sig.
$$
Hence the homeomorphisms
$$
(E_\sig)^\loga\simeq
(\toric_{\sig})^\loga\x_{\toric_\sig}E_\sig,\quad
E_\sig^\sharp\simeq
\abtoric_\sig\x_{\toric_\sig}E_\sig
$$
induce a homeomorphism
$$
(E_\sig)^\loga\simeq
(\bS^1\ox_\bZ\G(\sig)^\gp)\x E_\sig^\sharp.\tag1
$$
This homeomorphism is compatible with the actions of
$\sig_\bC$.
Here $z=x+iy\in\sig_\bC$ ($x,y\in\sig_\bR$) acts on
the left-hand side of (1) as the map $(E_\sig)^\loga\to
(E_\sig)^\loga$ induced by the action of $z$ on
$E_\sig$ (5.2.1) and on the right-hand side as
$(u,q,F)\mapsto (\be(x)u,\be(iy)q,\exp(-z)F)$
($u\in\bS^1\ox_\bZ\G(\sig)^\gp$, $q\in\abtoric_\sig$,
$F\in\Dc$).
The homeomorphism between the quotient spaces of
these actions is
$$
(\G(\sig)^\gp\bs D_\sig)^\loga\simeq
\G(\sig)^\gp\bs D_\sig^\sharp.\qed
$$
\medskip

\proclaim{Corollary 5.4.5}
The space $\G(\sig)^\gp\bs D_\sig$
is Hausdorff.
\endproclaim

The space $\G(\sig)^{\gp}\bs D_\sig^{\sharp}$ is Hausdorff by 5.3.6 (i) 
and the map
$\G(\sig)^{\gp}\bs D_\sig^{\sharp}\to \G(\sig)^{\gp}\bs D_{\sig}$ is proper and surjective by 5.4.4. Hence $\G(\sig)^\gp\bs D_\sig$
is Hausdorff by 5.2.4.7.

\proclaim{Proposition 5.4.6} 
The action of $\sig_\bC$ on $E_\sig$ $($resp.\ $E_{\sig,\val})$ is proper.
\endproclaim

The result for $E_\sig$ follows from Theorem 2.5.3 and Corollary 5.4.5 by 5.2.4.6. The result for $E_{\sig,\val}$ follows from it by 5.2.4.3 (ii).
%%%%%

\head
\S5.5. Study of $\G \bs D_{\Sig}$
\endhead

In this subsection, we prove Theorems 2.5.2, 2.5.4, 2.5.5, 2.5.6, and 2.6.6.
\proclaim{Lemma 5.5.1}
Let $X$ be a topological space with a continuous action of a discrete group $\G$ and let $Y$ be a set with an action of $\G$.
Let $f:X\to Y$ be a $\G$-equivariant surjective map.
Let $\G'$ be a subgroup of $\G$.
We introduce the quotient topologies of $X$ on $\G'\bs Y$ and on $\G\bs Y$.
Let $V$ be an open set of $\G'\bs Y$ and let $U$ be the inverse image of $V$ in $\G'\bs X$.
We assume moreover the three conditions {\rm(1)--(3)} below.
Then, $V\to\G\bs Y$ is a local homeomorphism.
\medskip

\noindent
{\rm (1)} $X\to\G\bs X$ is a local homeomorphism and
$\G\bs X$ is Hausdorff.
\smallskip

\noindent
{\rm (2)} $U\to V$ is proper.
\smallskip

\noindent
{\rm (3)} If $x\in X$ and $\g\in\G$, and if the images
of $\g x$ and $x$ in $\G'\bs Y$ are contained in $V$ and
they coincide, then $\g\in\G'$.
\endproclaim

{\it Proof.} 
This is \cite{KU09} Lemma 7.4.7. See \cite{KU09} for a proof.

\medskip

{\bf 5.5.2.}
{\it Proof of Theorem 2.5.4.} We prove that $\G(\sig)^{\gp}\bs D_\sig \to \G \bs D_\Sig$ is locally an isomorphism. 

We use Lemma 5.5.1 for $X=D_\Sig^\sharp$, $Y=D_\Sig$, $\G=\G$, $\G'=\G(\sig)^\gp$, $V=\G(\sig)^\gp\bs D_\sig$ and $U=\G(\sig)^\gp\bs D_\sig^\sharp$.
Theorem 5.3.6 for $D^\sharp_\Sig$ shows that the condition (1) in Lemma 5.5.1 is satisfied.
Theorem 2.5.6 for $\G(\sig)^\gp\bs D_\sig$, proved in \S5.4, shows that the condition (2) in 5.5.1 is satisfied.
Theorem 5.3.5 shows that the condition (3) in 5.5.1 is satisfied.
\qed
\medskip

{\bf 5.5.3.}
{\it Proof of Theorem 2.5.2.} 
We prove that $\G \bs D_\Sig$ is a log manifold.

This follows from Theorem 2.5.2 for
$\G(\sig)^\gp\bs D_\sig$ (5.4.3) by Theorem 2.5.4.
\qed
\medskip

{\bf 5.5.4.}
{\it Proof of Theorem 2.5.6.} We prove that $(\G \bs D_\Sig)^{\log} = \G \bs D_\Sig^{\sharp}$. 

By using Theorem 2.5.4, it is easily seen that the
canonical homeomorphisms $(\G(\sig)^\gp\bs
D_\sig)^\loga\simeq\G(\sig)^\gp\bs D_\sig^\sharp$
(Theorem 2.5.6 for $\G(\sig)^\gp\bs D_\sig$, proved in
\S5.4) glue uniquely to a homeomorphism
$(\G\bs D_\Sig)^\loga\simeq\G\bs D_\Sig^\sharp$.
\qed
\medskip

{\bf 5.5.5.} {\it Proof of Theorem 2.5.5.} 
We prove that $\G\bs D_\Sig$ is a Hausdorff space.

Replacing $\G$ by a subgroup of finite index, we may assume $\G$ is neat.
Then, by Theorem 2.5.6, the map $\G\bs D_\Sig^\sharp\to\G\bs D_\Sig$ is proper and surjective.
Since $\G\bs D_\Sig^\sharp$ is Hausdorff by Theorem 5.3.6 (i), it follows by 5.2.4.7 
that $\G\bs D_\Sig$ is Hausdorff.
\qed
\medskip

A local description of the log manifold $\G\bs D_\Sig$
is given by the following.

\proclaim{Theorem 5.5.6} Let $\Sig$ be a weak fan and let $\G$ be a neat subgroup of $G_\bZ$ which is strongly compatible with $\Sig$. 
Let $(\sig,Z)\in D_\Sig$ and let $p$ be its image in $\G\bs D_\Sig$.
\medskip

{\rm(i)} Let $F_{(0)}\in Z$.
Then there exists a locally closed analytic submanifold
$Y$ of $\Dc$ satisfying the following conditions {\rm(1)} and {\rm(2)}.
\smallskip

\noindent
{\rm{(1)}} $F_{(0)}\in Y$.
\smallskip

\noindent
{\rm{(2)}} The canonical map $T_Y(F_{(0)})\op\sig_\bC
\to T_{\Dc}(F_{(0)})$ is an isomorphism. Here $T_Y(F_{(0)})$ denotes the tangent space of $Y$ at $F_{(0)}$ and $T_{\Dc}(F_{(0)})$ denotes the tangent space of $\Dc$ at $F_{(0)}$. 
\medskip

{\rm(ii)} Let $F_{(0)}$ and $Y$ be as in {\rm(i)}.
Let $X:=\toric_\sig\x Y$, and let $S$ be the log manifold 
defined by
$$
S:=\{(q,F)\in X\;|\; N(F^p)\sub F^{p-1}\,
(\forall N\in\sig(q),\forall p\in\bZ)\},
$$
where $\sig(q)$ denotes the face of $\sig$ corresponding to $q$. 
Then, there exist an open neighborhood $U_1$ of
$(0,F_{(0)})$ in $S$ in the strong topology, an open
neighborhood $U_2$ of $p$ in $\G\bs D_\Sig$, and an
isomorphism $U_1\tra U_2$ of log manifolds
sending $(0,F_{(0)})$ to $p$.
\endproclaim

{\it Proof.}
This is the mixed Hodge version of \cite{KU09} 7.4.13 and proved similarly.
\qed

\medskip

{\bf 5.5.7.}
We prove Theorem 2.6.6. 
This is the mixed Hodge theoretic version of Theorem B in \cite{KU09}.
Since everything is prepared in \S2.5, especially Theorem 2.5.3,  the proof of 2.6.6  is exactly parallel to the pure case \cite{KU09} 8.2. 
That is, we first describe the functor represented by $E_{\sig}$. 
Second we describe the functor represented by $\G(\sig)^{\gp} \bs D_{\sig}$ by taking quotients where we use 2.5.3.
Finally, gluing them up, we show that the functor $\operatorname{LMH}_\Phi$ is represented by $\G\bs D_\Sig$ (cf.\ 2.6.5). 

\bigskip

\head
\S5.6. Proofs for Section 3
\endhead

Let the notation be as in \S3. 

\medskip

{\bf 5.6.1.} 
We reduce the results in \S3 to the results in \S2 in the following way. 
We define a modification $\Lam^e=(H_0^e, W^e, (\langle\;,\;\rangle^e_w)_w, (h^{e,p,q})_{p,q})$ of $\Lam=(H_0, W, (\langle\;,\;\rangle_w)_w, (h^{p,q})_{p,q})$ of 3.1.1 (1).  
Then, the results in \S3 are reduced to the results in \S2 for $\Lam^e$.

Take an integer $r$ such that $W_{2r}=0$.  
Let $H_{0,(w)}^e$ be $\G'$ if $w=2r-2$, 
be $\bZ$ if $w=2r$, and be $H_{0,(w)}$ if $w \neq 2r-2, 2r$. 
Let $H_0^e= \bigoplus_{w\in \bZ} H_{0,(w)}^e= \G' \oplus \bZ\oplus H_0$. 
Let the increasing filtration $W^e$ on $H_{0,\bR}^e$ be the evident one. 

Define the bilinear form $\langle\;,\;\rangle_w^e: H_{0,(w),\bR}^e\times H_{0,(w),\bR}^e \to \bR$ as follows. 
If $w \neq 2r-2, 2r$, it is just the bilinear form $\langle\;,\;\rangle_w$ on $H_{0,(w)}$. 
If $w=2r-2$, it is any fixed positive definite symmetric bilinear form $\bR \otimes \G' \x \bR \otimes \G' \to \bR$ which is rational over $\bQ$. 
If $w=2r$, it is the form $\bR\times \bR\to \bR,\;(x,y)\mapsto xy$. 
  
Let $h^{e, p,q}$ be $\text{rank}(\G') $ if $p=q=r-1$, be $1$ if $p=q=r$, and be $h^{p,q}$ otherwise.

We denote by $D^e$ the classifying space in \S2 for $\Lam^e=(H_0^e, W^e, (\langle\;,\;\rangle_w^e)_w, (h^{e,p,q})_{p,q})$. 

For a nilpotent cone $\sig$ in \S3, we regard $\sig$  as a nilpotent cone in \S2 for $\Lam^e$ in the following way. 
For $N\in \sig$, we regard it as a linear map $H^e_{0,\bR}\to H^e_{0,\bR}$ whose restriction to $H_{0,\bR}$ is  $N_{\fg}: H_{0,\bR}\to H_{0,\bR}$, whose restriction to $\bR \otimes \G'=H_{0,(2r-2),\bR}^e$ is $0$, and which sends  $1\in \bR=H_{0,(2r)}^e$ to the image of $N$ in $\bR \otimes \G'= H_{0,(2r-2),\bR}^e$ under $\sig\to \sig'\sub \sig'_\bR \simeq \bR \otimes \G'$ (3.1.1).

By this, a weak fan $\Sig$ in \S3 is regarded as a weak fan in \S2 for $\Lam^e$.

We denote  the objects $E_{\sig}$, $G_\bZ$  etc.\ of \S2 for $\Lam^e$, by  $E_{\sig}^e$, $G^e_\bZ$, etc.

$\G$ in $\S3$ is regarded as a subgroup of $G^e_{\bZ}$ in the following way. We regard $(a,b)\in \G$ ($a\in \G'$, $b\in G_{\bZ}$) as an element of $G^e_\bZ$ whose restriction to $H_0$ is $b$, whose restriction to $\G'=H_{0,(2r-2)}^e$ is the identity map, and which sends $1\in \bZ=H_{0,(2r)}^e$ to the sum of $1\in H_{0,(2r)}^e$ and the element $a$ of $\G'=H_{0,(2r-2)}^e$. $\Sig$ and $\G$ are compatible (resp.\ strongly compatible) in the sense of \S3 if and only if they are compatible (resp.\ strongly compatible) in the sense of \S2 for $\Lam^e$. The notation $\G(\sig)$ in \S3 and that in \S2 for $\Lam^e$ are compatible. 

We define a morphism $\sig'_\bC \times \Dc \to \Dc^e$ as follows. 
An element $a\in \sig'_\bC$ corresponds to a decreasing filtration $f(a)$ on $V_\bC:= \bC \otimes \G' \oplus \bC = H_{0,(2r-2),\bC}^e \oplus H_{0,(2r),\bC}^e$ given by 
$$
f(a)^{r+1}=0 \sub f(a)^r= \bC \cdot (a,1) \sub  f(a)^{r-1}= V_\bC, 
$$ 
where $a$ in $(a,1)$ is identified with an element of $\bC \otimes \G'$ via the canonical isomorphism $\sig'_\bC \simeq \bC \otimes \G'$. 
The map $\sig'_\bC \times \Dc\to \Dc^e$ sends an element $(a, F)$ of  $\sig'_\bC \times \Dc$  to the element $f(a) \oplus F$ of $\Dc^e$ defined as the direct sum of $f(a)$ and $F$ on $V_\bC\oplus H_{0,\bC} =H_{0,\bC}^e$. 

Note that, via the above morphism, $\sig'_\bC \times \Dc$ is a closed analytic submanifold of $\Dc^e$. 

\proclaim{Proposition 5.6.2} Let $S_0=E'_{\sig'}$ as in $3.1.1$--$3.1.3$. 
Then the following diagram is cartesian in $\cB(\log)$. 
$$
\CD
E_{S_0,\sig} @>>> E_{\sig}^e \\
@VVV@VVV \\
 \sig'_\bC \times \Dc @>>> \;\Dc^e.
 \endCD
$$
Here the lower horizontal arrow sends $(z, F)$ 
$(z\in \sig'_\bC, F\in \Dc)$ to $f(z)\oplus F$ as in $5.6.1$, the left vertical arrow sends $(s,z,q,F)$ to $(z, F)$ $(s\in S_0, z\in \sig'_\bC$, $q\in \toric_\sig$, $F\in \Dc)$, the right vertical arrow sends $(q, F)$ to $F$ $(q\in \toric_\sig$, $F\in \Dc^e)$, and the upper horizontal arrow sends 
$(s, z, q, F)$ to $(q, f(z)\oplus F)$. 
\endproclaim

This is seen easily.

\medskip

{\bf 5.6.3.} 
We give preliminaries to Propositions 5.6.4, 5.6.5.

Let $X$, $Y$ be topological spaces, let $r:X\to Y$ be a continuous map, and assume that, locally on $Y$, there is a continuous map $s: Y\to X$  such that $r\circ s$ is the identity map of $Y$. Let $B$ be a subspace of $Y$ and let $A=r^{-1}(B)$.

\medskip

{\bf 5.6.3.1.} Endow $A$ with the topology as a subspace of $X$. Then 
as is easily seen, the following two topologies on $B$ coincide. 
The image of the topology of $A$ under $r:A\to B$, and the inverse image of the topology of $Y$.  

\medskip

{\bf 5.6.3.2.}
Assume furthermore that $X$ (resp.\ $Y$) is endowed with a subsheaf $\cO_X$ (resp.\ $\cO_Y$) of rings  over  $\bC$ of the sheaf of $\bC$-valued continuous functions on $X$ (resp.\ $Y$). 
Assume that the stalks of $\cO_X$ (resp.\ $\cO_Y$) are local rings. Assume that $r$ respects these sheaves and assume that, locally on $Y$, we can take $s$ as above which respects these sheaves.

Let $\cO_A$ be the subsheaf of the sheaf of $\bC$-valued continuous functions on $A$ defined as follows. 
For a $\bC$-valued continuous function $f$ on an open set $U$ of $A$, $f$ belongs to $\cO_A$ if and only if, for any $x\in U$, there are an open neighborhood $U'$ of $x$ in $U$, an open neighborhood $V$ of $x$ in $X$ such that $U'\sub V$, and an element $g$ of $\cO_X(V)$  such that the restrictions of $f$ and $g$ to $U'$ coincide. 

Then it is easy to see the coincidence of the following two subsheaves (1) and (2) of the sheaf of $\bC$-valued continuous functions on $B$: Let $f$ be  
 a $\bC$-valued continuous function on an open set $U$ of $B$. $f$ belongs to the sheaf (1) if and only if $f\circ r$ belongs to $\cO_A(r^{-1}(U))$. 
$f$ belongs to the sheaf (2) if and only if, for any $y\in U$, there are an open neighborhood $U'$ of $y$ in $U$, an open neighborhood $V$ of $y$ in $Y$ such that $U'\sub V$, and an element $g$ of $\cO_Y(V)$  such that the restrictions of $f$ and $g$ to $U'$ coincide. 

Let $\cO_B$ be this sheaf on $B$.

\medskip

{\bf 5.6.3.3.}
Assume further that we have elements $f_1,\dots, f_n$ of $\cO_X(X)$ such that $A=\{x\in X\;|\;f_1(x)=\dots = f_n(x)=0\}$ and such that the kernel of $\cO_X|_A\to \cO_A$ is the ideal generated by $f_1,\dots, f_n$. 
Here $\cO_X|_A$ is the sheaf-theoretic inverse image of $\cO_X$ by $A \hra X$.
Then locally on $Y$, for $s$ as above, $B=\{y\in Y\;|\;f_1s(y)=\dots = f_ns(y)=0\}$ and the kernel of $\cO_Y|_B\to \cO_B$ is the ideal generated by $f_1\circ s,\dots, f_n\circ s$. 

\medskip

{\bf 5.6.3.4.}
Assume further that $Y$ is endowed with a log structure $M_Y$. Let $M_X$ (resp.\ $M_A$) be the inverse image of $M_Y$ on $X$ (resp.\ $A$) and assume that the structural map $M_A\to \cO_A$ of log structure is injective. 

Then the following two log structures (1) and (2) on $B$ coincide. 
(1) For an open set $U$ of $B$, $M_B(U)$ is the subset of $\cO_B(U)$ consisting of all elements whose pull back to $r^{-1}(U)$ belongs to the image of $M_A\to \cO_A$. 
(2) The inverse image of $M_Y$ on $B$.

\medskip

{\bf 5.6.3.5.}
Assume furthermore that $X$ and $Y$ belong to $\cB(\log)$. Let $H$ be a complex Lie group. Assume that 
an action of $H$ on $X$ in $\cB(\log)$ is given, and assume that $X$ is an $H$-torsor over $Y$ in $\cB(\log)$. Note that $A$ belongs to $\cB(\log)$. By the above, $B$ also belongs to $\cB(\log)$.  Assume that $H$ acts also on $A$ in $\cB(\log)$ and assume that this action is compatible with the action of $H$ on $X$. Then 
$A$ is an $H$-torsor over $B$ in $\cB(\log)$. 
\medskip

{\it Proof} of 5.6.3.5.
Considering locally on $B$, we may assume that we have $s: Y\to X$ as above. Hence we have $H\times Y \to X, \;(h, y) \mapsto hs(y)$. Let $t:X\to H$ be the composition $X\simeq H \times Y \to H$. Then we have a morphism $A \to H \times B, \;a\mapsto (t(a), r(a))$ which is the inverse map of the morphism $H \times B \to A, \; (h, b)\mapsto hs(b)$. These morphisms are inverse to each other in $\cB(\log)$. This proves that $A$ is an $H$-torsor over $B$ in $\cB(\log)$. 
\qed
\medskip

It follows that the diagram
$$
\CD 
A @>>> B\\
@VVV@VVV \\ 
X @>>> Y
\endCD
$$
in $\cB(\log)$ is cartesian. 
\medskip

\proclaim{Proposition 5.6.4}  
{\rm (i)} $\G(\sig)^{\gp}\bs D_{S_0, \sig}$ is a relative log manifold over $S_0$. 

\medskip

{\rm (ii)} $E_{S_0,\sig}$ is a relative log manifold over $S_0$. 
\medskip

{\rm (iii)} $E_{S_0, \sig}\to  \G(\sig)^{\gp}\bs D_{S_0, \sig}$ is a $\sig_\bC$-torsor in the category $\cB(\log)$. 
\endproclaim
{\it Proof.} 
We apply 5.6.3 by taking 
$$
X=E_\sig^e, \quad Y=\G(\sig)^{\gp}\bs D_\sig^e, \quad 
A=E_{S_0,\sig}, \quad B=\G(\sig)^{\gp}\bs D_{S_0,\sig}.
$$ 
By 5.6.3, $\G(\sig)^{\gp} \bs 
D_{S_0,\sig}\in \cB(\log)$ and we have (iii). By (iii), (i) and (ii) are equivalent. It is sufficient to prove that (ii). 
 This is reduced to the fact that the morphism $$\sig'_{\bC}\times \toric_{\sig}\to \toric_{\sig'},\quad (z, q) \mapsto \be(z)\cdot q(\gr^W)$$ is log smooth, which can be shown easily. 
 \qed
\medskip

In the following, $\Sig$ is a weak fan in the sense of \S3.

\proclaim{Proposition 5.6.5} 
{\rm (i)} $\G\bs D_{S_0,\Sig}$ is a relative log manifold over $S_0$. 

\medskip

{\rm (ii)} The morphism $ \G(\sig)^{\gp}\bs D_{S_0,\sig} \to  \G\bs D_{S_0, \Sig}$ in $\cB(\log)$ is locally an  isomorphism.

\medskip

{\rm (iii)} $\G \bs D_{S_0, \Sig}$ is Hausdorff. 

\endproclaim

{\it Proof.} 
By the next lemma, we can apply 5.6.3 by taking
$$
X= \tsize\bigsqcup_{\sig\in \Sig} \G(\sig)^{\gp}\bs D_{\sig}^e, \quad Y= \G \bs D_{\Sig}^e, \quad A= \tsize\bigsqcup_{\sig\in \Sig} \G(\sig)^{\gp}\bs D_{S_0,\sig},\quad B= \G \bs D_{S_0, \Sig}.
$$ 
By 5.6.3.1--5.6.3.4, $\G \bs D_{S_0,\Sig}$ belongs to $\cB(\log)$ and we have (ii) by 2.5.4. 
(ii) and 5.6.4 (i) imply (i). 
(iii) follows from the fact that $\G\bs D_\Sig^e$ is Hausdorff, proved in \S2.
\qed

\medskip

\proclaim{Lemma 5.6.6} 
For $\sig\in \Sig$, the following diagrams of sets are cartesian.
$$
\CD
 D_{S,\sig} @>>> D_{S, \Sig}\\
@VVV@VVV\\
 D_{\sig}^e @>>>  \;D_{\Sig}^e,
\endCD  \qquad\qquad 
\CD
\G(\sig)^{\gp}\bs D_{S,\sig} @>>> \G\bs D_{S, \Sig}\\
@VVV@VVV\\
\G(\sig)^{\gp}\bs D_{\sig}^e @>>> \;\G \bs D_{\Sig}^e.
\endCD
$$
\endproclaim

Recall that we endowed $\G \bs D_{S,\Sig}$ with the structure of an object of $\cB(\log)$ by regarding it as the fiber product of $S\to S_0 \leftarrow \G\bs D_{S_0,\Sig}$, where $S \to S_0$ is a strict morphism in $\cB(\log)$ as in 3.1.1 (4). Hence Theorem 3.2.8 follows from 
Proposition 5.6.5.

\medskip

{\bf 5.6.7.} 
We prove Theorem 3.3.3.
By \S2 for $\Lam^e$, $ \G \bs D_\Sig^e$ represents the moduli functor ${\text {\rm LMH}}_{\Phi^e}$ for $\Phi^e=(\Lam^e, \Sig, \G)$. From 5.6.2, we have that
the subspace $\G \bs D_{S_0, \Sig}$ of  
$ \G \bs D_\Sig^e$ represents the subfunctor of ${\text{\rm LMH}}_{\Phi^e}$  consisting of the classes of LMH which are the direct sum of an LMH $H$ of weight $>2r$ and an LMH $h$ of weight $\leq 2r$. 
 The projection $E_{S_0,\sig}\to
 S_0\sub \toric_{\sig'}\times \Dc(\gr^W)$ represents $h\oplus H\mapsto (h, H(\gr^W))$. 
 Hence over $S_0$, $\G\bs D_{S_0,\Sig}$ is the moduli space of this $H$ with given $H(\gr^W)$, and hence 
  represents 
 ${\text {\rm LMH}}_{Q, \G}^{(\Sig)}$. 
\qed

\bigskip

\head
\S6. N\'eron models
\endhead

In this section, we construct the N\'eron model and the connected N\'eron model in degeneration of Hodge structures. In the case of Hodge structures corresponding to abelian varieties of semi-stable reductions, these coincide with the classical N\'eron model and the connected N\'eron model in the theory of degeneration of abelian varieties. 
We expect that the constructions of such models are useful in the study of degeneration of intermediate Jacobians. 
In the case of degeneration of intermediate Jacobians of Griffiths 
\cite{G68a} and \cite{G68b}, our connected N\'eron model coincides with the model of Zucker \cite{Z76} (which need not be Hausdorff) modified by putting slits (then we have a Hausdorff space), and under a certain assumption, our N\'eron model coincides with the model constructed by Clemens \cite{C83}.

In \S6.4, we give the 
 proofs of the results in \cite{KNU10c} which were 
omitted there.

\head
\S6.1. Results on N\'eron models and connected N\'eron models
\endhead

In \S6.1--\S6.3, as in \S3, let $S$ be an object of $\cB(\log)$, and assume 
that for each $w\in \bZ$, we are given a polarized log Hodge structure $H_{(w)}$ on $S$. Assume $H_{(w)}=0$ for almost all $w$. 

\proclaim{Theorem 6.1.1} 
There is a relative log manifold $J_1$ $($resp.\ $J_0)$ over $S$ which is strict over $S$ and which represents the following functor on $\cB/S^{\circ}$\rom: 

\smallskip
 $\cB/S^{\circ}\ni S' \mapsto $ the set of all isomorphism classes of 
LMH H on $S'$ with $H(\gr^W_w)=H_{(w)}$ $(w\in \bZ)$ satisfying the following condition\rom: 
Locally on $S'$, there is a splitting $H_\bQ(\gr^W)\simeq H_\bQ$ $($resp.\ $H_\bZ(\gr^W)\simeq H_\bZ)$ of the weight filtration $W$ on the local system $H_\bQ$ $($resp.\ $H_\bZ)$ on $(S')^{\log}$. 
\endproclaim

Here $S^{\circ}$ denotes $S$ regarded as an object of $\cB$ by 
forgetting the log structure. 
In the description of the functor in 6.1.1, $S'\in \cB/S^{\circ}$ is endowed with the inverse image of the log structure of $S$. 

\medskip

{\bf 6.1.2.} 
We call the relative log manifold $J_1$ the {\it N\'eron model}, and $J_0$  the {\it connected N\'eron model}. 

In the case of degeneration of abelian varieties, these terminologies agree with the classical ones (see 7.1.4). 

\medskip

{\bf 6.1.3.} 
Let $w\in \bZ$, $w<0$. 
Assume $H_{(k)}=0$ for $k\neq w, 0$, and that $H_{(0)}=\bZ$. 
Let $H'=H_{(w)}$. 
Consider the exact sequence on $S^\loga$ 
$$
0 \to H'_\bZ \to H'_{\cO^\loga}/F^0 \to 
H'_\bZ \bs H'_{\cO^\loga}/F^0 \to 0.
$$
Descending this by $\tau : S^\loga \to S$, we have an exact sequence of sheaves of abelian groups on $\cB(\log)/S$
$$
0 \to \tau_*H'_\bZ\bs (H'_\cO/F^0H'_\cO)^{sG} \to {\cE}xt^1(\bZ, H')\to R^1\tau_*H'_\bZ.
$$ 
Here $(\;)^{sG}$ is the small Griffiths part, i.e., the part restricted by the small Griffiths transversality (cf. [KU09] 2.4.9), and ${\cE}xt^1(\bZ, H'):={\cE}xt^1_{\LMH/S}(\bZ, H')$ is the Ext-sheaf of the category of LMH which is identified with  the subgroup of $\tau_*(H'_\bZ \bs H'_{\cO^\loga}/F^0)^{sG}$ restricted by the condition of the admissibility of local monodromy. 
\medskip

\proclaim{Proposition 6.1.4} 
In the situation of $6.1.3${\rm :}

{\rm(i)} On $\cB/S^{\circ}$, the N\'eron model represents 
$$
\Ker({\Cal E}xt^1(\bZ, H')\to R^1\tau_*H'_\bZ \to R^1\tau_*H'_\bQ).
$$

\medskip

{\rm(ii)} On $\cB/S^{\circ}$, the connected N\'eron model represents 
$$
\Ker({\cE}xt^1(\bZ, H')\to R^1\tau_*H'_\bZ)
= \tau_*H'_\bZ\bs (H'_\cO/F^0H'_\cO)^{sG}.
$$

\endproclaim

{\bf 6.1.5.} 
{\it Remark 1.} 
6.1.4 (ii) shows that in the situation of 6.1.3, the connected N\'eron model is the Zucker model (\cite{Z76}) modified by putting slits induced by the small Griffiths transversality.

\medskip
{\it Remark 2.} 
Proposition 6.1.4 shows that in the situation of 6.1.3,  the N\'eron model and the connected N\'eron model are commutative group objects over $S$ which represent subgroup functors of
${\Cal E}xt^1(\bZ,H')$ on $\cB/S^\circ$. 
\medskip

For an object $S$ of $\cB(\log)$ and for an integer $r$, we say the log rank of $S$ is $\leq r$ if $\rank_{\bZ}((M_S^{\gp}/\cO_S^\times)_s)\leq r$ for any $s\in S$. As an example which appears in 6.1.7 below, if $S$ is a complex analytic manifold with a smooth divisor which gives the log structure of $S$, then the log rank of $S$ is $\leq 1$. 

\proclaim{Corollary 6.1.6} 
In the situation of $6.1.3$, assume furthermore $w=-1$ and the log rank of $S$ is $\leq 1$. 
Then, on $\cB/S^{\circ}$, the N\'eron model represents ${\Cal E}xt^1(\bZ, H')$.
\endproclaim

This is because under the assumption of 6.1.6, $W$ on the local system $H_\bQ$ of any LMH $H$ splits locally on $S$. 

\proclaim{Corollary 6.1.7} 
In the situation of $6.1.6$, assume that $S$ is a complex analytic manifold with a smooth divisor which gives the log structure of $S$. 
Then we have
$$
\align
&(\text{a section of the N\'eron model over $S$})\\
= \;&(\text{a normal function on $S^*$ which is admissible with respect to $S$}).
\endalign
$$

\endproclaim
(Cf. \cite{Sa96} for admissible normal functions.) 

\medskip

{\bf 6.1.8.} In what follows, the proof of 6.1.1 is given in \S6.2, and the proof of 6.1.4 is given in \S6.3. 
Note that for the proofs of 6.1.1 and 6.1.4, we may work locally on $S$. 
Locally on $S$, we construct in \S6.2 the N\'eron model and the connected N\'eron model by the method of \S3, by constructing suitable weak fans explicitly.

\vskip20pt

\head
\S6.2. Moduli spaces of LMH with $\bQ$-splitting local monodromy
\endhead

Let $S$ and $(H_{(w)})_w$ be as in 6.1.
We assume that we are given the data (1)--(5) in 3.1.1. (Note that in general, 
as we proved in 3.1.2, we have always such (1)--(5) locally on $S$.)

Let $\G_u=G_{\bZ,u}$ and let $\G$ be the corresponding group (3.1.7).

\medskip

In this subsection, we construct certain weak fans $\Sig_1$ and $\Sig_0$. 
We will see that $J_1$ in Theorem 6.1.1 is $J_{\Sig_1}$ and $J_0$ is $J_{\Sig_0}$ (3.3.2).
\medskip

\proclaim{Theorem 6.2.1} 
Let $ \sig'=\Hom(P, \bR_{\geq 0}^{\add})\to \fg_\bR(\gr^W)\to \fg_\bR$,  $a\mapsto a_{\fg}$, be the composite map as in $\S3$. 
Let $\Upsilon$ be a subgroup of $G_{\bQ,u}$. 
Let $\Sig(\Upsilon)$ be the set of cones $\sig_{\tau',\upsilon}:=\{(x, \Ad(\upsilon)x_{\fg})\;|\; x\in \tau'\}\subset \sig' \times \fg_\bR$, where $\tau'$ ranges over all faces of $\sig'$ and $\upsilon$ ranges over all elements of $\Upsilon$. 
Then $\Sig(\Upsilon)$ is a weak fan. 
\endproclaim

{\it Proof.} 
By construction, it is easy to see the rationality and the closedness under the operation of taking a face.
We examine the condition (1) in 3.1.6.

Let $\tau', \rho'\in \face(\sig')$ and let $\tau_\fg, \rho_\fg$ be their images in $\fg_\bR$, respectively. 
Let $\upsilon \in \Upsilon$ and let $\rho_{\fg,\upsilon}=\Ad(\upsilon)(\rho_\fg)$. 
Note that $\rho_{\fg,\upsilon}$ coincides with the image of $\sig_{\rho', \upsilon}$ in $\fg_\bR$.

Assume that $N$ is in the interior of $\sig_{\tau',1}$ and of $\sig_{\rho',\upsilon}$. 
Let $F\in \Dc$. 
Assume that $(\sig_{\tau',1}, F)$ and $(\sig_{\rho', \upsilon}, F)$ generate nilpotent orbits. 

We first claim that $N_\fg\upsilon =\upsilon N_\fg$.
In fact, $\gr^W(N_\fg) = \gr^W(\Ad(\upsilon)^{-1}N_\fg)$ because $\upsilon \in \Upsilon \sub G_{\bQ,u}$.
On the other hand, 
$N_\fg\in \tau_\fg$ and $\Ad(\upsilon)^{-1}N_\fg\in \rho_\fg$ are pure of weight $0$ with respect to $W$ under the fixed splitting of $W$. 
Hence $N_\fg = \Ad(\upsilon)^{-1}N_\fg$, that is, $N_\fg\upsilon =\upsilon N_\fg$.

Since $\face(\sig')$ is a fan and the image of $N$ in $\sig'$ is in the interior of $\tau'$ and of $\rho'$, we see $\tau'=\rho'$.

Since $(\sig_{\tau',1}, F)$ generates a nilpotent orbit, 
the action of $\tau_\fg$ on $H_{0,\bR}$ is admissible with respect to $W$. By 2.2.8, the adjoint action of $\tau_{\fg}$ on $\fg_\bR$ is admissible with respect to the filtration $W\fg_\bR$ on $\fg_\bR$ induced by $W$. 

Let $M=M(\Ad(N_\fg), W\fg_\bR)$.
Since $\upsilon - 1 \in (W\fg_\bR)_{-1}$ and $\Ad(N_\fg)(\upsilon - 1) = 0$, we see $\upsilon-1\in M_{-1}$ by 1.2.1.3.

Let $h\in \tau_\fg$.
Then $h \in M_{-2}$.
Hence, by the above result, $h\upsilon- \upsilon h=h(\upsilon-1)-(\upsilon-1)h\in M_{-3}$. 
Applying $\upsilon^{-1}$ from the right, we have $h-\upsilon h\upsilon^{-1} \in M_{-3}$.
Since $(\tau_\fg, F)$ and $(\tau_{\fg, \upsilon}, F)$ generate nilpotent orbits in the sense of \S2 (cf.\ 3.1.5), we have $h\in F^{-1}\fg_\bC$ and $\upsilon h\upsilon^{-1}\in F^{-1}\fg_\bC$, respectively (Griffiths tranversality). 
Hence $h-\upsilon h\upsilon^{-1}\in F^{-1}\fg_\bC$, and hence $h - \upsilon h \upsilon^{-1}\in M_{-3}\cap F^{-1}\fg_\bC\cap \bar F^{-1}\fg_\bC=0$. 
Here we use the fact that $(M, F\fg_\bC)$ is an $\bR$-mixed Hodge structure.
Thus we have $\tau_\fg=\tau_{\fg, \upsilon}=\rho_{\fg, \upsilon}$.
Here the second equality follows from $\tau' = \rho'$.
Hence $\sig_{\tau',1}=\sig_{\rho',\upsilon}$ and the condition (1) in 3.1.6 is verified.
\qed

\medskip

\cite{KNU10c} \S2 contains an analogous result (ibid. Theorem 2.1) and its proof by the method of \S2 in this paper. 
Theorem 6.2.1 above is its generalization by the method of \S3.

\medskip

{\bf 6.2.2.}
Recall that $\G_u= G_{\bZ,u}$. 
So $\G$ is a semi-direct product of $\G'$ in 3.1.1 (3) and $G_{\bZ, u}$. 

Take a subgroup $\Upsilon$ of $G_{\bQ,u}$.
Assume 

\medskip

(1) $G_{\bZ, u}\subset \Upsilon$.

\medskip

(2)  $\Int(\g) \Upsilon =\Upsilon$ for any $\g\in \G_G:=\Image(\G \to G_\bZ)$.

\medskip

$G_{\bZ,u}$ and $G_{\bQ,u}$ are examples of $\Upsilon$ which satisfies (1) and (2).

\medskip

\proclaim{Lemma 6.2.3}  
Under the assumption on $\Upsilon$ in $6.2.2$, $\Sig(\Upsilon)$ is strongly compatible with $\G$. 
\endproclaim 

We will deduce this from the following lemma. 

\medskip

\proclaim{Lemma 6.2.4} 
Let $\tau'$ be a face of $\sig'$ and let $\upsilon\in G_{\bQ,u}$. 
Let $\G'(\tau')$ be the face of $\Hom(P, \bN)$ corresponding to $\tau'$.  Then, the canonical map $\G(\sig_{\tau',\upsilon}) \to \G'(\tau')$ is of Kummer type, that is, it is injective and some power of any element of $\G'(\tau')$ belongs to its image. 
\endproclaim

{\it Proof.}
First, let $\g$ be an element of $\G(\sig_{\tau',\upsilon})$.
Let $\g_G$ be the image of $\g$ in $G_\bZ$.
Let $\g'$ be the image of $\g$ in $\G'(\tau')$.
Let $\g'_G$ be the image of $\g'$ in $G'_\bZ$.
$\g'_G$ is identified with its trivial extension in $G_\bZ$.
Since $\upsilon^{-1}\g_G\upsilon$ is also the trivial extension of $\g'_G$,
it coincides with $\g'_G$.
Hence $\g_G = \upsilon \g'_G\upsilon^{-1}$.
It implies that $\g_G$ is recovered from $\g'_G$, and hence $\g$ is recovered by $\g'$. 
This means that the canonical map $\G(\sig_{\tau',\upsilon}) \to \G'(\tau')$ is injective.

Next, let $\g'$ be an element of $\G'(\tau')$.
Let $\g'_G$ be the image of $\g'$ in $G'_\bQ \sub G_\bQ$.
Then, $\g'$ belongs to the image of the above map
if and only if $\upsilon \g'_G\upsilon^{-1}$ belongs to $\G_G= \Image(\G \to G_\bZ)$, which is equivalent to that $\upsilon \g'_G\upsilon^{-1}$ belongs to $G_{\bZ}$.
Hence, to prove that the canonical map is of Kummer type, 
it suffices to show that there exists $n\geq 1$ such that $\upsilon (\g'_G)^n\upsilon^{-1}$ belongs to $G_{\bZ}$ for any $\g' \in \G'(\tau')$.
But, we can take an $n$ such that $n\Ad(\upsilon)\log(\g'_G)$
sends $H_0$ into $H_0$ for any $\g'
\in \G'(\tau')$.
Multiplying $n$ if necessary, we may assume further that
$\Int(\upsilon)((\g'_G)^n)$ sends $H_0$ into $H_0$ for
any $\g' \in \G'(\tau')$.
This $n$ satisfies the desired condition. 
\qed
\medskip

{\it Proof of Lemma 6.2.3.}
First, we show that $\Sig:= \Sig(\Upsilon)$ is compatible with $\G$.
Let $\tau' \in \face(\sig')$ and $\upsilon \in \Upsilon$.
Let $\tau_1:= \sig_{\tau',\upsilon}$ as in 6.2.1.
Fix $\g \in \G$.
We have to show $\Ad(\g)(\tau_1) \in \Sig$.
Since $\G$ is a semi-direct product of $\G'$ and $G_{\bZ,u}$ (6.2.2), and since $G_{\bZ,u}$ is contained in $\Upsilon$ by the property (1) of $\Upsilon$, we may assume that $\g\in \G'$.
By the property (2) of $\Upsilon$, there exists $\upsilon' \in \Upsilon$ such that $\g_{G'}\upsilon=\upsilon'\g_{G'}$.
Then, $\Ad(\g)(\tau_1)=\Ad(\upsilon')(\sig_{\tau',\g})$, where 
$\sig_{\tau',\g}:=\{(x, \Ad(\g_{G'})x_{\fg})\;|\; x\in \tau'\}$.
Hence it is enough to show that $\sig_{\tau',\g}$ belongs to $\Sig$.
But, $\sig_{\tau',\g} = \sig_{\tau',1}$, and this certainly belongs to $\Sig$.
The compatibility of $\Sig$ and $\G$ is proved.

Next, we show that $\Sig$ and $\G$ are strongly compatible.
Let $\tau', \upsilon$, and $\tau_1$ be as in the beginning of this proof.
We have to show that any element of $\tau_1$ can be written as a finite $\bR_{\ge0}$-combination of elements $N \in\tau_1$ satisfying $\exp(N) \in \G$.
Any element of $\tau_1$ can be written as an $\bR_{\ge0}$-combination of $N_1,\ldots,
N_m \in \tau_1$ satisfying $\exp(N_{j,\fg}) \in \upsilon\G'_G\upsilon^{-1}$ for all $1 \le j \le m$, where $\G'_G=$ Image\,($\G' \to G'_{\bZ})$. 
  On the other hand, by 6.2.4, for some $n\geq 1$, the image of $\exp(nN_j)$ in $\G'(\tau')$ belongs to the image of the canonical map $\G(\tau_1) \to \G'(\tau')$ for all $j$.
By the proof of 6.2.4, this implies that $\exp(nN_j)$ belongs to $\G(\tau_1)$ for all $j$.
Since any element of $\tau_1$ can be written as an $\bR_{\ge0}$-combination of $nN_1,\ldots, nN_m$, this implies the desired statement. 
\qed
\medskip

As in 3.3.2, for a weak fan $\Sig$ which is strongly compatible with $\G$, let
$J_{\Sig}=\G \bs D_{S,\Sig}$. 

\medskip

We say a morphism $f:X\to Y$ in $\cB(\log)$ is of Kummer type if, for any $x\in X$, the homomorphism $(M_Y/\cO^\times_Y)_y\to (M_X/\cO^\times_X)_x$ of fs monoids is of Kummer type, where $y=f(x)$. 

\medskip

\proclaim{Proposition 6.2.5} 
$($Characterization of $J_{\Sig}$ for $\Sig=\Sig(\Upsilon)$ with $\Upsilon$ as in $6.2.2.)$

\medskip

{\rm(i)} $J_{\Sig}$ is of Kummer type over $S$.

\medskip

{\rm(ii)} For $S'$ of Kummer type over $S$, a morphism $S'\to J_{\Sig}$ over $S$ corresponds in one-to-one manner to an isomorphism class of LMH $H$ over $S'$ with $H(\gr^W_w)=H_{(w)}$ $(w\in \bZ)$ such that, locally on $S'$, there is a splitting $H_\bQ(\gr^W)\simeq H_\bQ$ of the filtration $W$ on the local system $H_\bQ$ on $(S')^{\log}$ satisfying the following condition{\rm: }
For any local isomorphism $H_\bZ(\gr^W)\simeq H_0(\gr^W)$ on $(S')^{\log}$ which belongs to the given $\G'$-level structure on $H_\bZ(\gr^W)$, the composition $H_\bQ\simeq H_\bQ(\gr^W)\simeq H_{0,\bQ}(\gr^W) = H_{0,\bQ}$ sends $H_\bZ$ onto $\upsilon H_0$ for some $\upsilon\in \Upsilon$.
\endproclaim

{\it Proof.}
(i) is by 6.2.4. 
  We prove (ii). 
Since $\G$ is a semi-direct product of $\G'$ and $G_{\bZ,u}$ (6.2.2), any $H \in \LMH_Q(S')$ has automatically a unique $\G$-level structure whose $\gr^W$ belongs to the given $\G'$-level structure. 
For $t \in (S')^\loga$, take a representative $\mu_t : H_{\bZ,t} @>\sim>> H_0$ of such $\G$-level structure, and put $\mu'_t:= \mu_t(\gr^W)$.

The condition that $H$ belongs to $\LMH_{Q,\G}^{(\Sig)}(S')$ (3.3.2), 
which $J_{\Sig}$ represents by 3.3.3, 
 is equivalent to that, for any $s \in S'$ and any $t \in \tau^{-1}(s) \sub (S')^\loga$, there exist a face $\rho'$ of $\sig'$ and $\upsilon \in \Upsilon$ such that the image of the homomorphism 
$$
f : \pi_1^+(s^{\log}) @>>> \sig' \times \fg_\bR
$$
in 3.3.2 is contained in $\sig_{\rho',\upsilon}$.
Here $\pi_1^+(s^{\log}):= \Hom((M_{S'}/\cO_{S'}^\times)_s, \bN)\sub \pi_1(s^{\log})$.
If this is the case, $a:= \mu_t^{-1}\upsilon^{-1}\mu'_t : H_{\bQ,t}(\gr^W) @>\sim>> H_{\bQ,t}$ extends to a local section over the inverse image by $\tau : (S')^\loga @>>> S'$ of some open neighborhood of $s$ in $S'$. 
  Since the composition $\mu'_t a^{-1} = \upsilon \mu_t$ sends $H_{\bZ,t}$ 
onto $\upsilon H_0$, this $a$ satisfies the condition in (ii). 

Conversely, if there is such a splitting $a : H_\bQ(\gr^W) @>\sim>> H_\bQ$ satisfying the condition in (ii), the homomorphism $\pi_1^+(s^{\log}) @>>> \Aut(H_{\bQ,t})$ factors into
$$
\pi_1^+(s^{\log}) 
@>>> \tsize\prod_w \Aut(H_{\bQ,t}(\gr^W_w))
@>\Int(a)>> \Aut(H_{\bQ,t}).
$$
Define $\upsilon:= \mu'_t a^{-1}{\mu}_t^{-1}$. 
 Then, $\upsilon$ belongs to $\Upsilon$ by the property of $a$, 
and the homomorphism 
$\pi_1^+(s^{\log}) @>>> \Aut(H_{0,\bQ})$ is the composite of $\pi_1^+(s^{\log}) @>>> \prod_w \Aut(H_{0,(w),\bQ}(\gr^W_w))$ and $\Int(\upsilon) : \prod_w \Aut(H_{0,(w),\bQ}(\gr^W_w)) \hookrightarrow 
\Aut(H_{0,\bQ}) @>>> \Aut(H_{0,\bQ})$. 
Thus we see that $f$ is the composite of $f' : \pi_1^+(s^{\log}) @>>> \sig' \times \fg'_\bR$ and $1_{\sig'} \times \Ad(\upsilon) : \sig' \times \fg'_\bR \hra \sig' \times \fg_\bR @>>> \sig' \times \fg_\bR$.
Let $\rho'$ be the smallest cone in $\Sig' = \face(\sig')$ which contains the image of $f'$.
($\sig' \times_{\fg'_\bR} \fg'_\bR = \sig'$.)
Then, the image of $f$ is contained in $\sig_{\rho',\upsilon} \in \Sig(\Upsilon)$.
\qed

\medskip

{\bf 6.2.6.} 
We now consider the N\'eron model and the connected N\'eron model.

Let 
$$
\Sig_0:=\Sig(G_{\bZ,u}),
$$
$$
\Sig_1:=\{\sig_{\tau',\upsilon}\in \Sig(G_{\bQ,u})\;|\;\G(\sig_{\tau',\upsilon})\to \G'(\tau')\; \text{is an isomorphism}\}.
$$

Then $\Sig_1$ is a weak fan which is strongly compatible with $\G$, and 
$J_{\Sig_1}$ coincides with the open set of $J_{\Sig(G_{\bQ,u})}$ consisting of all points at which $J_{\Sig(G_{\bQ,u})}\to S$ is strict. 

\medskip

\proclaim{Proposition 6.2.7} 
{\rm(i)} 
$\Sig_0\sub \Sig_1$. In particluar, $J_{\Sig_0}$ is 
 an open subspace of $J_{\Sig_1}$. 

\medskip

{\rm(ii)} 
$J_{\Sig_1}$ $($resp.\  $J_{\Sig_0}$$)$ has the property of $J_1$ $($resp.\  $J_0$$)$ in Theorem $6.1.1$.

\endproclaim

{\it Proof.} (i) is easy to see. (ii) follows from 3.2.8 (i) and 
6.2.5.
\qed
\medskip

This proposition completes the proof of Theorem 6.1.1.

\vskip20pt

\head
\S6.3. Case of $\cE xt^1$
\endhead

Let the notation be as in \S6.2.

We discuss more on the objects in \S6.2 in the case where, for some $w<0$, 
$H_{(k)}=0$ for $k\neq w, 0$, and $H_{(0)}=\bZ$. 

\medskip

{\bf 6.3.1.}
Let $H'_0=\gr^W_w(H_0)$. So $H_0(\gr^W)=H'_0\oplus \bZ$. 
 
Hence, 
$$
G_{\bZ,u}\simeq H'_0,\quad G_{\bQ,u}\simeq H'_{0, \bQ},
$$ 
and $\Upsilon$ satisfying (1) and (2) in 6.2.2 corresponds to a subgroup $B$ of $H'_{0, \bQ}$ containing $H'_0$ which is stable under the action of $\G'$ (3.1.7). 

The correspondence is as follows.
Let $h \in H'_{0,\bQ}$.
Then the corresponding element of $G_{\bQ, u}$ sends $e=1 \in \bZ$
to $e + h$.

\medskip

{\bf 6.3.2.} 
In the situation 6.3.1, 

\medskip

(i) $\Upsilon= G_{\bZ,u}$ corresponds to $B=H'_0$.

\medskip

(ii) $\Upsilon=G_{\bQ,u}$ corresponds to $B=H'_{0,\bQ}$.

\medskip

By the construction of the N\'eron model and the connected N\'eron model in 6.2.6--6.2.7, Proposition 6.1.4 follows from the following proposition.

\proclaim{Proposition 6.3.3} 
On the category of Kummer objects in $\cB(\log)$ over $S$, $J_{\Sig}$ for $\Sig=\Sig(\Upsilon)$ represents
$$
\Ker({\Cal E}xt^1(\bZ, H')\to R^1\tau_*H'_\bZ \to R^1\tau_*\tilde B),
$$
where we denote by $\tilde B$ the sub local system of $H'_\bQ$ corresponding to the subgroup $B \sub H'_{0,\bQ}$.
\endproclaim

\demo{Proof}
Let $S' \in \cB(\log)$ be of Kummer type over $S$. Then 
for an $H \in \LMH_Q(S')$, the image of $H$ in $R^1\tau_*\tilde B$ vanishes if and only if locally on $S'$, the extension $0\to H'_\bZ \to H_\bZ \to \bZ\to 0$ of local systems on $(S')^{\log}$ splits after pushed out by $H_\bZ'\to \tilde B$. The last condition is equivalent to the condition that locally on $S'$, there is a 
splitting $H_\bQ(\gr^W) \simeq H_\bQ$ of $W$ on $H_\bQ$ such that the local isomorphism $H_\bQ\simeq H_{0,\bQ}$ as in 6.2.5 (ii) sends $H_\bZ$ to $\bZ e + B$, i.e., onto $\upsilon H_0$ for some $\upsilon \in \Upsilon$. By 6.2.5 (ii), this is equivalent to $H\in J_{\Sig}(S')$. 
\qed
\enddemo

\proclaim{Proposition 6.3.4} 
Assume that $S$ is of log rank $\leq 1$. Assume $\G'\simeq \bZ$, and let $\g'$ be a generator of $\G'$. 
Then $\Sig_1$ in $6.2.6$ coincides with $\Sig(\Upsilon_1)$, where $\Upsilon_1$ is the subgroup of $G_{\bQ,u}$ corresponding to the subgroup
$$
B_1:=\{b\in H'_{0, \bQ} \,|\,\g' b-b \in H'_{0}\}
$$
of $H'_{0,\bQ}$. 
\endproclaim

\demo{Proof} Replacing $\g'$ by $(\g')^{-1}$ if necessary, we may assume $\g'\in \Hom(P, \bN)\sub \G'$. 
Let $N'$ be the logarithm of the image of $\g'$ in $\fg'_\bQ$. 
By abuse of notation, we denote its trivial extension in $\fg_\bQ$ also by $N'$. 
Let $\upsilon \in G_{\bQ,u}$, and let $b=\upsilon e-e\in H'_{0,\bQ}$. Then
the homomorphism $\G(\bR_{\geq 0}\Ad(\upsilon)N')\to \G'$ is an isomorphism if and only if
$\upsilon \g'_G\upsilon^{-1}\in G_\bZ$. Since $\upsilon \g'_G \upsilon^{-1}e- e = b-\g' b$, the last condition is equivalent to $b-\g'b\in H_0$, i.e., to $\upsilon \in \Upsilon_1$. 
\qed
\enddemo

\medskip

{\bf 6.3.5.}
We review some relationships with other works by summarizing 
\cite{KNU10c} 4.18. 
(In \cite{KNU10c} 4.18, \lq\lq 4.4, 4.15'' should be replaced by 
\lq\lq 4.15''.)

  Assume that we are in the situation of the beginning of this subsection. 

  If $S$ is a disk and $w=-1$, 
Green-Griffiths-Kerr 
\cite{GGK10} constructed a N\'eron model, which is 
homeomorphic to ours (6.1.1, 6.3.4, \cite{KNU10a} 8.2) 
as proved by Hayama \cite{H.p}.

  If $S$ is of higher dimensional, 
Brosnan-Pearlstein-Saito \cite{BPS.p} constructed a generalization of 
the N\'eron model of Green-Griffiths-Kerr.
In \cite{Scn.p}, Schnell constructed a connected N\'eron model.
  The relationships between theirs and the (connected) N\'eron model in this 
paper are not known.

 \medskip

 {\bf 6.3.6.}
Though we assumed $H_{(0)} = \bZ$ in this \S6.3, it is not essential. 
In fact, we can study the case where there are only two non-trivial $\gr^W_w$, say $\gr^W_a$ and $\gr^W_b$ with $a<b$ by using its special case where $b=0$ 
and $H_{(0)} = \bZ$, and by the canonical isomorphism 
$$
\cE xt^1(\gr^W_b, \gr^W_a)=\cE xt^1(\bZ, (\gr^W_b)^*\otimes \gr^W_a).
$$
\medskip

\vskip 20pt

\head
\S6.4. The proofs for the previous announcements
\endhead

  Here we give the proofs of the results in \cite{KNU10c} which were 
omitted there.

Fix $\Lambda=(H_0, W, (\langle\;,\;\rangle_w)_w, (h^{p,q})_{p, q})$ as in 2.1.1. 
In this subsection, we use the method of \S2. 
\medskip

{\bf 6.4.1.} 
  In \cite{KNU10c}, we constructed some N\'eron models by 
the method of \S2 in this paper as follows. 
  Let $\Sig'_w$ be a weak fan on $\gr^W_w$ for each $w$. 
  Let $\G'_w$ be a subgroup of $G_{\bZ}(\gr^W_w)$ for each $w$, 
which is strongly compatible with $\Sig'_w$. 
  Let $\Sig':= \prod_w \Sig'_w$ and $\G' = \prod_w \G'_w$. 
  Let $\Upsilon$ be a subgroup of $G_{\bQ,u}$. 
  Let $\Sig=\Sig_{\Upsilon}$ be the weak fan consisting of all 
$\Upsilon$-translations of the trivial extension of $\Sig'$ (\cite{KNU10c} 2.1). 
  Let $\G$ be the inverse image of $\G'$ in $G_{\bZ}$. 

  Assume 

\smallskip

(1) $G_{\bZ, u} \subset \Upsilon$, and 

\smallskip

(2) $\g\Upsilon\g^{-1} = \Upsilon$ for any $\g \in \G$. 

\smallskip

  Then, $\Sig$ is strongly compatible with $\G$ (\cite{KNU10c} 4.2). 
  This is shown in the same way as 6.2.3. (6.2.3 is shown by using 6.2.4. 
Similarly, as said in loc.\ cit., this strong compatibility is shown by using 
the fact that for any $\sig' \in \Sig'$ and any $\upsilon \in G_{\bQ,u}$, 
the map $\G(\Ad(\upsilon)(\sig)) \to \G'(\sig')$ is of Kummer type, 
where $\sig$ 
is the trivial extension of $\sig'$ (ibid.\ 4.3), which is an 
anlogue of 6.2.4, and shown in the same way as 6.2.4).

\medskip

{\bf 6.4.2.}
  Let $D':=D(\gr^W)$. 
  Let $S$ be an object of $\cB(\log)$. 
  Assume that a morphism $S \to \G' \bs D'_{\Sig'}$ is given. 
  Then, we can obtain a moduli space of LMH with the given $\gr^W$ on $S$ as the fiber product $J_{\Sig}$ 
of $S\to \G'\bs D'_{\Sig'}\leftarrow \G \bs D_{\Sig}$. 

\medskip

{\bf 6.4.3.}
  More precisely, $J_{\Sig}$ represents the functor associating with 
$T\in \cB(\log)/S$ the set of isomorphism classes of an LMH $H$ on $T$ 
endowed with a $\G$-level structure $\mu$ 
with the given $\gr^W$ satisfying the following conditions (1) and (2). 

\smallskip

(1) Locally on $T$, there is a splitting 
of the local system $a : H_{\bQ}(\gr^W) \simeq H_{\bQ}$ such that 
for any isomorphism $b : H_{\bZ} \to H_0=H_0(\gr^W)$ belonging to $\mu$, 
we have $b a \gr^W(b)^{-1} \in \Upsilon$. 

\smallskip

(2) For any $t \in T^{\log}$, any element of the smallest $\sig' \in 
\Sig'$ whose $\exp$ contains 
the image of the induced map $\pi_1^+(\tau^{-1}
\tau(t)) \to \prod_w \Aut(H_0(\gr^W_w))$ satisfies Griffiths transversality 
with respect to $\gr^W(b_t)a_t^{-1}(\bC \otimes_{\cO^{\log}_{T,\tau(t)}}
F_{\tau(t)})$ on $H_0$, where $F$ is the Hodge filtration of $H$. 

\smallskip

We do not use this fact in this paper and omit the proof. 

\medskip

{\bf 6.4.4.}
  In the following, we assume that $\Sig'=\face(\prod \sig'_w)$ for some $\sig'_w 
\in \Sig'_w$. 
  Let $\sig$ be the trivial extension of $\sig':=\prod \sig'_w$. 
  Let $\Upsilon_1$ be the subset of $G_{\bQ,u}$ consisting of all 
elements $\upsilon$ such that 
$\G(\Ad(\upsilon)(\sig)) \to \G'(\sig')$ is an isomorphism. 
  Assume further that there are only two non-trivial $\gr_k^W$.
  Then, as \cite{KNU10c} 4.6 (i) said, 
$\Upsilon_1$ forms a subgroup and satisfies (1) and (2) in 6.4.1. 
  This is seen as follows. 

Let $a, b$ be integers with $a<b$, and assume $\gr^W_w=0$ unless $w = a, b$.
 In the following, we use the same symbol for an element of
$G_{\bR}(\gr^W)$ and its trivial extension to $G_{\bR}$ by abuse
of notation.
 Let $\upsilon \in G_{\bQ,u}$.
 As seen as in the proof of 6.2.4, $\upsilon$ belongs $\Upsilon_1$ if and only if
$\upsilon \g' \upsilon^{-1}$ belongs to $G_{\bZ}$ for any $\g' \in \G'(\sig')$.
 Since $\g' \in G_{\bZ}$, the latter condition is satisfied by any
$\upsilon \in G_{\bZ,u}$, that is, $G_{\bZ,u} \subset \Upsilon_1$.

 Next, to see that $\Upsilon_1$ forms a subgroup of $G_{\bQ,u}$,
we describe the above condition as follows.
 Let $q$ be the $\Hom(\gr^W_b,\gr^W_a)$-component of $\upsilon$.
 Let $\g'_a$ and $\g'_b$ be the $G_{\bZ}(\gr^W_a)$-component and
$G_{\bZ}(\gr^W_b)$-component of $\g'$ respectively.
 Then, the $\Hom(\gr^W_b,\gr^W_a)$-component of 
$\upsilon \g' \upsilon^{-1}$ is $q\g'_b - \g'_a q$.
 Hence, the above condition is equivalent to the condition that
the element $q\g'_b - \g'_a q$ of $\Hom(\gr^W_b, \gr^W_a)$
sends $(H_0 \cap W_b)/(H_0 \cap W_{b-1})$ into
$(H_0 \cap W_a)/(H_0 \cap W_{a-1})$.
 If $q_1, q_2 \in G_{\bQ,u}$ satisfy this last condition, then,
$q_1-q_2$ also satisfies it.
 Since $G_{\bQ,u}$ is naturally isomorphic to the additive group
$\Hom(\gr^W_b, \gr^W_a)$, our $\Upsilon_1$ is a subgroup.

 Finally, let $\g \in \G$, and we prove $\g \Upsilon_1 \g^{-1} =
\Upsilon_1$.
 Let $\upsilon \in \Upsilon_1$.
 Let $\g' \in \G'(\sig')$.
 We have to prove $(\g\upsilon \g^{-1})\g'(\g\upsilon^{-1}\g^{-1}) \in G_{\bZ}$.
 Since $\G$ is a semi-direct product of $\G'$ and
$G_{\bZ,u}$ and since any element of $G_{\bZ,u}$ commutes with
$\upsilon$ (here we use the assumption that $\gr^W_w=0$ unless $w = a, b$ again),
we may assume that $\g$ is in $\G'$.
 Then, $\g^{-1}\g'\g \in \G'$ and $\log(\g^{-1}\g'\g)=\Ad(\g^{-1})\log \g'
\in \Ad(\g^{-1})(\sig')=\sig'$, where the last equality is by the fact
that $\Sig'=\text{face}\,(\sig')$ and $\G'$ are compatible.
 Hence $\g^{-1}\g'\g \in \G'(\sig')$.
 Since $\upsilon \in \Upsilon_1$, this implies that
$\upsilon(\g^{-1}\g'\g)\upsilon^{-1}\in G_{\bZ}$, and hence
$\g\upsilon \g^{-1}\g'\g\upsilon^{-1}\g^{-1} \in G_{\bZ}$, which completes the proof.

\smallskip

{\bf 6.4.5.}
  Assume that we are in the situation of 6.4.4. 
  Let $\Upsilon$ be as in 6.4.1.  If $\Upsilon \subset \Upsilon_1$, 
then, $\G \bs D_{\Sig} \to \G' \bs D'_{\Sig'}$ is strict.  
  As is said in \cite{KNU10c} 4.6 (iii), this projection is a 
relative manifold with slits in the sense of loc.\ cit. 
  In fact, by 3.4.2, this is reduced to 3.2.8 (i) and the fact that 
any strict relative log manifold is a relative manifold with slits. 

\medskip

{\bf 6.4.6.}
  We continue to assume that we are in the situation of 6.4.4.
  Assume further that there exists a strict morphism $S \to \G' \bs D'_{\Sig'}$.  
  Let $\Sig_{\Upsilon_1}$ be the weak fan corresponding to $\Upsilon_1$ (6.4.1). 
  By 6.4.5, $J_{\Sig_{\Upsilon_1}} \to S$ is strict and a relative manifold with 
slits. 

  As is said in \cite{KNU10c} 4.8, the following (i)--(iii) hold.
\smallskip

(i) For an object $T$ of $\cB/S^{\circ}$, there is a natural functorial injection $\iota(T) : \Mor_{S^\circ}(T, J_{\Sig_{\Upsilon_1}}) \to 
\{H\,|\,\text{LMH on $T$ with the given $\gr^W$ %of type $\Phi'$ 
such that $W$ on $H_\bQ$ splits locally on $T$}\}.$ 
%Here $\Phi'=(\Phi'_w)_{w\in\bZ}$ with $\Phi'_w=(\Lam'_w, \Sig'_w, \G'_w)$. 

\smallskip

(ii) If $\dim(\sig')=1$, then $\iota$ is bijective. 

\smallskip

(iii) Assume that $s\in S$ is given such that 
$\pi_1^+(s^{\log}) \to \G'(\sig')$ is bijective. 
Let $T \to S$ be a strict morphism, and let $t \in T$ be a 
point lying over $s$.  
  Then, the stalk of $\iota$ at $t$ is bijective. 

\smallskip

(i) is proved similarly as 6.2.5.

 We prove (ii) and (iii).
 Let $T \to S$ be a strict morphism from an object $T$ in $\cB(\log)$.
 Let $(H, b)$ be an LMH over $T$ with polarized graded quotients endowed
with a $\G$-level structure.
 Assume that its $\gr^W$ is the given one and that $W$ splits locally on $T$. 
 To prove (ii) (resp.\ (iii)), it is enough to show that,
under the assumption of (ii) (resp.\ (iii)), $(H,b)$ satisfies
the conditions (1) and (2) in 2.6.2 (resp. after replacing $T$ with
an open neighborhood of $t$).

 First, 
(1) is always satisfied.

 We consider the condition (2).
 Let $t \in T$ lying over $s \in S$ 
and $t' \in T^{\log}$ lying over $t$.
 Let $\tilde b_{t'}$ be a representative of the germ of $b$ at $t'$.
 Since $T \to S \to \G' \bs D'_{\Sig'}$ 
is strict, the map $f': \pi_1^+(t^{\log})
@>\simeq>> \pi_1^+(s^{\log}) \to \G'(\sig')$
is injective and its image is a face. 
 On the other hand, 
since the $W$ on the local system splits over $\bQ$, 
there exists $\upsilon \in G_{\bQ,u}$ such that the image
of the map $f: \pi_1^+(t^{\log}) \to \Aut(H_0)$ is contained in
$\exp(\Ad(\upsilon)(\sig'))$.
 We observe that if $f'$ is surjective,
then $\upsilon \in \Upsilon_1$.
 This is because $f'$ factors as
$\pi_1^+(t^{\log}) \overset f \to
\to \G(\Ad(\upsilon)(\sig')) \hookrightarrow \G'(\sig')$
so that $\G(\Ad(\upsilon)(\sig')) \hookrightarrow \G'(\sig')$ becomes an
isomorphism.

 Now we assume that $S$ is of  log rank $\leq 1$. 
  Then, the map $f'$ is either trivial or surjective.
 Hence, in both cases, the image of $f$ is contained in
$\exp(\Ad(\upsilon)(\sig'))$ for some $\upsilon \in \Upsilon_1$.
 Thus the former half of (2) is satisfied.
 Further,
$\tilde b_{t'}(\bC \otimes_{\cO_{T,{t'}}^{\log}} F_{t'})$ satisfies
the Griffiths transversality with respect to each element of the image of
$f$.
  This implies that it generates a $\tau$-nilpotent
orbit, where $\tau$ is the smallest cone of $\Sig_{\Upsilon_1}$ such that $\exp(\tau)$
contains the image of $f$ (that is, $\tau$ is $\{0\}$ or
$\Ad(\upsilon)(\sig'))$.
 Therefore, the latter half of (2) is also satisfied, which completes the
proof of (ii).

 Finally, suppose that we are in the situation of (iii).
 Then, at least, the image of $f$ at the prescribed $t$ is contained in
$\exp(\Ad(\upsilon)(\sig'))$ for some $\upsilon \in \Upsilon_1$ because
$f'$ at $t$ is bijective by the assumption.
 Since, for any point $u$ of $T$ which is sufficiently near $t$,
the map $f$ at $u$ naturally factors as
$\pi_1^+(u^{\log}) \to \pi_1^+(t^{\log}) \to \Aut(H_0)$, the image of
$f$ at $u$ is contained in
$\exp(\Ad(\upsilon)(\tau'))$, where $\tau'$ is the smallest face of $\sig'$
such that $\exp(\tau')$ contains the image of $f'$ at $u$.
 Hence, after replacing $T$, the former part of (2) is satisfied.
 In fact, $\Ad(\upsilon)(\tau')$ is the smallest cone which has the above
property that its exp contains the image of $f$ at $u$.
 Since the Griffiths transversality is satisfied with respect to each
element of this cone, the latter half of (2) is also satisfied, which
completes the proof of (iii).

\medskip

{\bf 6.4.7.}
  Let the assumption be as in 6.4.6.
  Let $\Upsilon_1$ be as in 6.4.4.  
  From now on, as in 6.3, assume further 
 $\gr^W_0=\bZ$ with the standard polarization, 
and there is $w<0$ such that $\gr^W_k=0$ for $k\neq 0, w$.

Then, $\Upsilon_1$ corresponds to the subgroup 
$B_1=\{x\in H'_{0, \bQ}\;|\; \g x-x\in H'_{0}\;\text{for all}\;\g\in \G'(\sig')\}$ 
of $H'_{0, \bQ}$. 
  This is seen similarly as in the proof of 6.3.4.

Let $\Upsilon$ be as in 6.4.1 and assume $\Upsilon \sub \Upsilon_1$.
Let $B$ be the subgroup of $H'_{0,\bQ}$ corresponding to $\Upsilon$.
  Then, on $\cB/S^{\circ}$, 
$J_{\Sig}$ represents the functor 
$\{$LMH with the given $\gr^W$ 
whose 
image in $R^1\tau_*H'_\bZ$ belongs to the kernel of $R^1\tau_*H'_\bZ\to R^1\tau_*B\}$ (\cite{KNU10c} 4.11).
  This is easily seen in the same way as 6.3.3.

\medskip

{\bf 6.4.8.}
  Let the assumption be as in 6.4.6. 
  In \cite{KNU10c} 4.4, we called
$J_{\Sig}$ the connected N\'eron model when $\Upsilon = G_{\bZ,u}$. 
  Also, in \cite{KNU10c} 4.7, we called $J_{\Sig_{\Upsilon_1}}$ the 
N\'eron model. 

  As seen easily, 
the connected N\'eron model in this sense is the connected N\'eron model 
in this paper (6.1.1). 
  However, the N\'eron model in this sense is not necessarily the 
N\'eron model in this paper (6.1.1). 
  Their relationship is as follows. 
  Consider the weak fan  
$\Sig_1 := \{\Ad(\upsilon)(\tau')\, |\, \upsilon \in G_{\bQ, u}, \tau' \in \Sig',$ 
such that $\G(\Ad(\upsilon)(\tau')) \simeq \G'(\tau')\}$. 
  Then $\Sig_{\Upsilon_1}$ is a subfan of $\Sig_1$. 
  Let $J_{\Sig_1}$ be the fiber product of $S \to \G' \bs D'_{\Sig'}\gets 
\G \bs D_{\Sig}$. 
  Then, this $J_{\Sig_1}$ is the N\'eron model in this paper, and 
$J_{\Sig_{\Upsilon_1}}$ is an open subspace of $J_{\Sig_1}$ having a 
group structure. 

\medskip

{\bf 6.4.9.}
Finally, we remark that, though the assumption of the 
existence of a strict morphism $S \to \G' \bs D'_{\Sig'}$ in 6.4.6 is rather restrictive, 
the following holds.
Let $S$ be an object of $\cB(\log)$ whose log rank is $\leq 1$.
  We assume that we are given a 
family $Q=(H_{(w)})_{w \in \bZ}$ of PLH's over $S$.  
  Then, as was stated in \cite{KNU10c} 4.14, 
locally on $S$, we can always take $\Sig'$, $\G'$, and 
a strict morphism $S\to \G'\bs D'_{\Sig'}$ such that 
the pull-back of the canonical family of PLH's on $\G' \bs D'_{\Sig'}$ 
is isomorphic to $Q$.
(The assumption in \cite{KNU10c} 4.14 was inaccurate.)
  This is shown in \cite{KU09} 4.3. 

  We remark also that, for a general $S$, 
the conclusion of 6.4.7 holds even if the strictness is 
weakened into the condition that $S \to \G' \bs D'_{\Sig'}$ is log 
injective.  
  Here a morphism $f: X \to Y$ in $\cB(\log)$ is called {\it log injective} 
if 
the induced map $f^{-1}(M_Y/\cO_Y) \to M_X/\cO_X$ is injective. 

\vskip20pt

\head
\S7. Examples and complements
\endhead

In 7.1, basic examples, including log abelian varieties and 
log complex tori, are described.
In 7.2, we discuss by examples why weak fans are essential. 
In 7.3, we propose some problems on the existence of complete weak fans. 
In 7.4, we gather some statements on the existences of period maps. 

\head
\S7.1. Basic examples
\endhead

{\bf 7.1.1.}  $0\to \bZ(1)\to H\to \bZ\to 0$.

\medskip

We take $\Lambda=(H_0, W, (\langle\;,\;\rangle_w)_w, (h^{p,q})_{p,q})$ in 2.1.1 as follows. 
Let 
$$
H_0=\bZ^2=\bZ e_1\oplus \bZ e_2,
$$
$$
W_{-3}=0 \sub W_{-2}= W_{-1}=\bR e_1 \sub W_0=H_{0,\bR},
$$
$$
\langle e_1',e_1'\rangle_{-2} =1, \quad \langle e_2',e_2'\rangle_0=1,
$$
where for $j=1$ (resp.\  $j=2$), $e_j'$ denotes the image of $e_j$ in $\gr^W_{-2}$ (resp.\   $\gr^W_0$),
$$
h^{0,0}=h^{-1,-1}=1, \quad h^{p,q}=0\;\text{for all the other}\; (p,q).
$$
Then $D=\bC$, where $z\in \bC$ corresponds to $F=F(z)\in D$ defined by
$$
F^1=0\sub  F^0=\bC \cdot (ze_1+e_2)\sub F^{-1}=H_{0,\bC}.
$$
This $D$ appeared in \cite{KNU09} and \cite{KNU.p} as Example I. 

Let 
$$
\G=G_{\bZ,u}= \pmatrix 1&\bZ \\ 0 & 1\endpmatrix\sub \Aut(H_0, W, (\langle\;,\;\rangle_w)_w),
$$
$$
\Sig=\{\sig, -\sig, \{0\}\}\quad \text{with}\;\; \pm \sig=\pmatrix 0&\pm \bR_{\geq 0}\\
0&0\endpmatrix\subset \fg_\bR.
$$
Then $\Sig$ is a fan and is strongly compatible with $\G$. 
We have 
$$
D=\bC \to \G \bs D = \bC^\times \sub \G \bs D_{\Sig} = \bP^1(\bC),
$$
where 
$\bP^1(\bC)$ is endowed with the log structure corresponding to the divisor $\{0,\infty\}$. 
Here $\bC \to \bC^\times$ corresponding to $D\to \G\bs D$ is $z\mapsto e^{2\pi iz}$,
and $0$ (resp.\  $\infty$) $\in \bP^1(\bC)$ corresponds to the class of the nilpotent orbit  $(\sig, Z)$ (resp.\  $(-\sig, Z)$), where 
$Z=D$. We have $D_{\Sig,\val}=D_{\Sig}$. 
We have 
$$
D_{\Sig}^{\sharp}=D_{\Sig,\val}^{\sharp}
= D_{\SL(2)}^I=D_{\SL(2)}^{II}= \bR\times [-\infty, \infty]
$$
as topological spaces, where the second equality 
is given by the CKS map (\S4.3), the last equality is given in \cite{KNU.p} 3.6.1 Example I, $D=\bC$ is embedded in $\bR\times [-\infty, \infty]$ by $x+iy\mapsto (x, y)$ ($x,y\in \bR$), and the projection $D_{\Sig}^{\sharp}\to \G\bs D_{\Sig}$ sends $(x,y)\in \bR\times [-\infty,\infty]$ to $e^{2\pi i(x+iy)}\in \bP^1(\bC)$. 

Next we treat this example by the formulation in \S3. 
Let $S$ be an object of $\cB(\log)$ and let $Q=(H_{(w)})_w$ with
$$
H_{(-2)}=\bZ(1), \quad H_{(0)}=\bZ, \quad H_{(w)}=0\;\;\text{for}\;w\neq  -2, 0.
$$
Then $\text{LMH}_Q={\cE}xt^1_{LMH}(\bZ, \bZ(1))$ is identified with the following functor. 
For an object $T$ of $\cB(\log)$ over $S$, $\text{LMH}_Q(T)=\G(T, M_T^{\gp})$. 
Here $a \in \G(T, M_T^{\gp})$ corresponds to the class of the following LMH $H$ on $T$. 
Define the local system $H_\bZ$ of $\bZ$-modules on $T^{\log}$ by the commutative diagram of exact sequences
$$
\matrix 0 & \to & \bZ(1) & \to &H_{\bZ}& \to & \bZ &\to & 0 \\&&
\Vert &&\downarrow && \downarrow && \\ 
0 &\to& \bZ(1)& \to& \cL_T& @>\exp>>& \tau^{-1}(M_T^{\gp}) &\to &0\endmatrix
$$
of sheaves on $T^{\log}$, where $\tau$ is the canonical map $T^{\log}\to T$, the bottom row is the exact sequence of the sheaf of logarithms $\cL_T$ (cf.\ \cite{KN99}, \cite{KU09} 2.2.4), and the right vertical arrow $\bZ\to \tau^{-1}(M_T^{\gp})$ is the homomorphism which sends $1$ to $a^{-1}$.
Let $H$ be the 
LMH $(H_\bZ, W, F)$, where
$$W_{-3}=0\sub W_{-2}= W_{-1}=\bR(1) \sub W_0=H_\bR$$
and $F$ is the decreasing filtration on $\cO_T^{\log}\otimes H_{\bZ}$ defined by
$$
F^1=0\subset 
F^0=\Ker(c:\cO_T^{\log}\otimes H_\bZ \to \cO_T^{\log})\subset 
F^{-1}=\cO_T^{\log}\otimes H_{\bZ}
 $$
with $c$ being the $\cO_T^{\log}$-homomorphism induced by $H_{\bZ}\to \cL_T\subset \cO_T^{\log}$. 
Take $\sig'=\{0\}$ in \S3.1. 
The above $\Sig$ is regarded as a fan in \S3, and the above $\G$ is regarded as $\G$ in $\S3$. $\G \bs D_{S,\Sig}$ in \S3 is identified with $S \times \G \bs D_{\Sig}\simeq S \times \bP^1(\bC)$ with $\G\bs D_{\Sig}$ as above. 
The inclusion map $\text{LMH}_Q^{(\Sig)}(T)=\Mor_{\cB(\log)}(T,  \bP^1(\bC))
\sub \text{LMH}_Q(T)=\G(T, M_T^{\gp})$ is understood by the fact that $\bP^1(\bC)$  with the above log structure represents the subfunctor $T\mapsto \G(T, M_T\cup M_T^{-1})$ of $T\mapsto \G(T, M_T^{\gp})$. 

\medskip

{\bf 7.1.2.} $0\to H^1(E)(1) \to H \to \bZ\to 0$ ($E$ a degenerating elliptic curve).

\medskip

We take $\Lambda=(H_0, W, (\langle\;,\;\rangle_w)_w, (h^{p,q})_{p,q})$ in 2.1.1 as follows. 
Let 
$$
H_0=\bZ^3= \tsize\bigoplus_{j=1}^3 \bZ e_j,
$$
$$
W_{-2}=0 \sub W_{-1}=\bR e_1 +\bR e_2\sub W_0=H_{0,\bR},
$$
$$
\langle e_2',e_1'\rangle_{-1}  =1, \quad \langle e_3',e_3'\rangle_0=1,
$$
where for $j=1,2$ (resp.\  $j=3$), $e_j'$ denotes the image of $e_j$ in $\gr^W_{-1}$ (resp.\   $\gr^W_0$),
$$
h^{0,0}=1, \quad h^{0,-1}=h^{-1,0}=1, \quad h^{p,q}=0\;\text{for all the other}\; (p,q).
$$
Then $D=\fh \times \bC$, where $\fh$ is the upper half plane, and $(\tau, z)\in \fh\times \bC$ corresponds to the following $F(\tau,z)\in D$. 
  For $\tau, z \in \bC$, we define $F=F(\tau, z) \in \check D$ by 
$$
F^1=0\sub F^0=\bC \cdot (\tau e_1+e_2)+ \bC \cdot (ze_1+e_3) \sub F^{-1}=H_{0,\bC}.
$$
This $D$ appeared in \cite{KNU09} and \cite{KNU.p} as Example II. 

Consider
$$
G_{\bZ,u}=\pmatrix 1 & 0 & \bZ\\ 0 & 1 & \bZ \\ 0 & 0 & 1
\endpmatrix \subset 
\G:= \pmatrix 1 & \bZ & \bZ\\ 0 & 1 & \bZ \\ 0 & 0 & 1
\endpmatrix\subset G_\bZ.
$$
The quotients 
$G_{\bZ,u}\bs D$ and $\G \bs D$ are \lq\lq universal elliptic curves" over 
$\frak h$ and over $\bZ \bs \frak h\simeq \Delta^*=\Delta\smallsetminus \{0\} $ ($\Delta:=\{q\in \bC\;|\;|q|<1\}$), respectively, and described as 
$$
\CD
G_{\bZ, u}\bs D = (\fh\times
\bC)/\sim\;\; @>>> \G \bs D=\bigcup_{q\in \Delta^*} \bC^\times/q^{\bZ} \\
@VVV@VVV\\
D(\gr^W_{-1}) = \fh 
@>>> \;\G' \bs D(\gr^W_{-1}) = \Delta^*.
\endCD
$$
Here, in the above diagram, the notation is as follows. 
$\bZ\bs \fh\simeq  \Delta^*$ is given by $\tau\mapsto e^{2\pi i\tau}$. 
The equivalence relation $\sim$ on $\fh\times \bC$ is defined as follows. 
For $\tau, \tau' \in \fh$ and $z, z'\in \bC$, $(\tau, z) \sim (\tau',z')$ if and only if $\tau=\tau'$ and $z\equiv z' \mod \bZ\tau+\bZ$. $D(\gr^W_{-1})$ is the $D$ for $\Lambda(\gr^W_{-1})$. 
$\G'$ denotes the image of $\G$ in $\Aut(\gr^W_{-1})$ and is identified with $\bZ$. 
The fiber $\bC/(\bZ \tau+\bZ)$ of $(\fh \times \bC)/\sim\;\to \fh$ over $\tau\in \fh$ is identified with the fiber $\bC^\times/q^{\bZ}$ of $\bigcup_q \bC^\times/q^{\bZ}\to \Delta^*$ over $q=e^{2\pi i\tau}\in \Delta^*$, via the isomorphism $\bC/(\bZ \tau+\bZ)\simeq \bC^\times/q^{\bZ},\;z\bmod (\bZ\tau+\bZ)\mapsto e^{2\pi iz}\bmod q^{\bZ}$. 
 
For $n\in \bZ$, define $N_n\in \fg_\bQ$ by
$$
N_n(e_1)=0, \;N_n(e_2)=e_1, \; N_n(e_3)=ne_1.
$$
Define the fans $\Sig$ and $\Sig_0$ by 
$$
\Sig=\{\sig_{n,n+1}, \;\sig_n\;  \; (n\in \bZ), \;\{0\}\}\quad \text{with}\;\;\sig_{n,n+1}=\bR_{\geq 0}N_n+\bR_{\geq 0} N_{n+1}, \;\;\sig_n=\bR_{\geq 0} N_n,
$$
$$
\Sig_0= \{\sig_n \;(n\in \bZ),\;\{0\}\}.
$$
Then $(\Sig, \G)$ is strongly compatible and $(\Sig_0,\G)$ is also strongly compatible. 
Let $\sig'\sub \fg_\bR(\gr^W_{-1})$ be the cone $\bR_{\geq 0}N'$, where $N'$ is the linear map $\gr^W_{-1}\to \gr^W_{-1}$ which sends $e_1'$ to $0$ and $e_2'$ to $e_1'$. 
We have a commutative diagram of fs log analytic spaces
$$
\matrix 
\G \bs D=\bigcup_{q\in \Delta^*} \bC^\times/q^{\bZ}  &\subset& 
\G\bs D_{\Sig_0}=(\G \bs D) \cup \bC^\times&\subset& 
\G\bs D_\Sig= \\
&&&&(\G \bs D) \cup (\bP^1(\bC)/(0 \sim \infty))\\
\downarrow & & \downarrow && \downarrow \\ 
&&&&\\
\G'\bs D(\gr^W_{-1})= \Delta^* & \subset & \G' \bs D(\gr^W_{-1})_{\sig'}=\Delta&=& \;\G' \bs D(\gr^W_{-1})_{\sig'}.
\endmatrix   
$$
Here the fiber of the middle (resp.\  right) vertical arrow over $0\in \Delta$ is $\bC^\times$ (resp.\  the quotient $\bP^1(\bC)/(0\sim \infty)$ of $\bP^1(\bC)$ obtained by identifying $0$ and $\infty$). 
The log structures of $\G\bs D_{\Sig_0}$ and of $\G\bs D_{\Sig}$ are given by the inverse images of $0\in \Delta$ which are divisors with normal crossings. 
The element $t\in \bC^\times$ in the fiber over $0\in \Delta$ corresponds to the class of the nilpotent orbit $(\sig_0, Z)$ with $Z=F(\bC, z)$ for $z\in\bC$ such that $t=e^{2\pi iz}$.
The element $0=\infty\in \bP^1(\bC)/(0\sim \infty)$ in the fiber over $0\in \Delta$ corresponds to the class of the nilpotent orbit $(\sig_{0,1}, Z)$ with $Z=F(\bC,\bC)$. 

The infinity of the topological spaces $D_{\Sig_0}^{\sharp}$, $D_{\Sig}^{\sharp}$, $D_{\Sig_0,\val}^{\sharp}$, $D_{\Sig,\val}^{\sharp}$, the projection $D_{\Sig}^{\sharp}\to \G\bs D_{\Sig}$, and the CKS map $D_{\Sig,\val}^{\sharp}\to D_{\SL(2)}$ are described as follows. We have a homeomorphism  $\bZ\times \bR \times  \bC\tra D_{\Sig_0}^{\sharp}\smallsetminus D, \; (n,x,z)\mapsto p_n(x,z)$, where 
$$
p_n(x,z)=\lim_{y\to \infty} F(x+iy, z+iny).
$$ 
We have a homeomorphism $\bZ\times \bR\times\bR\tra D_{\Sig}^{\sharp}\smallsetminus D_{\Sig_0}^{\sharp}, \;(n,x,a)\mapsto p_{n,n+1}(x,a)$, where
$$
p_{n,n+1}(x,a)=\lim_{y,\,y'\to \infty} F(x+i(y+y'), a+in(y+y')+iy').
$$
We have $D_{\Sig_0}^{\sharp}= D_{\Sig_0,\val}^{\sharp}$. 
The points of $D_{\Sig,\val}^{\sharp}\smallsetminus D_{\Sig_0,\val}^{\sharp}$  are written uniquely as either $p_r(x,z)$ with $r\in \bQ\smallsetminus \bZ$, $x\in \bR$ and $z\in \bC$, or $p_r(x,a)$ with $r\in \bR\smallsetminus \bQ$ and $x,a\in \bR$, or $p_{r,+}(x,a)$ with $r\in \bQ$ and $x, a\in \bR$, or $p_{r,-}(x,a)$ with $r\in \bQ$ and $x,a\in \bR$. 
Here
$$
p_r(x,z)=\lim_{y\to\infty} F(x+iy, z+iry) \quad (r\in \bR, x\in \bR, z\in \bC)
$$
(in the case $r\notin \bQ$, we have $p_r(x,z)=p_r(x, \Re(z))$), 
$$
p_{r,\pm}(x,a)=\lim_{y',\,y/y'\to \infty} F(x+iy, a+iry\pm iy').
$$ 
For $n\in \bZ$, the projection $D_{\Sig,\val}^{\sharp}\to D_{\Sig}^{\sharp}$ sends $p_r(x,z)$ ($r\in \bR$,  $n<r<n+1$, $x\in \bR$, $z\in \bC$) to $p_{n,n+1}(x, \Re(z))$, and sends $p_{r,+}(x,a)$ ($n\leq r<n+1$, $x,a\in\bR$) and $p_{r,-}(x,a)$ ($n<r\leq n+1$, $x,a\in \bR$) to $p_{n,n+1}(x,a)$.  
The projection $D^{\sharp}_{\Sig}\to \G\bs D_{\Sig}$ sends $p_n(x,z)$ ($n\in \bZ, x\in \bR, z\in \bC$) to $e^{2\pi iz}\in \bC^\times$ in the fiber over $0\in \Delta$, and sends $p_{n,n+1}(x,a)$ to $0=\infty\in \bP^1(\bC)/(0\sim \infty)$. 

As in \cite{KNU.p} 3.6.1 Example II, we have $D_{\SL(2)}^I=D_{\SL(2)}^{II}$ as topological spaces. For the admissible set $\Psi$ of weight filtrations on $H_{0,\bR}$ (4.1.2) consisting of single filtration 
$$
W'_{-3}=0 \sub W'_{-2}=W'_{-1}=\bR e_1 \sub W'_0=H_{0,\bR},
$$
the open set $D_{\SL(2)}^I(\Psi)$ (4.1.2) of $D_{\SL(2)}^I$ has a description as a topological space
$$
D_{\SL(2)}^I(\Psi)= \spl(W) \times \bar \frak h = \bR^2 \times \bar \frak h,
$$
where $\bar \frak h=\{x+iy\;|\;x\in \bR,\, 0<y\leq \infty\}$, and we identify $s\in \spl(W)$ with $(a,b)\in \bR^2$ such that $s(e'_3) = ae_1 + be_2 + e_3$. 
The embedding $D\to D_{\SL(2)}$ sends $F(\tau,z)$ ($\tau\in \fh$, $z\in \bC$) to $(\Re(z)-\Re(\tau)\Im(z)/\Im(\tau), - \Im(z)/\Im(\tau), \tau)\in \bR^2\times \fh$. 
The CKS map $D_{\Sig,\val}^{\sharp}\to D_{\SL(2)}$ has image in $D_{\SL(2)}^I(\Psi)$ and sends $p_r(x,z)$ $(r\in \bR, x\in \bR, z\in \bC)$ to $(\Re(z)-rx, -r, x+i\infty)\in \bR^2\times \bar \fh$, and sends $p_{r,+}(x,a)$ and $p_{r,-}(x,a)$ ($r\in \bQ$, $x,a\in \bR$) to $(a-rx, -r, x+i\infty)\in \bR^2\times \bar \fh$. 

Using the formulation in \S3, let $S=\Delta=\G'\bs D(\gr^W_{-1})_{\sig'}$ 
with the log structure at $0\in \Delta$, and consider $Q=(H_{(w)})_w$, where $H_{(-1)}$ is the universal log Hodge structure of weight $-1$ on $S$ (so $H_{(-1)}$ is of rank $2$), $H_{(0)}=\bZ$, and $H_{(w)}=0$ for $w\neq -1, 0$.
Then the above $\sig'$ is regarded as $\sig'$ in \S3. 
Since $\sig'\to \fg_\bR(\gr^W_{-1})$ is injective,  the above $\Sig$ and $\Sig_0$ are regarded as fans in \S3, and the above $\G$ is regarded as $\G$ in \S3. $(\Sig, \G)$ and also $(\Sig_0,\G)$ are strongly compatible. 
We have $\G\bs D_{S,\Sig}=\G\bs D_{\Sig}$, $\G \bs D_{S,\Sig_0}=\G \bs D_{\Sig_0}$. 
Furthermore, $\Sig_0$ is identified with $\Sig_j$ in \S6 for $j=0,1$ ($\Sig_0$ and $\Sig_1$ in \S6 coincide in this case). 
That is, $\G\bs D_{\Sig_0}\to S$ is the connected N\'eron model and is also the N\'eron model. 

\medskip

{\bf 7.1.3.} 
$0 \to H^1(E)(b) \to H \to \bZ\to 0$ ($b\geq 2$) ($E$ a degenerating elliptic curve).

Let $b \geq 2$. 
We take $\Lambda=(H_0, W, (\langle\;,\;\rangle_w)_w, (h^{p,q})_{p,q})$ in 2.1.1 as follows. 
Let 
$$
H_0=\bZ^3=\tsize\bigoplus_{j=1}^3 \bZ e_j,
$$
$$
W_{-2b}=0 \sub W_{1-2b}=W_{-1}=\bR e_1 +\bR e_2\sub W_0=H_{0,\bR},
$$
$$
\langle e_2',e_1'\rangle_{1-2b}  =1, \quad \langle e_3',e_3'\rangle_0=1,
$$
where for $j=1,2$ (resp.\  $j=3$), $e_j'$ denotes the image of $e_j$ in $\gr^W_{1-2b}$ (resp.\   $\gr^W_0$),
$$
h^{0,0}=1, \quad h^{1-b,-b}=h^{-b,1-b}=1, \quad h^{p,q}=0\;\text{for all the other}\; (p,q).
$$
Then $D=\fh \times \bC^2$, 
where $(\tau, z,w)\in \fh\times \bC^2$ corresponds to the following $F(\tau,z,w)\in D$.  
For $\tau,z,w \in \bC$, we define $F=F(\tau,z,w)\in \Dc$  by
$$
F^1=0\sub F^0= F^{2-b} =\bC\cdot (ze_1+we_2+e_3) \sub F^{1-b} = F^{2-b} +\bC\cdot (\tau e_1+e_2) \sub F^{-b}=H_{0,\bC}.
$$
This $D$ in the case $b=2$ appeared in \cite{KNU09} and \cite{KNU.p} as Example III. 

Like in 7.1.2, consider
$$
\G= \pmatrix 1 & \bZ & \bZ\\ 0 & 1 & \bZ \\ 0 & 0 & 1 \endpmatrix\subset G_\bZ, \quad  
N_n =
\pmatrix 0 & 1 & n\\ 0 & 0 & 0 \\ 0 & 0 & 0 \endpmatrix\in \frak g_\bQ\quad (n\in\bZ),
$$
$$
\Sig_0=\{\sig_n\; (n\in \bZ), \{0\}\}\quad \text{with}\;\sig_n=\bR_{\geq 0}N_n.
$$
Then $(\Sig_0, \G)$ is strongly compatible. 
We have $\G(\sig_0)=\exp(\bZ N_0)\simeq \bZ$. 
We have a commutative diagram in $\cB(\log)$ 
$$
\matrix 
\G(\sig_0)^{\gp} \bs D= \Delta^* \times \bC^2   & \subset&  \G(\sig_0)^{\gp} \bs D_{\sig_0}=\{(q, z, w)\in \Delta \times
\bC^2\;|\; \text{$w=0$ if $q=0$}\} \\
&&\\
\downarrow & & \downarrow  \\ 
&&\\
\G' \bs D(\gr^W_{1-2b})= \Delta^* & \subset & \;\G' \bs D(\gr^W_{1-2b})_{\sig'}= \Delta,
\endmatrix   
$$
where the isomorphism $\G(\sig_0)^{\gp}\bs D\simeq \Delta^* \times \bC^2$ sends the class of $F(\tau,z,w)$ to $(e^{2\pi i\tau}, z-\tau w, w)$. 
Here $\G'$ denotes the image of $\G$ in $\Aut(\gr^W_{1-2b})$, which is isomorphic to $\bZ$, and $\sig'$ denotes the image of $\sig_0$ in $\fg_\bR(\gr^W_{1-2b})$. 
The element $(0, z, 0)\in \Delta\times \bC^2$ corresponds to the class of the nilpotent orbit $(\sig_0, Z)$, where $Z=F(\bC,z,0)$.

Note that $\G(\sig_0)^{\gp}\bs D_{\sig_0}\to \G \bs D_{\Sig_0}$ is locally an isomorphism (Theorem 2.5.4), and is surjective. 
Hence, by the above, we have a local description of $\G \bs D_{\Sig_0}$. 
The fiber of $\G\bs D_{\Sig_0}\to \G'\bs D(\gr^W_{1-2b})_{\sig'}=\Delta$ over $q\in \Delta$ is described, as a group object, as follows. 
If $q\neq 0$, the fiber is isomorphic to $(\bC^\times)^2/\bZ$, where 
$a\in \bZ$ acts on $(\bC^\times)^2$ as $(u,v) \mapsto (uv^a, v)$. 
Here, for any fixed $\tau$ with $q=e^{2\pi i\tau}$, the isomorphism  is given by sending the class of $F(\tau, z, w)$ in the fiber to the class of 
$(e^{2\pi iz}, e^{2\pi iw})\in (\bC^\times)^2$. 
The fiber over $0\in \Delta$ is isomorphic to $\bC^\times$. 
Here the class of the nilpotent orbit $(\sig_0, F(\bC, z, 0))$ in the fiber is sent to $e^{2\pi iz}\in \bC^\times$. 

If we define $\Sig=\{\sig_{n,n+1}, \;\sig_n\;  \; (n\in \bZ), \;\{0\}\}$ with $\sig_{n,n+1}=\bR_{\geq 0}N_n+\bR_{\geq 0} N_{n+1}$ just as in 7.1.2, we have $D_\Sig=D_{\Sig_0}$ in this case.

We have a homeomorphism $\bZ\times \bR\times \bC\tra D_{\Sig_0}^{\sharp}\smallsetminus D$, $(n,x,z)\mapsto p_n(x,z)$, where 
$$
p_n(x,z)= \lim_{y\to \infty}  F(x+iy, z, -n).
$$
The projection $D_{\sig_0}^{\sharp}=D\cup \{p_0(x,z)\;|\;x\in \bR, z\in \bC\} \to \G(\sig_0)^{\gp}\bs D_{\sig_0} \sub 
\Delta \times \bC\times \bC$ sends $p_0(x,z)$ to $(0,z,0)$. 
We have $D_{\Sig_0}^{\sharp}=D_{\Sig_0,\val}^{\sharp}$. 

We have $D_{\SL(2)}^I=D_{\SL(2)}^{II}$. 
As in \cite{KNU.p} 3.6.1 Example III, if $b=2$, for the admissible set  $\Psi$ of weight filtrations on $H_{0,\bR}$ (4.1.2) consisting of single filtration 
$$
W'_{-5}=0 \sub W'_{-4}=W'_{-3}=\bR e_1 \sub W'_{-2}=H_{0,\bR},
$$
the open set $D_{\SL(2)}^I(\Psi)$ (4.1.2) of $D_{\SL(2)}^I$ has a description as a topological space
$$
D_{\SL(2)}^I(\Psi)= \spl(W) \times \bar \fh \times \bar L=\bR^2\times \bar \fh\times \bar L,
$$
where $L=\bR^2$ and $\bar L=L\cup \{\infty v\;|\;v\in L\smallsetminus \{0\}\}$ with $\infty v:=\lim_{t\to\infty} tv$. 
In this $b=2$ case, the CKS map $D_{\Sig_0,\val}^{\sharp}\to D_{\SL(2)}$ has the image in $D_{\SL(2)}^I(\Psi)$ and sends $p_n(x,z)$ to
$$
(\Re(z), \;-n\;, x+i\infty\;, \Im(z)\;,0)\in \bR^2\times \bar\fh\times L.
$$

Using the formulation in \S3, let $S=\G'\bs D(\gr^W_{1-2b})_{\sig'}
=\Delta$ with the log structure at $0\in \Delta$, and consider $Q=(H_{(w)})_w$, where $H_{(1-2b)}$ is the universal log Hodge structure of weight $1-2b$ on $S$ (so $H_{(1-2b)}$ is of rank $2$), $H_{(0)}=\bZ$, and $H_{(w)}=0$ for $w\neq 1-2b,0$.
Then the above $\sig'$ is regarded as $\sig'$ in \S3. 
Since $\sig'\to \fg_\bR(\gr^W_{1-2b})$ is injective, the above  $\Sig_0$ is regarded as a fan in \S3, and the above $\G$ is regarded as $\G$ in \S3. We have $\G \bs D_{S,\Sig_0}=\G \bs D_{\Sig_0}$. 
Furthermore, $\Sig_0$ is identified with $\Sig_j$ in \S6 for $j=0,1$ ($\Sig_0$ and $\Sig_1$ in \S6 coincide also in this case). 
That is, $\G\bs D_{\Sig_0}\to S$ is the connected N\'eron model and is also the N\'eron model. 

\medskip

{\bf 7.1.4.} 
Log abelian varieties and log complex tori.

This is a short review of \cite{KKN08} from the point of view of this paper. 

\medskip

Let $S$ be an object of $\cB(\log)$. 
A log complex torus $A$ over $S$ is a sheaf of abelian groups on the category $\cB(\log)/S$ which is written in the form 
$A=\cE xt^1_{\text{LMH}}(\bZ,H')$ for some log Hodge structure  $H'$ on $S$ of weight $-1$ whose Hodge filtration $F$ satisfies $F^1H'_{\cO}=0$ and 
$F^{-1}H'_{\cO}=H'_{\cO}$. 
Thus,  $A=\text{LMH}_Q$ for $Q=(H_{(w)})_w$, where $H_{(-1)}=H'$, $H_{(0)}=\bZ$, and $H_{(w)}=0$ for $w\neq -1,0$. 
This $H'$ is determined by $A$ (see  \cite{KKN08}). Note that we do not assume 
that $H'$ is polarizable. 
(The definition of $\text{LMH}_Q$ in \S3 did not use the polarization.) 
We call a log complex torus $A$ a log abelian variety over $S$ if the pull back of $H'$ to the log point $s$ is polarizable for any $s\in S$. For example, 
$\text{LMH}_Q$ for the $Q$ 
in the end of 7.1.2 is a log abelian variety over $S=\Delta$. 

\medskip

In  \cite{KKN08}, for a cone decomposition $\Sig$ in the sense of ibid. 5.1.2, a representable subfunctor $A^{(\Sig)}$ of a log complex torus $A$ is constructed. 
If $S$ is an fs log analytic space, $A^{(\Sig)}$ is also an fs log analytic space which is log smooth over $S$. 
We explain that this $A^{(\Sig)}$ is a special case of $\text{LMH}_Q^{(\Sig)}$ of \S3.
  In the following, for simplicity, we assume that $H'$ is polarizable and hence $A$ is a log abelian 
variety.  For the general case, we have to extend the formulation in \S3 
a little. 

Let $H'$, $Q$, $A=\cE xt^1(\bZ, H')=\text{LMH}_Q$ be as above. 
Then, locally on $S$, there are finitely generated free $\bZ$-modules $X$ and $Y$ and an exact sequence 
$$
0\to \Hom(X, \bZ)\to H'_\bZ \to Y\to 0\tag1
$$ 
of sheaves on $S^{\log}$ such that the canonical map
$$
F^0H'_{\cO} \oplus (\cO_S\otimes_{\bZ} \Hom(X, \bZ))\to  H'_\cO
$$
is an isomorphism. 
Locally on $S$, we can take $H_{0,(-1)}$ in 3.1.1 
and an exact sequence $0\to \Hom(X, \bZ)\to H_{0,(-1)}\to Y\to 0$ which corresponds to the above exact sequence (1). 
Let $\sig'$ be as in \S3.1. Then a cone decomposition $\Sig$ in the sense of \cite{KNN08} 5.1.2
 is naturally regarded as 
 a 
 rational fan $\Sig$ in $\sig' \times \Hom(X, \bR)$, and then 
 regarded as a fan $\Sig$ in \S3 by identifying $\sig'\times \Hom(X, \bR)$ with a part of $\sig'\times_{\fg_\bR'} \fg_{\bR,}$ by the following embedding. 
We send $(a, l)\in \sig'\times \Hom(X, \bR)$ to an element $(a, N)$ of $\sig'\times_{\fg'_\bR} \fg_\bR$, 
where $N$ sends $1\in \bZ=H_{0,(0)}$ to $l\in \Hom(X, \bR)\subset H_{0,(-1),\bR}$ and $N$ induces the linear map $H_{0,(-1),\bR}\to H_{0,(-1),\bR}$ given by the action of $a$. 
  Then, we have $A^{(\Sig)}=\text{LMH}_Q^{(\Sig)}$. 

\medskip

We consider the relation with classical N\'eron models of abelian varieties. Let $S$ be a smooth algebraic curve over $\bC$, let $I$ be a finite subset of $S$, and let $B$ be an abelian variety over the function field of $S$ which has good reductions at any point of $U:=S\smallsetminus I$ and semi-stable reduction at any point in $I$. Let $H'_U$ be the variation of Hodge structure on $U$ of weight $-1$ associated to $B$. Since $B$ is polarizable, $H'_U$ is polarizable. Furthermore, by the assumption of semi-stability, $H'_U$ has unipotent local monodromy at any point of $I$. Hence $H'_U$ extends to a polarizable log Hodge structure $H'$ of weight $-1$ on $S$ by the nilpotent orbit theorem of Schmid (\cite{Scm73}) (which is interpreted in terms of log Hodge structures as in \cite{KMN} Proposition 2.5, \cite{KU09} 2.5.14). Hence we have the log abelian variety $A$ over $S$ corresponding to $H'$. 
The usual (analytic) connected N\'eron model of $B$ over $S$ is the subgroup functor $G$ of $A$ in  \cite{KKN08}. 
It is clear that this $G$ coincides with  the connected N\'eron model in \S6. 
The usual  (analytic) N\'eron model of $B$ over $S$ is 
the fiber product of $A\to A/G\leftarrow \G(S, A/G)$.  
This coincides with the N\'eron model in \S6 as is seen from \S6.2 and \S6.3. 
\medskip

{\bf 7.1.5.} 
We consider the degeneration of an intermediate Jacobian. 
Let 
$$
E=\tsize\bigcup_{q\in\Delta^*} \;\bC^\times/q^{\bZ}\to \Delta^*
$$ 
be the family of elliptic curves which appeared in 7.1.2, and consider the family 
$$
X:=E^3\to S^*:=(\Delta^*)^3.
$$ 
Let $S=\Delta^3$ and endow $S$ with the log structure associated to the divisor $S\smallsetminus S^*$.  
Consider the connected N\'eron model $J_0$ over $S$ (\S6.1) of the second intermediate Jacobian $J^2(X/S^*)$ of $X$ over $S^*$.  It is the connected N\'eron model for the case
$Q=(H_{(w)})_w$ with $H_{(0)}=\bZ$, $H_{(-1)}$ is the unique extension on $S$ as a log Hodge structure of weight $-1$ (obtained by the nilpotent orbit theorem of Schmid \cite{Scm73}) of the polarizable variation of Hodge structure $H^3(X/S^*)(2)$ on $S^*$, and $H_{(w)}=0$ for $w\neq -1, 0$.  Here $H^3(X/S^*)$ is the third higher direct image of the constant Hodge structure $\bZ$ on $X$ under $X\to S^*$. 

As a group object, the fiber of $J_0\to S$ over $(q_1,q_2, q_3)\in \Delta^3$ has the following structure. 
If $q_1q_2q_3\neq 0$, the fiber (which is the second intermediate Jacobian of a product of three elliptic curves) is a quotient of $\bC^{10}$ by a discrete group which is isomorphic to $\bZ^{20}$. 
If one of $q_j$ is zero and two of them are non-zero, the fiber is isomorphic to a quotient of $\bC^9$ by a discrete subgroup which is isomorphic to $\bZ^{14}$. 
If two of $q_j$ are zero and one of them is non-zero, the fiber is isomorphic to a quotient of $\bC^8$ by a discrete subgroup which is isomorphic to $\bZ^{10}$. 
If $q_1=q_2=q_3=0$, the fiber is isomorphic to a quotient of $\bC^7$ by a discrete subgroup which is isomorphic to $\bZ^7$.

This can be seen as follows. 
Since the connected N\'eron model is the slit Zucker model, the fiber has the form $V/L$, 
where $V$ is the small Griffiths part of the fiber of the vector bundle $H_{(-1),\cO}/F^0H_{(-1),\cO}$ and $L$ is the stalk of $\tau_*(H_{(-1),\bZ})$ ($\tau : S ^\loga \to S$). 
Since the fiber is Hausdorff (Theorem 3.2.8 (iii)), $L$ must be discrete in $V$. 
Let $H'$ be the unique extension on $\Delta$ of $H^1(E/\Delta^*)(1)$ on $\Delta^*$, as a log Hodge structure of weight $-1$. 
Then we have
$$
H_{(-1)}\simeq (\tsize\bigotimes_{j=1}^3  \text{pr}_j^*(H'))(-1) \oplus  \text{pr}_1^*(H')^{\oplus2}\oplus  \text{pr}_2^*(H')^{\oplus2}\oplus  \text{pr}_3^*(H')^{\oplus 2}.
$$ 
From this, we have that if $q_1=0$ and $q_2q_3\neq 0$ (resp. $q_1=q_2=0$ and $q_3\neq 0$, resp. $q_1=q_2=q_3=0$), 
the dimension of $V$ is $3+2+2+2=9$ (resp.\  $2+2+2+2=8$, resp. $1+2+2+2=7$), and the rank of $L$ is $4+2+4+4=14$ (resp.\  $2+2+2+4=10$, resp. $1+2+2+2=7$). 

The N\'eron model coincides with the connected N\'eron model in this case. 

\bigskip

\head
\S7.2. The necessity of weak fan
\endhead 

We explain by examples that a weak fan (not a fan) is necessary to construct the connected N\'eron model. 

\medskip

{\bf 7.2.1.}  
We give an example in which the weak fan $\Sig_0$ in \S6, which is used for the construction of the connected N\'eron model, is not a fan. 

In \S6, assume $H_{0,(-2)}=L\otimes L$ with $L=\bZ^2=\bZ e_1+\bZ e_2$,  $H_{0,(0)}=\bZ$, $H_{0,(w)}=0$ for $w\neq  -2, 0$, 
$\sig'=\bR_{\geq 0}N'_1+\bR_{\geq 0}N'_2$, where $N'_j : H_{0,(-2),\bQ}\to H_{0,(-2),\bQ}$ ($j=1,2$) are as follows:  $N'_1=N'\otimes 1$, $N'_2=1\otimes N'$ with $N':L\to L$ being the homomorphism defined by $N'(e_1)=0$, $N'(e_2)=e_1$.
Let $e$ be the base $1$ of $H_{0,(0)}$, let $\tau=\bR_{\geq 0}N_1+\bR_{\geq 0}N_2$ with $N_j\in \fg_\bQ$ being the extension of $N'_j$ such that $N_j(e)=0$. 
Then the weak fan  $\Sig_0=\{\Ad(\g)(\alpha)\;| \; \g\in G_{\bZ,u}, \,\alpha\;\text{is a face of}\;\tau\}$ in \S6 is not a fan.
To see this, let $\g_{m,n}\in G_{\bZ,u}$ for integers $m,n>0$ be the element which fixes all elements of $H_{0,(-2)}$ and which sends $e$ to $e+me_1\otimes e_2-ne_2\otimes e_1$. 
Then $\tau\cap \Ad(\g_{m,n})(\tau) =\bR_{\geq 0}\cdot (mN_1+nN_2)$ (note that $\Ad(\g_{m,n})(mN_1+nN_2)=mN_1+nN_2$) and this intersection is not a face of $\tau$.

For a later use in 7.3.3, we remark that if $\tau_1$ is a subcone of $\tau$ of rank $2$, there is no fan $\Sig$ which is stable under the adjoint action of $G_{\bZ,u}$ such that $\tau_1\sub \sig$ for some $\sig\in \Sig$. 
In fact, the interior of $\tau_1$ contains $mN_1+nN_2$ and $m'N_1+n'N_2$ for some integers $m,n,m',n'>0$ such that $mn'-m'n\neq 0$. If such a fan $\Sig$ exists, the smallest $\sig\in \Sig$ such that $\tau_1\sub \sig$ exists, and any interior point of $\tau_1$ is an interior point of $\sig$.  Hence $\sig$ and $\Ad(\g_{m,n})(\sig)$ contains a common interior point $mN_1+nN_2$ and hence we should have $\Ad(\g_{m,n})(\sig)=\sig$ and hence  $\Ad(\g_{m,n})^k(\sig)=\sig$ for all $k\in \bZ$. Let $N_0\in \fg_\bQ$ be the element which kills $H_{0,(-2)}$ and which sends $e$ to $e_1\otimes e_1\in H_{0,(-2)}$. 
Then $\Ad(\g_{m,n})^k(m'N_1+n'N_2) =m'N_1+n'N_2+k (m'n-mn')N_0$. But since $\sig$ is a finitely generated sharp cone, it can not contain the set $\{m'N_1+n'N_2+k (m'n-mn')N_0\;|\;k\in \bZ\}$. 

This $((H_{0,(w)})_w, \sig')$ appears in \S3 in the case $S$ and $(H_{(w)})_w$ are as follows: $S=\Delta^2$ with the log structure associated to the divisor $S\smallsetminus S^*$ with normal crossings, where $S^*=(\Delta^*)^2$. $H_{(0)}=\bZ$. $H_{(-2)}$ is the extension of $\text{pr}_1^*H^1(E/\Delta^*)(1) \otimes \text{pr}_2^*H^1(E/\Delta^*)(1)$ on $S^*$ to $S$ as a log Hodge structure of weight $-2$, where $E= \bigcup_{q\in \Delta^*} \bC^\times/q^{\bZ}$ as in 7.1.5 and $H^1(E/\Delta^*)$ is the VHS on $\Delta^*$ obtained as the first higher direct image of the constant Hodge structure $\bZ$ on $E$ under $E \to \Delta^*$. 
$H_{(w)}=0$ for $w\neq -2, 0$. 
\medskip

{\bf 7.2.2.} 
We show that in the example 7.1.5, we need a weak fan, not a fan, to construct the connected N\'eron model of the second intermediate Jacobian.  Let $X=E^3\to S^*=(\Delta^*)^3\sub S=\Delta^3$ be as in 7.1.5. Let $L=\bZ^2=\bZ e_1+\bZ e_2$ and
 $N':L\to L$ be as in 7.2.1.
We take
$H_{0,(0)}=\bZ$, $H_{0,(-1)}=L\otimes L\otimes L\oplus L^{\oplus 2}\oplus L^{\oplus 2}\oplus L^{\oplus 2}$, $H_{0,(w)}=0$ for $w\neq -1,0$, and regard $L$ as a stalk of $H^1(E/S^*)(1)_{\bZ}$. 
Define $N'_j\in \fg'_\bQ$ ($j=1,2,3$) as $N'_1=(N'\otimes 1\otimes 1)\oplus (N')^{\oplus 2}\oplus 0\oplus 0$, $N'_2=(1\otimes N' \otimes 1)\oplus 0 \oplus (N')^{\oplus 2} \oplus 0$, $N'_3=(1\otimes 1\otimes N') \oplus 0\oplus 0\oplus (N')^{\oplus 2}$. Extend $N'_j$ to $N_j\in \fg_\bQ$ by $N_j(e)=0$, where $e$ is $1\in \bZ=H_{0,(0)}$.   
Take $P=\bN^3$,  take the inclusion morphism $S\to \toric_{\sig'}=\Spec(\bC[P])_{\an}=\bC^3$, and take the homomorphism $\G'=\Hom(P^{\gp}, \bZ)=\bZ^3 \to \fg'_\bR$ given by $(N_1', N_2', N'_3)$. 
We show that the weak fan $\Sig_0$ is not a fan. 
We show more strongly that if $\tau_1$ is a finitely generated subcone of $\tau:=\bR_{\geq 0} N_1+\bR_{\geq 0}N_2+\bR_{\geq 0}N_3$ of rank $3$, there is no fan $\Sig$ which is stable under the adjoint action of $G_{\bZ,u}$ such that $\tau_1\sub \sig$ for some $\sig\in \Sig$. 
The interior of  $\tau_1$ has elements of the form $mN_1+nN_2+\ell N_3$ and $m'N_1+n'N_2+\ell N_3$ with $m,n,m',n',\ell\in \bZ$ such that $mn'-m'n\neq 0$.
Let $\g$ be the element of $G_{\bZ,u}$ which fixes every element of $H_{0,(-1)}$ and which sends $e\in H_{0,(0)}$ to $e+(me_1\otimes e_2\otimes e_1-ne_2\otimes e_1\otimes e_1,\,0,0,0)$.
Then $\Ad(\g)(mN_1+nN_2+\ell N_3)=mN_1+nN_2+\ell N_3$, $\Ad(\g)^k(m'N_1+n'N_2+\ell N_3)= m'N_1+n'N_2+\ell N_3+ k(m'n-mn')N_0$, where $N_0(H_{0,(-1)})=0$ and 
$N_0(e)=(e_1\otimes e_1\otimes e_1,\,0,0,0)$. 
Hence, by the same argument as in 7.2.1, we can deduce that there is no such fan $\Sig$. 
 
\medskip

{\bf 7.2.3.} 
In \cite{KNU10c} 4.13, in a certain case, we gave  a necessary condition for $\Sig_0$ in \S6 to be a fan.  

\bigskip

\head
\S7.3. Completeness of weak fan
\endhead 

We discuss the completeness of a weak fan $\Sig$.
We explain that, in the pure case, the conjecture in \cite{KU09} \lq\lq Added in the proof" after 12.7.7 is false. 
Here we reformulate such a conjecture by replacing  fan by weak fan (generalizing it to the mixed Hodge case).

\medskip

{\bf 7.3.1.} 
In the situation of \S2 (resp.\ \S3), we say a weak fan $\Sig$ is {\it complete } (resp.\ {\it relatively complete}) if the following condition (1) is satisfied.

\smallskip

(1) For any rational nilpotent cone $\tau$ such that a $\tau$-nilpotent orbit exists, there is a finite rational subdivision $\{\tau_j\}_{1\leq j\leq n}$ of $\tau$ having the following property. For each $1\leq j\leq n$, there is an element $\sig$ of $\Sig$ such that $\tau_j\sub \sig$. 

\medskip

{\bf 7.3.2.} 
Assume that we are in the situation of \S2. 
In the pure case, for a fan, this completeness is weaker than the completeness defined in \cite{KU09} 12.6.1.

In \cite{KU09} 12.6.8, the following relation of completeness and  compactness was proved for the completeness of a fan in the sense of \cite{KU09} 12.6.1.
The same proof works for the present completeness of a weak fan.
Note that if $\G$ is a subgroup of finite index of $G_\bZ$ and $\Sig$ is a weak fan which  is compatible with $\G$, then $\Sig$ is strongly compatible with $\G$. 

\proclaim {Theorem}
Assume that we are in the pure case and that we are in the classical situation \cite{KU09} $0.4.14$.
Let $\Sig$ be a weak fan in $\fg_\bQ$, let $\G$ be a subgroup of $G_\bZ$
of finite index, and assume that $\G$ is compatible with $\Sig$.
Then, $\G\bs D_\Sig$ is compact if and only if $\Sig$ is complete.
\endproclaim

In the pure case and in the classical situation \cite{KU09} 0.4.14, a complete fan was constructed in  \cite{AMRT75}.  

In general, we make a 

\proclaim
{Conjecture} 
{\rm There is a complete weak fan which is compatible with $G_{\bZ}$ (and hence with any subgroup of $G_{\bZ}$).} 
\endproclaim

\medskip

{\bf 7.3.3.} 
Assume that we are in the situation of \S2. 
We show that even in the pure case, the  conjecture in 7.3.2 becomes not true if we replace weak fan by fan.

Consider the case $\Lambda=(H_0, W, (\langle\;,\;\rangle_w)_w, (h^{p,q})_{p,q})$, where 
$$
H_0=(L\otimes L) \oplus \text{Sym}^2(L)\quad \text{with}\;\;L=\bZ^2=\bZ e_1\oplus \bZ e_2,
$$
$$
W_1=0\sub W_2=H_{0,\bR},
$$
$\langle\;,\;\rangle_2$ is induced from the anti-symmetric $\bZ$-bilinear form
$L\times L\to \bZ$ which sends $(e_2,e_1)$ to $1$, and 
$$
h^{2,0}=h^{0,2}=2, \quad h^{1,1}=3, \quad h^{p,q}=0\;\;\text{for all the other}\;(p,q).
$$

Let $N:L\to L$ be the homomorphism defined by $N(e_1)=0$, $N(e_2)=e_1$. 
Define $N_1,N_2, N_3\in \fg_\bQ$ by
$$
N_1= (N \otimes 1) \oplus 0, \quad N_2=(1\otimes N) \oplus 0,\quad N_3=0\oplus \text{Sym}^2(N),
$$ 
and let $\tau=\bR_{\geq 0}N_1+\bR_{\geq 0}N_2+\bR_{\geq 0}N_3$. 

Note that there exists a $\tau$-nilpotent orbit. In fact, for any decreasing filtration $F$ on $L_\bC$ such that $F^0=L_\bC$, $\dim_{\bC}(F^1)=1$, and $F^2=0$, 
$(N_1, N_2, N_3, (F\otimes F)\oplus \text{Sym}^2(F))$ generates a nilpotent orbit. 
(This nilpotent orbit appears in the degeneration of $H^1(E_1)\otimes H^1(E_2) \oplus \text{Sym}^2(H^1(E_3))$, where $E_j$ are degenerating elliptic curves.)  
We show that if $\tau_1$ is a finitely generated subcone of $\tau$ of rank $3$, there is no fan $\Sig$ which is compatible with $G_{\bZ}$ such that  $\tau_1\sub \sig$ for some $\sig\in \Sig$. 

Assume that such a $\tau_1$ and such a fan $\Sig$ exist. 
Then the interior of the cone $\tau_1$ contains elements of the form $mN_1+nN_2+\ell N_3$ and $m'N_1+n'N_2+\ell N_3$ with $m,n,m',n',\ell\in \bZ$ such that $mn'-m'n\neq 0$. 
Let $\g$ be the element of $G_\bZ$ which does not change $(e_1\otimes e_1, 0)$, $(e_2\otimes e_2, 0)$, $(0, e_1e_2)$, $(0, e_1^2)$ and for which  
$$
\align
&\g (0, e_2^2)=(0, e_2^2) + (m e_1\otimes e_2-n e_2\otimes e_1,0),\\
&\g(e_1\otimes e_2,0)= (e_1\otimes e_2,0) + (0, ne_1^2),\\
&\g(e_2\otimes e_1, 0)=(e_2\otimes e_1,0)+(0, me_1^2).
\endalign
$$
Then $\Ad(\g)(mN_1+nN_2+\ell N_3)=mN_1+nN_2+\ell N_3$, 
$\Ad(\g)^k(m' N_1+n'N_2+\ell N_3)= m'N_1+n'N_2+\ell N_3+ k(m'n-mn')N_0$, where $N_0(0, e_2^2)=(e_1\otimes e_1, 0)$ and $N_0$ kills all $(e_j\otimes e_k,0)$ and kills $(0,e_1^2)$ and $(0, e_1e_2)$. 
By the argument in 7.2.1, we see that  $\Sig$ as above does not exist. 

\medskip 

{\bf 7.3.4.} 
We explain variants of the conjecture in 7.3.2. 
Let $I$ be the set of all rational nilpotent cones $\tau$ for which a $\tau$-nilpotent orbit exists. 

\smallskip

Conjecture a.
{\rm There is a weak fan $\Sig$ which is compatible with $G_{\bZ}$ such that $\bigcup_{\tau\in I} \; \tau=  \bigcup_{\sig\in \Sig} \; \sig$.}

\smallskip

Conjecture b.
{\rm There is a weak fan $\Sig$ which is compatible with $G_{\bZ}$ such that $\bigcup_{\tau\in I} \; \tau\sub  \bigcup_{\sig\in \Sig} \; \sig$.}

\smallskip

Clearly we have the implications
$$
\text{(Conjecture in 7.3.2)} \Rightarrow \text{(Conjecture b)} \Leftarrow \text{(Conjecture a)}.
$$

If weak fan in Conjecture a is replaced by fan, it becomes the conjecture stated after \cite{KU09} 12.7.7 (in the pure case). 
We show that the last conjecture fails. 
We show more strongly that Conjecture b becomes not true if weak fan in Conjecture b is replaced by fan (even in the pure case).

Consider the example in 7.3.3. 
Let $a_1, a_2, a_3$ be positive 
real numbers which are linearly independent over $\bQ$.  
If a fan 
$\Sig$ as in Conjecture b exists, there is $\sig\in \Sig$ such that $a_1N_1+a_2N_2+a_3N_3\in \sig$. 
Let $\tau_1=\tau\cap \sig$, where  $\tau=\bR_{\geq 0}N_1+\bR_{\geq 0}N_2+\bR_{\geq 0}N_3$. 
Since $\tau_1$ is rational and contains $a_1N_1+a_2N_2+a_3N_3$ with  $a_1,a_2,a_3$  linearly independent over $\bQ$, $\tau_1$ is of rank $3$. 
But this contradicts what we have  seen in 7.3.3.

 \medskip
 
 {\it Question.} 
In case where the conjecture in 7.3.2 is valid, we may ask if we can further 
require that, for any cone $\sig$ in the weak fan in the conjecture, 
a $\sig$-nilpotent orbit exists.  
Note that
a delicate example constructed by Watanabe (\cite{W08}) is harmless 
to this question.

 \medskip

{\bf 7.3.5.} 
Next, assume that 
we are in the situation of \S3.

Note that if $\G_u$ is of finite index in $G_{\bZ,u}$ and $\Sig$ is compatible with $\G$, then  $\Sig$ is strongly compatible with $\G$ (here $\G_u$ and $\G$ 
are as in 3.1.7). 

If $\G_u$ is of finite index in $G_{\bZ,u}$ and $\Sig$ is compatible with $\G$, then the relative completeness of $\Sig$ is related to the properness of $\G \bs D_{S, \Sig}\to S$ in certain cases (\cite{KNU10a} \S5 Proposition). 

\proclaim{Conjecture} 
{\rm There is a relatively complete weak fan which is  compatible with any $\G$ as in 3.1.7.}
\endproclaim

This conjecture becomes not true if we replace weak fan by fan. 
In the cases 7.2.1 and 7.2.2, there is no relatively complete fan which is stable under the adjoint action of $G_{\bZ,u}$,  as is seen from the arguments in 7.2.1 and 7.2.2. 
\medskip

{\bf 7.3.6.} 
We show that  if  the following conditions (1) and (2) are satisfied, then the conjecture in 7.3.5 is true
 even if we replace weak fan by fan. 
 
\smallskip

(1) There are only two $w$ such that $H_{(w)}\neq 0$. 

\smallskip

(2) $\G'\simeq \bZ$ (cf.\ 3.1.1 (3)). 

\smallskip

The following construction is a variant of that in \cite{KNU10a} \S3.

Assume (1) and (2), and assume $H_{(a)}\neq 0$, $H_{(b)}\neq 0$, $a<b$. 
Note that the fs monoids $P$ and $\Hom(P, \bN)$ are isomorphic to $\bN$. Let $\g'$ be a generator of $\Hom(P, \bN)\subset \G'$, and let $N'_a: H_{0,(a),\bQ}\to H_{0,(a),\bQ}$ and $N'_b: H_{0, (b),\bQ}\to H_{0,(b),\bQ}$ be the actions of $\log(\g')$ on $H_{0,(a),\bQ}$ and on $H_{0,(b),\bQ}$, respectively. 
Let $V= \Hom_{\bQ}(H_{0,(b),\bQ}, H_{0,(a),\bQ})$, and define $\bQ$-subspaces $X$ and $Y$ of $V$ as follows:
$$
\align
&X:=\{h\in V\;|\; hM(N'_b, W(\gr^W_b))_w\subset M(N'_a, W(\gr^W_a))_{w-2}\;\text{for all}\;w\} + \{hN'_b-N'_ah\;|\;h\in V\},\\
&Y:=\{h\in X\;|\; hN'_b=N'_ah\}\sub X.
\endalign
$$
  Let $\g_a:=\gr^W_a(\g')$ and $\g_b:=\gr^W_b(\g')$. 
Fix a finitely generated
 $\bZ$-submodule $L$ of $V$ which is $\g_a$-stable and $\g_b$-stable such that 
 $$L\supset \Hom_{\bZ}(H_{0,(b)}, H_{0, (a)}) + \{hN'_b-N'_ah\;|\; h\in \Hom_{\bZ}(H_{0,(b)}, H_{0,(a)})\}.$$ 
Fix a $\bZ$-homomorphism $s: (X\cap L)/(Y\cap L)\to X\cap L$ which is a section of the projection $X\cap L\to(X\cap L)/(Y\cap L)$. 
For an element $x$ of $X/Y$, let $d(x)$ be the order of the image of $x$ in $X/((X\cap L)+Y)$ which is a torsion element. 
Fix a $\bZ$-basis $(e_j)_{1\leq j\leq m}$ of $Y\cap L$. 
Let $v$ be the image of $\g'$ under the canonical injection $\Hom(P, \bN)\to \sig'$, 
which is a generator of $\sig'\simeq \bR_{\geq 0}$.

For $x\in X/Y$ and $n=(n_j)_{1\leq j\leq m}\in \bZ^m$, let $\sig(x,n)$ be nilpotent cone in $\sig' \times_{\fg'_\bR}\fg_\bR$ 
generated by all elements of the form $(v,N)$, where $N$ is an element of 
$\fg_\bR$ satisfying the following (3)--(5). 
\smallskip

(3) The $\Hom_{\bQ}(H_{0,(a),\bQ}, H_{0,(a), \bQ})$-component of $N$ coincides with $N'_a$. 
\smallskip

(4) The $\Hom_{\bQ}(H_{0,(b),\bQ}, H_{0,(b), \bQ})$-component of $N$ coincides with $N'_b$. 

\smallskip

(5) The $\Hom_{\bQ}(H_{0,(b),\bQ}, H_{0,(a), \bQ})$-component of $N$ has the form
$$
s(x)+\frac{1}{d(x)}\tsize\sum_{j=1}^m c_je_j\quad \text{with}\;\; n_j\leq c_j\leq n_j+1\;\;\text{for all}\;\;j.
$$
\smallskip

Let 
$$
\Sig=\{\text{face of}\; \sig(x,n)\;|\;x\in X/Y, n\in \bZ^m\}.
$$ 
Then we have:

\medskip

\noindent
$\Sig$ is a relatively complete fan which is compatible with any $\G$ as in 3.1.7. 

\medskip

First we prove the compatibility. 
  For simplicity, we assume that either $N'_a$ or $N'_b$ is non-trivial and 
regard $\sig' \times_{\fg'_\bR}\fg_\bR \sub \fg_\bR$ via the second projection. 
  The case $N'_a=N'_b=0$ is similar and easier. 

  Let $\g$ be any element of $\G$. 
It is sufficient to prove $\Ad(\g)\sig(x, n)\in \Sig$ for any $x\in X/Y$ and $n\in \bZ^m$. 
 Write the $\Hom_{\bQ}(H_{0,(b),\bQ}, H_{0,(a), \bQ})$-component of
$\g^{-1}$ by $h$. 

  Let
$N \in \sig(x, n)$ be an element satisfying (3)--(5). 
Then, since $\g_b$ and $N'_b$ commute, 
the $\Hom_{\bQ}(H_{0,(b),\bQ}, H_{0,(a), \bQ})$-component of 
$\g N\g^{-1}$ is 
$$
\g_a s(x)\g_b^{-1}+ \g_a(N'_ah-hN'_b)+ d(x)^{-1}\tsize\sum_{j=1}^m c_je_j. \tag6
$$ 
  Write 
$\g_a s(x)\g_b^{-1}+ \g_a(N'_ah-hN'_b) = s(y)+ z$ with $y\in X/Y$ and $z\in Y$.   Since $L$ is $\g_a$-stable and $\g_b$-stable, we have $d(y)=d(x)$ and 
 $d(x)z\in Y\cap L$. 
Write $d(x)z= \sum_{j=1}^m m_je_j$ with $m_j\in \bZ$, and let $m=(m_j)_j\in \bZ^m$. 
Then, $(6)$ equals to 
$$
s(y)+ d(y)^{-1}\tsize\sum_{j=1}^m (c_j + m_j)e_j.
$$
  This means 
$$
\Ad(\g)\sig(x, n) = \sig(y, n+m).
$$ 

Next, the relative completeness can be proved by noticing the following two points (7) and (8). Let $R$ be the set of all elements $N$ of $\fg_\bR$ 
such that the $\Hom_\bR(H_{0,(a),\bR}, H_{0,(a),\bR})$-component of $N$ coincides with $N'_a$ and the $\Hom_\bR(H_{0,(b),\bR},H_{0,(b),\bR})$-component of $N$ coincides with $N'_b$. For $N\in R$, let $N_u$ be the  $\Hom_\bR(H_{0,(b),\bR}, H_{0,(a),\bR})$-component of $N$. Then:
 
 \medskip
 
 (7) Let $N\in R$. Then $M(N, W)$ exists if and only if $N_u$ belongs to the $\bR$-subspace $\bR\otimes_{\bQ} X$ of $\Hom_\bR(H_{0,(b),\bR}, H_{0,(a),\bR})$.

\medskip

(8) Let $N_1, N_2\in R$ and assume $(N_1)_u, (N_2)_u\in \bR\otimes_{\bQ} X$ in
$\Hom_\bR(H_{0,(b),\bR}, H_{0,(a),\bR})$. 
Then, 
 $N_1N_2=N_2N_1$ if and only if $(N_1)_u-(N_2)_u$ belongs to the $\bR$-subspace  $\bR\otimes_{\bQ} Y$ of 
 $\Hom_\bR(H_{0,(b),\bR}, H_{0,(a),\bR})$.

 \medskip
 
{\bf 7.3.7.} {\it Remark} 1. 7.3.6 gives a proof of \cite{KNU10a} \S3, and hence a proof of the theorem in \cite{KNU10a} \S2. (In fact, for the proofs, we have to show that the fan $\Sig$ defined above satisfies the admissibility in the sense of \cite{KNU10a}, i.e., the condition (2) in 1.2.6. But this can be checked easily.)

\medskip
{\it Remark} 2.   Here is a correction to \cite{KNU10a}.
  The submodule $L$ in ibid.\ \S3 should have been imposed to be 
$\g'$-stable as in 7.3.6 above.  
  Without this condition, the resulting $\Sig$ was not necessarily compatible with $\G$. 

\medskip

{\bf 7.3.8.} 
In the part of 7.1.1 concerning the situation of \S3, the fan $\Sig$ is relatively complete. In the part of 7.1.2 concerning  the situation of \S3, the fan $\Sig$ is relatively complete. In the part of 7.1.3 concerning the situation of \S3, the fan $\Xi$ (3.1.6) 
is relatively complete. In 7.1.5, there is no relatively complete fan as is seen by the argument in 7.2.2.

\bigskip

\head 
\S7.4. Extended period maps
\endhead

\medskip

In this subsection, we see extensions of period maps.

\proclaim{Theorem 7.4.1}
Let $S$ be a connected, log smooth, fs log analytic space, and let $U$ be the open subspace of $S$ consisting of all points of $S$ at which the log structure of $S$ is trivial.
Let $H$ be a variation of mixed Hodge structure on $U$ with polarized graded quotients for the weight filtration, and with unipotent local monodromy along $S\smallsetminus U$. 
Assume that $H$ extends to a log mixed Hodge structure on $S$ $($that is, $H$ is admissible along $S\smallsetminus U$ as a variation of mixed Hodge structure$)$. 
Fix a base point $u\in U$ and let $\Lambda=(H_0,W, (\lan\;,\;\ran_w)_w, (h^{p,q})_{p,q})$ be $(H_{\bZ,u}, W, (\lan\;,\;\ran_{w,u})_w, (\text{the Hodge numbers of $H$}))$.
Let $\G$ be a subgroup of $G_\bZ$ which contains the global monodromy group $\Image(\pi_1(U,u)\to G_\bZ)$ and assume that $\G$ is neat. 
Let $\vf: U\to \G\bs D$ be the associated period map. 
Then{\rm:}

\medskip

{\rm (i)} Let $S_{\val}^{\log}$ be the topological space over $S^{\log}$ defined as in \cite{KU09} $3.6.26$, which contains $U$ as a dense open subspace. Then 
the map $\vf:U\to\G\bs D$ extends uniquely to a continuous map
$$
S_\val^\loga\to\G\bs D_{\SL(2)}^I.
$$

{\rm(ii)} For any point $s\in S$, there exist an open neighborhood
$V$ of $s$, a log modification $V'$ of $V$ $(\text{\cite{KU09}}\ 3.6.12)$, a
commutative subgroup $\G'$ of $\G$, and a fan $\Sig$ in $\fg_\bQ$
which is strongly compatible with $\G'$ such that the period map
$\vf|_{U\cap V}$ lifts to a morphism $U\cap V \to \G'\bs D$ which
extends uniquely to a morphism $V'\to\G'\bs D_\Sig$ of log manifolds.
$$
\matrix
U & \supset & U\cap V & \sub & V' \\
&&&&\\
{{}^\vf}\downarrow\;\;\; & & \downarrow &  & \downarrow\\
&&&&\\
\G \bs D & \leftarrow  & \G' \bs D &  \sub  & \G' \bs D_\Sig.
\endmatrix
$$
Furthermore, we have$:$
\medskip

{\rm(ii-1)} Assume $S\smallsetminus U$ is a smooth divisor. 
Then we can take $V=V'=S$ and $\G'=\G$. 
That is, we have a commutative diagram 
$$
\matrix U&\sub & S \\ 
&&&&\\
{{}^\vf}\downarrow\;\;\;&&\downarrow\\
&&&&\\
\G\bs D&\sub&\G\bs D_{\Sig}.\endmatrix
$$

{\rm(ii-2)} Assume $\G$ is commutative.
Then we can take $\G'=\G$.

\medskip

{\rm (ii-3)} 
Assume that $\G$ is commutative and that the following condition {\rm(1)} is satisfied.
\smallskip

{\rm(1)} There is a finite family $(S_j)_{1\leq j\leq n}$ of connected
locally closed analytic subspaces of $S$ such that $S=\tCu_{j=1}^n S_j$
as a set and such that, for each $j$, the inverse image of the sheaf
$M_S/\cO^\x_S$ on $S_j$ is locally constant.
\smallskip

Then we can take $\G'=\G$ and $V=S$. 
\endproclaim

Note that in (ii), we can take a fan $\Sig$ (we do not need a weak fan). 
\medskip

{\it Proof.} 
(i) (resp.\  (ii))  is the mixed Hodge version of Theorem 
0.5.29 (resp.\  Theorem 4.3.1)  of \cite{KU09}. 
By using the results in \S2 and \S4, the proof goes exactly in the same way as in the pure case treated in \cite{KU09}. 
\qed

\medskip

We deduce the following Theorem 7.4.2 from Theorem 7.4.1 (i).

\proclaim{Theorem 7.4.2}
Let the assumption be as in $7.4.1$. 
Let $\tilde H_\bR$ be the unique extension of $H_\bR$ as a local system on $S_{\val}^{\log}$. 
Then the canonical continuous splitting of the weight filtration $W$ of $H_\bR$ $(4.1.1)$ extends unquely to a continuous splitting of the weight filtration of $\tilde H_\bR$. 

\endproclaim

Theorem 7.4.2 was proved also by Brosnan and Pearlstein (\cite{BP.p} Theorem 2.21)  by another method. 

\medskip

{\it Proof of Theorem 7.4.2.} As in \cite{KNU.p}, the map $\spl_W: D\to \spl(W)$ 
in 4.1.1, where $\spl(W)$ denotes the set of all splittings of $W$,  extends to a continuous map
$D_{\SL(2)}^I\to \spl(W)$. 
Hence the period map $S_{\val}^{\log}\to \G\bs D_{\SL(2)}^I$ in Theorem 7.4.1 (i) induces a continuous map
$u:S_{\val}^{\log}\to \G\bs \spl(W)$. Since the map $\spl(W) \to \G\bs \spl(W)$ is a local homeomorphism, $u$ lifts locally on $S_{\val}^{\log}$ to a 
continuous map to $\spl(W)$. This shows that locally on $S_{\val}^{\log}$, the canonical continuous splitting of $W$ of $H_\bR$ given on $U$ extends  to a continuous splitting of $W$ of $\tilde H_\bR$. 
 Since the local extensions are unique, this proves Theorem 7.4.2.
 \qed

\enddocument